\documentclass{amsart}%

\usepackage{amsfonts}
\usepackage{amsmath}
\usepackage{amssymb}
\usepackage{amsthm}
\usepackage{booktabs}
\usepackage{mathtools, comment}
\usepackage[dvipsnames]{xcolor}
\usepackage[colorlinks, linkcolor=BrickRed,citecolor=Violet]{hyperref}
\usepackage[noabbrev,nameinlink]{cleveref}
\usepackage[shortlabels]{enumitem}

\usepackage{pdfpages}

\usepackage[margin= 0.9 in]{geometry}
\usepackage{varwidth}
\usepackage{subcaption}
\usepackage{multirow}
\usepackage{rotating}
\usepackage[boxsize=1em]{ytableau}

\usepackage{bbm}
\definecolor{viridian}{HTML}{0AAEAB}
\definecolor{darkcyan}{HTML}{099A98}
\definecolor{sapphire}{HTML}{1B4DC2}
\definecolor{plum}{HTML}{9948AB}
\definecolor{rose}{HTML}{C10091}
\definecolor{pumpkin}{HTML}{E47604}
\definecolor{sunflower}{HTML}{F6AE2D}
\definecolor{maize}{HTML}{FDE34F}
\definecolor{saffron}{HTML}{DAA520}
\definecolor{auburn}{HTML}{A73937}
\definecolor{darkred}{rgb}{0.7,0,0}

\newcommand\viridian[1]{{\color{viridian}#1}}
\newcommand\sapphire[1]{{\color{sapphire}#1}}
\newcommand\plum[1]{{\color{plum}#1}}

\newcommand\pumpkin[1]{{\color{pumpkin}#1}}

\newcounter{r}
\newcounter{s}
\newcommand\Tableau[1]{
        \setcounter{r}{0}
        \setcounter{s}{0}
        \foreach \x [count = \c from 1] in {#1} {
		\foreach \y [count = \d from 1] in \x{
			\node at (\d-.5,\c-.5) {\small$\y$};
			\draw (\d,\c) to (\d,\c-1);
			{\ifnum\d=1
				\draw (0,\c) to (0,\c-1);
				\fi}
			\setcounter{r}{\d}
		}
		{\ifnum\c=1
			\draw (0,0)--(\value{r},0);
			\fi}
		\draw(0,\c) to (\value{r},\c);
		\setcounter{s}{\c}}}
\newcommand\TAB[2][.3]{\TIKZ[scale=#1]{\Tableau{#2}}}

\usepackage{tikz}
\usetikzlibrary{calc, shapes,arrows,positioning,cd, patterns}
\tikzset{>=stealth',
  head/.style = {fill = white, text=black},
  plaque/.style = {draw, rectangle, minimum size = 10mm},
  pil/.style={->,thick},
  junct/.style = {draw,circle,inner sep=0.5pt,outer sep=0pt, fill=black}
  }
\tikzstyle{opLabel}=[fill=white, inner sep=2pt, rounded corners, pos=.3]
\tikzstyle{edgeLabel}=[fill=white, inner sep=2pt, rounded corners]
\tikzstyle{v}=[draw, fill =black, circle, inner sep=0pt, minimum size=2pt] 
\tikzstyle{xMapsto}=[line width=.6pt, |-{Classical TikZ Rightarrow[length=1mm]}]

\tikzstyle{xRightArrow}=[line width=.6pt, -{Classical TikZ Rightarrow[length=1mm]}]
\tikzstyle{intEdge}=[Bracket-Bracket, shorten <= -.5pt, shorten >= -.5pt]

\tikzstyle{TabNode}=[inner sep=2pt, fill=white, rounded corners]
\tikzstyle{LabEdge}=[inner sep=1.75pt, fill=white, rounded corners]

\newcommand\TIKZ[2][]{\begin{tikzpicture}[baseline={([yshift=-.8ex]current bounding box.center)}, #1]#2\end{tikzpicture}}

\def\circasterisk{\TIKZ{\node[circle, draw, inner sep=1pt] at (1,.5) {$*$};}}

\definecolor{ssyt1}{HTML}{BEE8D8}
\definecolor{ssyt2}{HTML}{FFF2AA}
\definecolor{ssyt3}{HTML}{C7BBFF}
\definecolor{ssyt4}{HTML}{FFBFF7}
\definecolor{bssyt1}{HTML}{1EA974}
\definecolor{bssyt2}{HTML}{D5A618}
\definecolor{bssyt3}{HTML}{4A2DCB}
\definecolor{bssyt4}{HTML}{B81FA4}
\def\PartDirection{1}
\newcommand\SSYTcolors[1]{       
        \foreach \x [count = \c from 1] in {#1} {
		\foreach \y [count = \d from 1] in \x{
			\draw[ssyt\y, fill = ssyt\y, rounded corners=0] (\d,\c) rectangle (\d-1,\c-1);}}}
\newcommand\cTAB[3][.35]{\TIKZ[scale=#1, yscale=\PartDirection]{\SSYTcolors{#3}\Tableau{#2}}}	

\newcommand\Bracketed[3]{
\draw [rounded corners, |-|, black!50](#1, -.6) to (#1,-1.1-#3*.3) to  (#2, -1.1-#3*.4) to (#2,-.6);
}

\newcommand\Path[1]{
        \setcounter{r}{0}
        \setcounter{s}{0}
        \foreach \x [count = \c from 1] in {#1}{\addtocounter{r}{2*\x - 3}
        		\setcounter{s}{\c}};
	\pgfmathsetmacro\lambdaone{-.5*\value{r}+.5*\value{s}}
	\pgfmathsetmacro\lambdatwo{.5*\value{r}+.5*\value{s}}
	\begin{scope}[black!20, thin]
	\draw (0,0) to (\lambdaone, -\lambdaone);
	\foreach \d in {\lambdaone, ..., \lambdatwo}{\draw (\d, -\d) to +(\lambdatwo,\lambdatwo);}
	\foreach \d [evaluate = \d as \h using \d-\lambdaone, count = \c from 1] in {1, ..., \lambdatwo}
		{\draw (\d, -\d) to (2*\d,0) to +(-\h,\h);}
	\end{scope}
	\draw[thick] (0,0) 
		\foreach \x [count = \c from 1] in {#1}{to coordinate (c\c) ++(1,2*\x - 3)};
		}

\newcommand\Part[1]{
        \setcounter{r}{1}
	 \foreach \x in {#1}{
 	{\ifnum\value{r}=1
		\draw (0,\value{r}-1)--(\x,\value{r}-1); 
		\fi}
	\draw (0,\value{r}) to (\x,\value{r});
   	\foreach \y in {0, ..., \x} {\draw (\y,\value{r})--(\y,\value{r}-1);}
	\addtocounter{r}{1}
 }}
 
\colorlet{leftcolor}{bssyt1}
\colorlet{rightcolor}{bssyt2}
\newcommand\leftcolor[1]{{\color{leftcolor}#1}}
\newcommand\rightcolor[1]{{\color{rightcolor}#1}}

\newcommand{\defncolor}{\color{sapphire!70!black}}
\newcommand{\defn}[1]{{\defncolor\emph{#1}}} 
\newcommand{\CS}{{\mathsf{CS}}}
\newcommand{\wt}{{\mathsf{wt}}}
\newcommand{\Des}{{\mathsf{Des}}}
\newcommand{\DE}{{\mathsf{DE}}}
\newcommand{\evac}{{\mathsf{evac}}}
\newcommand{\row}{{\mathsf{row}}}

\newcommand{\std}{{\mathsf{std}}}
\newcommand{\cycle}{\mathsf{cycle}}
\newcommand{\rev}{\mathsf{rev}}
\newcommand{\LL}{\mathcal{L}}
\newcommand{\len}{\mathsf{len}}
\newcommand{\SYT}{\mathsf{SYT}}
\newcommand{\SSYT}{\mathsf{SSYT}}

\newcommand{\nl}{n_{\mathsf{left}}}
\newcommand{\nr}{n_{\mathsf{right}}}
\newcommand{\jdt}{\mathsf{jdt}}
\newcommand{\shape}{\mathsf{shape}}
\newcommand{\destd}{\mathsf{destd}}

\def\ZZ{\mathbb{Z}}

\newcommand{\xdownarrow}[1]{%
  {\left\downarrow\vbox to #1{}\right.\kern-\nulldelimiterspace}
}

\usepackage[colorinlistoftodos]{todonotes}

\newtheorem{theorem}{Theorem}[section]

\newtheorem{corollary}[theorem]{Corollary}
\newtheorem{lemma}[theorem]{Lemma}

\newtheorem{proposition}[theorem]{Proposition}

\theoremstyle{definition}
\newtheorem{definition}[theorem]{Definition}

\newtheorem{remark}[theorem]{Remark}
\newtheorem{example}[theorem]{Example}

\newtheorem{axiom}[theorem]{Axiom}

\numberwithin{equation}{section}

\title{Crystal skeletons: Combinatorics and axioms}

\author[Brauner]{Sarah Brauner}
\address[S.\ Brauner]{Department of Mathematics, University of Pennsylvania, Philadelphia, PA, USA}
\email{sarahbrauner@gmail.com}
\urladdr{\href{https://www.sarahbrauner.com/}{https://www.sarahbrauner.com/}}

\author[Corteel]{Sylvie Corteel}
\address[S.\ Corteel]{Department of Mathematics, University of California, Berkeley, CA, USA and CNRS, IMJ-PRG, Sorbonne Universit\'e, France}
\email{corteel@berkeley.edu}

\author[Daugherty]{Zajj Daugherty}
\address[Z.\ Daugherty]{Department of Mathematics and Statistics, Reed College, 3203 SE Woodstock Blvd, Portland, OR 97202-8199, USA}
\email{zdaugherty@reed.edu}
\urladdr{\url{https://people.reed.edu/~zdaugherty/}}

\author[Schilling]{Anne Schilling}
\address[A. Schilling]{Department of Mathematics, University of California, One Shields
Avenue, Davis, CA 95616-8633, U.S.A.}
\email{anne@math.ucdavis.edu}
\urladdr{\href{http://www.math.ucdavis.edu/~anne}{http://www.math.ucdavis.edu/~anne}}

\begin{document}

\begin{abstract}
  Crystal skeletons were introduced by Maas-Gari\'epy in 2023 by contracting quasi-crystal components
  in a crystal graph. On the representation theoretic level, crystal skeletons model the expansion of
  Schur functions into Gessel's quasisymmetric functions. Motivated by questions of Schur positivity,
  we provide a combinatorial description of crystal skeletons, and prove many
  new properties, including a conjecture by Maas-Gari\'epy that crystal skeletons generalize dual equivalence graphs.
  We then present a new axiomatic approach to crystal skeletons. We give three versions of the axioms based on 
  $GL_n$-branching, $S_n$-branching, and local axioms in analogy to the local Stembridge axioms for crystals
  based on novel commutation relations.
\end{abstract}

\keywords{Crystal graphs, Lusztig involution, branching rules, dual equivalence graphs, Stembridge axioms}
\maketitle
\tableofcontents
\section{Introduction}
\label{section.introduction}

Crystal graphs provide combinatorial tools to study the representation theory of Lie algebras
(see~\cite{BumpSchilling.2017} for details). For instance, crystals are well-behaved with respect
to taking tensor products and hence can be used to give combinatorial interpretations for 
Littlewood--Richardson coefficients. In type $A$, the character of an irreducible crystal
$B(\lambda)$ of highest weight $\lambda$ is the Schur function $s_\lambda$ (see \S\ref{section.crystal}).

It is an important problem in representation theory and algebraic combinatorics to deduce the Schur
function expansion of a symmetric function whose expansion in terms of Gessel's fundamental
quasisymmetric function~\cite{Gessel.1984}  $F_\alpha$ is known. For example, combinatorial expressions 
for the quasisymmetric expansion of LLT polynomials, modified Macdonald polynomials~\cite{HHL.2005},
characters of higher Lie modules (or Thrall's problem)~\cite{GR.1993} or the plethysm of two Schur
functions~\cite{LW.2012} exist, yet their Schur expansions are in general still elusive.
It is thus desirable to develop methods to deduce the Schur expansions from these quasisymmetric
expansions. Some recent algebraic approaches include~\cite{ELW.2010,GR.2018,Gessel.2019,OSSZ.2024}.
Whereas Schur functions are characters of irreducible crystals in type $A$, Gessel's
fundamental quasisymmetric functions are characters of \defn{quasi-crystals}~\cite{MG.2023, CMRR.2025,CMRR.2026}, 
which are certain subcomponents of a crystal (see \S\ref{section.QS}). We exploit this fact in this paper by
providing a representation theoretic approach to this problem.

In~\cite{MG.2023}, Maas-Gari\'epy introduced the \defn{crystal skeleton} $\mathsf{CS}(\lambda)$ by
contracting the quasi-crystals in $B(\lambda)$ to a vertex. Since there is a unique standard tableau
$T\in \mathsf{SYT}(\lambda)$ in each quasi-crystal in $B(\lambda)$, it is natural to label the vertices 
of the crystal skeleton by standard tableaux. The crystal skeleton construction is the crystal analogue
of Gessel's formula~\cite{Gessel.1984}
\begin{equation}
\label{equation.s in F}
    s_\lambda = \sum_{T \in \mathsf{SYT}(\lambda)} F_{\mathsf{Des}(T)},    
\end{equation}
where $\mathsf{Des}(T)$ is the descent composition of the standard tableau $T$ (see \S\ref{section.tableaux}).
As such, crystal skeletons have the potential to serve as a powerful tool in deriving Schur expansions from quasisymmetric expansions.
This is supported by Theorem \ref{thm:DEinsideCS}, affirmatively answering a conjecture of Maas-Gari\'epy~\cite{MG.2023}, stating that
crystal skeletons generalize the \defn{dual equivalence graphs}. Dual equivalence graphs were developed by Assaf~\cite{Assaf.2007,Assaf.2015} 
and Roberts~\cite{Roberts.2014, Roberts.2014a} specifically as a paradigm for Schur positivity.
See~\cite{AssafBilley.2012, Blasiak.2016,BlasiakFomin.2017} for further applications.

\begin{figure}[!h]
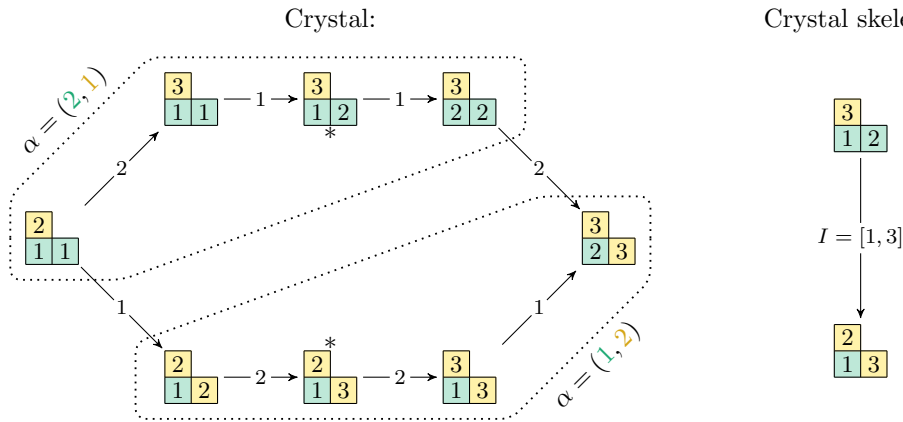

    \centering
\begin{tabular}{c@{\qquad\qquad}c}
 Crystal: & 
 Crystal skeleton:\\[5pt]
\TIKZ[scale=1.85]{
\draw[rounded corners, line width=.7pt,  dotted] (-.3, 0) to (-.3, -.3) to (.5,-.3) to (3.4,.75) to (3.4,1.3) to (.7,1.3) to 
 node[above, sloped]{$\alpha = (\leftcolor{2},\rightcolor{1})$} ++(-1,-1) to (-.3,0);
 \begin{scope}[shift={(4,0)}, scale = -1]
 \draw[rounded corners, line width=.7pt,  dotted] (-.3, 0) to (-.3, -.3) to (.5,-.3) to (3.4,.75) to (3.4,1.3) to (.7,1.3) to 
 	node[below, sloped]{$\alpha = (\leftcolor{1},\rightcolor{2})$} ++(-1,-1) to (-.3,0);
 \end{scope}
 \node[TabNode] (211)  at (0,0) {\cTAB{{1,1},{2}}{{1,1},{2}}};
 \node[TabNode] (212) at (1,-1) {\cTAB{{1,2},{2}}{{1,2},{2}}};
\node[TabNode] (311) at (1,1) {\cTAB{{1,1},{3}}{{1,1},{2}}};
\node[TabNode] (213) at (2,-1) {\cTAB{{1,3},{2}}{{1,2},{2}}};
\node[TabNode] (312) at (2,1) {\cTAB{{1,2},{3}}{{1,1},{2}}};
\node[TabNode] (322) at (3,1) {\cTAB{{2,2},{3}}{{1,1},{2}}};
\node[TabNode] (313) at (3,-1) {\cTAB{{1,3},{3}}{{1,2},{2}}};
 \node[TabNode] (323) at (4,0) {\cTAB{{2,3},{3}}{{1,2},{2}}};
\begin{scope}[every node/.style={LabEdge}]
\draw[->, shorten <=-2pt, shorten >=-2pt] (211) to node {\footnotesize$1$} (212);
\draw[->] (211) to node {\footnotesize$2$} (311);
\draw[->] (311) to node {\footnotesize$1$} (312);
\draw[->] (212) to node {\footnotesize$2$} (213);
\draw[->] (213) to node {\footnotesize$2$} (313);
\draw[->] (312) to node {\footnotesize$1$} (322);
\draw[->, shorten <=-2pt, shorten >=-2pt] (322) to node {\footnotesize$2$} (323);
\draw[->] (313) to node {\footnotesize$1$} (323);
\end{scope}
\node[below, shift={(0,-.25)}] at (312) {$*$};
\node[above, shift={(0,.25)}] at (213) {$*$};
 }
&
\TIKZ[scale=2]{
\node[TabNode] (312) at (0,.75) {$\cTAB{{1,2},{3}}{{1,1},{2}}$};
\node[TabNode] (213) at (0,-.75) {$\cTAB{{1,3},{2}}{{1,2},{2}}$};
\draw[->] (312) to node[LabEdge] {\footnotesize$I=[1,3]$} (213);
}
\end{tabular}
    \caption{Left: Crystal $B(2,1)$ of type $A_2$ with two quasi-crystal components indicated with dotted lines and standard tableaux indicated by $*$.
    The descent composition is denoted $\alpha$.
    Right: Corresponding crystal skeleton.}
    \label{figure.B21}
\end{figure}

A crystal, its quasi-crystal components and the corresponding crystal skeleton are shown in Figure~\ref{figure.B21}.

Our goal in this paper is to characterize the crystal skeleton both combinatorially 
and axiomatically in analogy to the local Stembridge axioms for crystals~\cite{Stembridge.2003}. 
Stembridge axioms have played a crucial role in crystal theory and have facilitated proofs of Schur positivity using crystals.
For example, in~\cite{MS.2016} the Schur expansion of Stanley symmetric functions was analyzed by defining a crystal structure
on the combinatorial objects underlying Stanley symmetric functions (decreasing factorizations of a permutation). The crystal 
structure was proved using Stembridge's axioms. Our new axioms for crystal skeletons will have similar applications for Schur positivity in cases where the 
quasisymmetric expansion is known.

We summarize these characterizations below. 
An extended abstract with summaries of our results is available~\cite{BCDS24}.

\subsection{Combinatorics of the crystal skeleton}\label{subsection.intro.combinatorics}
\begin{figure}[t]
\begin{center}
\scalebox{1}{\includegraphics{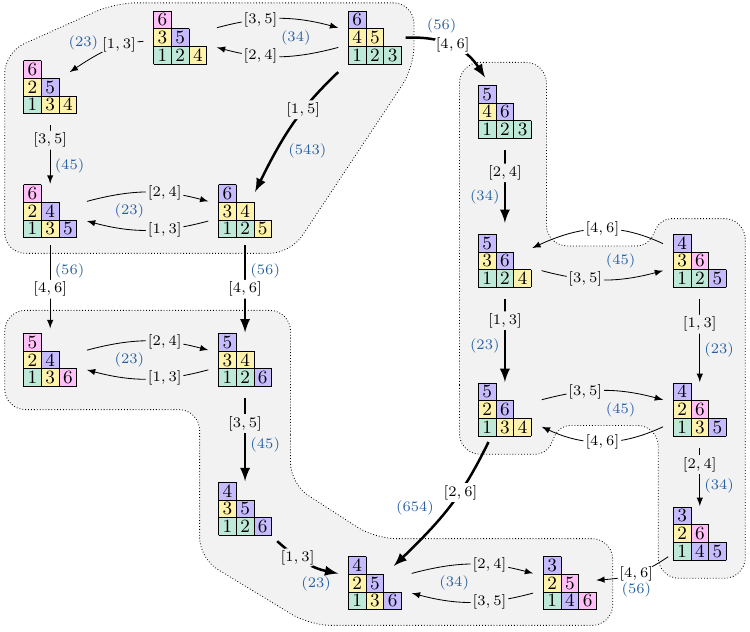}}
\end{center}
\caption{The crystal skeleton $\CS(3,2,1)$. The edge labels in terms of intervals and cycles are defined in Sections~\ref{section.Dyck pattern} 
and~\ref{section:cycles}. This example is further decorated by tableaux colored to indicate the descent compositions as in 
Example~\ref{example.descentcomposition}; thick arrows to indicate the crystal $B(3,2,1)_3$ in Theorem~\ref{theorem.B short}; 
and gray components to indicate the branching in Theorem~\ref{theorem.branching}.
\label{figure.CS321}}
\end{figure}

The original definition of the crystal skeleton by Maas-Gari\'epy comes from contracting edges of the crystal. Here, our work 
provides a self-contained combinatorial description of the crystal skeleton that does not reference the crystal. We highlight the combinatorial 
properties proved in~\S\ref{section.CS} and~\S\ref{section.CS properties} with the following example of the crystal skeleton $\CS(3,2,1)$ in 
Figure~\ref{figure.CS321}.

\subsubsection{Vertices and edges of the crystal skeleton}
The vertices of $\CS(\lambda)$ can be labeled in two ways:
\begin{enumerate}
    \item by the set $\SYT(\lambda)$, which is shown in Figure \ref{figure.CS321} for $\SYT(3,2,1)$; and  
    \item by compositions $\alpha$ of $n=|\lambda|$, where $\alpha = \Des(T)$ is the \defn{descent composition} of the tableau $T \in \SYT(\lambda)$. 
    This is shown in Figure \ref{figure.CS321} by the coloring of the tableaux, where each color represents a part of the composition, and the 
    composition order is given by ordering the numbers $1, \ldots, n$.
\end{enumerate}
We provide two natural ways to label the edges in $\CS(\lambda)$:
\begin{enumerate}
    \item We show in \S\ref{section.Dyck pattern} that the edges of $\CS(\lambda)$ can be labeled by certain odd length intervals $I$ called
    \defn{Dyck pattern intervals}. These intervals can be described in terms of Dyck paths; the interval $I$ corresponds to a Dyck path with $|I|+1$ steps.
    \item In \S\ref{section:cycles} we show that the edges of $\CS(\lambda)$ can also be labeled by certain decreasing \defn{cycles} 
    $\cycle(T_I)$, where $I$ is the Dyck pattern interval and $T_I$ is the restriction of $T$ to the letters in $I$. In particular, if there is an edge between 
    $T$ and $T'$ in $\CS(\lambda)$, then 
    \[ \cycle(T_I) \cdot T = T'.\]
\end{enumerate}
The above marks the first explicit descriptions of the edge labels of $\CS(\lambda)$, which are ambiguous in \cite{MG.2023}. 

We then determine the relationship between the edge and vertex labels in $\CS(\lambda)$ by describing how the descent compositions 
between adjacent vertices of $\CS(\lambda)$ differ. See Theorem \ref{theorem.descent composition} for a complete description. 

\begin{theorem}[Theorem \ref{theorem.descent composition} summary]\label{theorem.descentcomposition.introduction}
    Suppose there is an edge in $\CS(\lambda)$ labeled by 
    \[ 
	\big(T,\alpha = \Des(T)\big) \xrightarrow{~~I~~} \big(T',\beta = \Des(T')\big),
\]
where $\alpha$ has $\ell$ parts. Then $\beta$ has length $\ell-1,\ell$ or $\ell+1$, and there are combinatorial conditions on $T$ and $I$ that 
determine $\beta$.
\end{theorem}

\subsubsection{Subgraphs of the crystal skeleton}
We next show that $\CS(\lambda)$ contains several interesting subgraphs. 

First, we answer affirmatively a conjecture of Maas-Gari\'epy~\cite{MG.2023} by proving that crystal skeletons generalize the \defn{dual equivalence graphs} 
developed by Assaf~\cite{Assaf.2007,Assaf.2015} and Roberts~\cite{Roberts.2014, Roberts.2014a}. Write $\DE(\lambda)$ for the dual 
equivalence graph corresponding to the partition $\lambda$.

\begin{theorem}[Theorem \ref{thm:DEinsideCS}]
The dual equivalence graph $\mathsf{DE}(\lambda)$ is a subgraph of the crystal skeleton $\mathsf{CS}(\lambda)$. 
\end{theorem}

In particular, $\DE(\lambda)$ is obtained from $\CS(\lambda)$ by including edges with intervals $I$ of length $3$, or equivalently, cycles of 
length two.  For example, the graph $\DE(3,2,1)$ can be seen in Figure~\ref{figure.CS321} by considering only the edges with intervals of length 3.

Second, we show that $\CS(\lambda)$ has surprising \defn{self-similarity} properties. Let $\lambda \vdash n$ be a partition of $n$. 
Given $T \in \SYT(\lambda)$ and an interval $[a,b] \subseteq [1,n]$, the skew tableau $T_{[a,b]}$ obtained by restricting $T$ to $[a,b]$
can be straightened to a straight shaped tableau via \defn{jeu de taquin}, which we write as $\jdt(T_{[a,b]})$. 

\begin{theorem}[Theorem \ref{theorem.self similar}]
Let $T\in \mathsf{SYT}(\lambda)$ and $[a,b]\subseteq [1,n]$ be an interval. Set $\mu=\mathsf{shape}(\jdt(T_{[a,b]}))$.
Then $\mathsf{CS}(\mu)$ is a subgraph of $\CS(\lambda)$ (up to relabelling of the edges). 
\end{theorem}

\noindent As a special case of Theorem \ref{theorem.self similar}, we obtain elegant branching properties mirroring those of the symmetric group.

\begin{theorem}[Theorem \ref{theorem.branching}]\label{theorem.branching.introduction}
For $\lambda \vdash n$, denote by $\CS(\lambda)_{[1,n-1]}$ the restriction of the crystal skeleton graph $\CS(\lambda)$ 
by removing all edges labeled by $I$ such that $n\in I$. Then we have a graph isomorphism
\[
	\CS(\lambda)_{[1,n-1]} \cong \bigoplus_{\lambda^-} \CS(\lambda^-),
\]
where the sum is over all $\lambda^-$ such that $\lambda/\lambda^-$ is a skew shape with a single box.
\end{theorem}

\noindent
Theorem \ref{theorem.branching} thus implies that each $\CS(\lambda^-)$ is a subgraph of $\CS(\lambda)$. This is shown in Figure~\ref{figure.CS321} 
by the three components in gray, obtained by removing the edges labeled by intervals containing $6$. These components are isomorphic to $\CS(2,2,1)$, 
$\CS(3,1,1)$, and $\CS(3,2)$, moving from the bottom to the top, counter-clockwise. In all cases, the isomorphism sends the vertex $T \in \CS(\lambda)$ 
to $T_{[1,n-1]}$.

The third surprising fact is that the crystal itself is a subgraph of the crystal skeleton, despite the fact that the crystal skeleton is obtained from a crystal by
contracting quasi-crystal components. The length (or number of parts) of a partition $\lambda$ is denoted by $\ell(\lambda)$.

\begin{theorem}[Theorem~\ref{theorem.B short}]\label{theorem.B.short.introduction}
The induced subgraph of the crystal skeleton $\CS(\lambda)$ of vertices with descent composition of length $\ell(\lambda)$ is isomorphic 
to the crystal graph $B(\lambda)$ of type $A_{\ell(\lambda)-1}$.
\end{theorem}

\noindent
The crystal $B(3,2,1)$ of type $A_2$ (or equivalently $B(2,1)$ of type $A_2$) is a subgraph of $\CS(3,2,1)$ as indicated by the 
bold edges in Figure~\ref{figure.CS321}.

\subsubsection{Symmetry in the crystal skeleton}
Finally, we show that $\CS(\lambda)$ satisfies similar symmetry properties to crystals.  It is well-known that any crystal $B(\lambda)$ 
is invariant under the \defn{Lusztig involution}, which is a map 
\[ 
	\LL \colon B(\lambda) \longrightarrow B(\lambda)
\]
sending the highest weight of the crystal to the lowest and reversing (and relabeling) edges.

We show that this map descends to the crystal skeleton, defining a Lusztig involution $\LL_n$ on $\CS(\lambda)$ by mapping the vertex 
$T\in \SYT(\lambda)$ to $\LL(T)$ and an edge 
\[ 
	T \stackrel{~I~}{\longrightarrow} T' \quad \quad \textrm{ to } \quad \quad \LL(T') \stackrel{~I^\LL~}{\longrightarrow} \LL(T), 
\] 
where $I^\LL = [n+1-b,n+1-a]$ if $I=[a,b]$.

\begin{corollary}[Corollary \ref{corollary.lusztig}]\label{corollary.lusztig.intro}
Let $\lambda \vdash n$. The crystal skeleton $\mathsf{CS}(\lambda)$ is invariant under the Lusztig involution $\mathcal{L}_n$, that is,
\[
	\CS(\lambda) \cong \mathcal{L}_n(\CS(\lambda)).
\]
\end{corollary}

\subsection{Axioms for crystal skeletons}\label{section.intro.axioms}
The next major focus of our work is to characterize crystal skeleton graphs with axioms, analogously to the Stembridge axioms for 
crystal graphs~\cite{Stembridge.2003}.
We give three different axiomatic characterizations of the crystal skeleton $\CS(\lambda)$, described below. 
\begin{theorem}[Summary of \S\ref{section.axioms}]
There are three axiomatic descriptions that uniquely characterize crystal skeleton graphs:
\begin{enumerate}
    \item $GL_n$-Axioms, found in \S\ref{section.GLn axioms};
    \item $S_n$-Axioms, found in \S\ref{section.axiom Sn}; and 
    \item Local Axioms, found in \S\ref{section.axiom local}.
\end{enumerate}
\end{theorem}

The common thread of each of these perspectives is to prove that the combinatorial descriptions of $\CS(\lambda)$ established above 
in~\S\ref{subsection.intro.combinatorics} completely determine the crystal skeleton. All three axiom sets take as input a graph with vertices 
labeled by compositions of $n = |\lambda|$ and edges labeled by odd 
length intervals. In all three cases, there are several local axioms governing the behavior of the edge and vertex labels, including satisfying 
Theorem~\ref{theorem.descentcomposition.introduction}; see also
\ref{axiom.intervals},
\ref{axiom.outgoing},
and \ref{axiom.labels}.
We highlight below the differences between the axiom sets. 

\smallskip
\paragraph{\emph{$GL_n$-Axioms}}
Let $G$ be a graph with vertex set labeled by compositions of $n$ and edges labeled by odd-length intervals of the set $[n]$. Our first axiomatic 
description of $\CS(\lambda)$ in~\S\ref{section.GLn axioms} requires that the input graph $G$ satisfies a local \emph{fan} condition 
(see \ref{axiom.fans}), as well as:
\begin{itemize}
    \item  \emph{Strong Lusztig involution} (\ref{axiom.lusztig}) as in Corollary~\ref{corollary.lusztig.intro} on both $G$ and the restriction of $G$ to 
    $G_{[1,n-1]}$,  where $G_{[1,n-1]}$ is the subgraph obtained by removing edges of $G$ whose interval label contains $n$:
    \[ \mathcal{L}_{n-1}(G_{[1,n-1]}) \cong G_{[1,n-1]}.\]
 \item \emph{Top Subcrystal} (\ref{axiom.crystal}): $G$ must contain a subgraph isomorphic to a crystal $B(\lambda)$ as in 
 Theorem~\ref{theorem.B.short.introduction}, obtained from the subgraph of $G$ formed from vertices whose composition labels are of minimal length.
\end{itemize}
We call these axioms $GL_n$-axioms because restricting compositions by length mirrors $GL_n$ branching rules in representation theory. 

\smallskip
\paragraph{\emph{$S_n$-Axioms}}
We give our second axiomatic description of $\CS(\lambda)$ in~\S\ref{section.axiom Sn}, where we require that the input graph $G$ satisfies the 
same local fan condition as the $GL_n$-axioms (\ref{axiom.Sn fans}), as well as:
\begin{itemize}
    \item \emph{Lusztig involution} (\ref{axiom.Sn lusztig}):  $\mathcal{L}_n(G) \cong G$ as in Corollary \ref{corollary.lusztig.intro};
    \item \emph{Branching} (\ref{axiom.Sn branching}): $G$ must satisfy a branching recurrence analogous to 
    Theorem~\ref{theorem.branching.introduction};
    \item \emph{Connectivity} (\ref{axiom.Sn connectivity}): $G$ must satisfy certain connectivity conditions.
\end{itemize}
These axioms are called $S_n$-Axioms because the graphs must satisfy branching that mirrors the $S_n$ representation theoretic branching rules. 

\smallskip
\paragraph{\emph{Local Axioms}}
Finally, our third axiomatic description of the crystal skeleton in~\S\ref{section.axiom local} most closely mirrors the Stembridge axioms for crystal 
graphs. This set of axioms requires that the input graph $G$ satisfy \emph{local commutation relations} (\ref{axiom.local commutations}) 
(see also Theorem~\ref{theorem.commutation}) similar to 
the Stembridge commutation rules for crystals and certain conditions on string lengths (\ref{axiom.local string}).

\subsection{Outline}

In Section~\ref{section.background}, we provide the background on crystals, quasi-crystals and dual equivalence graphs.
In Section~\ref{section.CS}, we give the definition and our new combinatorial description of the crystal skeleton.
In Section~\ref{section.CS properties}, we prove various properties of the crystal skeleton.
The axioms for crystal skeletons are discussed in Section~\ref{section.axioms}.
We conclude in Section~\ref{section.two row} by analyzing the crystal skeleton for two-row partitions in more detail to
illustrate our results.

\subsection*{Acknowledgements}

The authors would like to thank the organizers of the conference at the Banff International Research Station 
entitled ``Community in Algebraic and Enumerative Combinatorics'' in January 2024, where this work began
as a research project proposed by the last author. We would like to thank Zach Hamaker, Florence Maas-Gari\'epy,
Ant\'onio Malheiro, Rosa Orellana, Joseph Pappe, Franco Saliola, Mary Claire Simone, and Mike Zabrocki
for fruitful discussions. 

SB was partially supported by the NSF MSPRF DMS-2303060 and NSF Grant DMS-2554084.
SC was partially supported by NSF grant DMS-2054482 and ANR grant  ANR-19-CE48-0011.
AS was partially supported by  NSF grant DMS--2053350 and Simons Foundation grant MPS-TSM-00007191.

\section{Background: Crystals, quasi-crystals, and dual equivalence graphs}
\label{section.background}

In this section, we review some basics about tableaux, crystals, quasi-crystals and dual equivalence graphs, which will be used
in this paper.

\subsection{Tableaux}
\label{section.tableaux}

Let $\mathsf{SSYT}(\lambda)_n$ be the set of \defn{semistandard Young tableaux} of shape $\lambda$ over the alphabet $[n]:=\{1,2,\ldots,n\}$ and 
$\SYT(\lambda)$ be the set of \defn{standard Young tableaux} of shape $\lambda$. 

We use French notation for partitions and tableaux, where the sizes of the rows weakly decrease from bottom to top. Work of Robinson--Schensted--Knuth 
defines a bijection $\mathsf{RSK}$ between words in the alphabet $[n]$ and pairs $(P,Q)$ of a semistandard tableau $P$ over $[n]$ and a standard tableau 
$Q$ of the same shape. If $\mathsf{RSK}(w) = (P,Q)$, then $P$ is known as the \defn{insertion tableau} and $Q$ as the \defn{recording tableau}.
To indicate the dependence on the word $w$, we also write $P(w)$ and $Q(w)$ for the insertion and recording tableau of $w$. See for 
example~\cite[Chapter 7.11]{Stanley.EC2} and~\cite[Chapter 3]{Sagan.2001} for more details about the Robinson--Schensted--Knuth (RSK) correspondence.

There are two elementary Knuth equivalence relations for letters $a,b,c$
\begin{equation}
\label{equation.Knuth}
\begin{split}
	acb & \equiv_1 cab \qquad \text{if $a\leqslant b<c$,}\\
	bac & \equiv_2 bca \qquad \text{if $a<b\leqslant c$.}
\end{split}
\end{equation}
Two words $w$ and $v$ are \defn{Knuth equivalent}, denoted $w\equiv_K v$, if $w$ can be transformed to $v$ by a sequence of type 1 and 2 Knuth relations
on adjacent letters. It is well-known that $P(w)=P(v)$ if and only if $w \equiv_K v$, that is, $w$ and $v$ are Knuth equivalent.

A semistandard tableau $b \in \SSYT(\lambda)_n$ gives rise to several combinatorial objects: 
\begin{itemize}
    \item  The (row) \defn{reading word} $\mathsf{row}(b)$ is the word obtained from $b$ by reading rows left to right, top to bottom. 
    \item The \defn{weight} $\mathsf{wt}(b)$ is the tuple $(\alpha_1,\alpha_2,\ldots,\alpha_n)$ with $\alpha_j$ the number of letters $j$ in $b$.
    \item The \defn{standardization} $\mathsf{std }(b)$
is obtained from $b$ and $\wt(b)$ by replacing the letters $i$ in $b$ from left to right by 
\[ \alpha_1+\alpha_2+\cdots+\alpha_{i-1}+1,\ldots,
\alpha_1+\alpha_2+\cdots+\alpha_i \textrm{ for all }1\leqslant i \leqslant n.\]  
The standardization of a word $w$ can be defined similarly. By construction, $\std(b) \in \SYT(\lambda)$ and $\std(w)$ is a permutation.
\end{itemize}

\begin{example}  
\[
         \text{Suppose} \quad
	b=
	\raisebox{0.3cm}{\begin{ytableau}
	4\\
	2&4\\
	1&3&3
	\end{ytableau}}.  \qquad \text{Then} \qquad \std(b) =  \raisebox{0.3cm}{\begin{ytableau}
	5\\
	2&6\\
	1&3&4
	\end{ytableau}}
\]
and $\row(b) = 424133$, $\wt(b) = (1,1,2,2)$. Note that $\row(\std(b)) = 526134 = \std(\row(b))$.
\end{example}

For a standard tableau $T \in \mathsf{SYT}(\lambda)$, the letter $i$ is a \defn{descent} if the letter $i+1$ is in a higher row of the tableau 
(in French notation). Denote the descents of $T$ by $d_1<d_2<\cdots<d_k$. The \defn{descent composition} is defined as
\begin{equation}
\label{equation.descent composition}
	\mathsf{Des}(T) = ( d_1, d_2 - d_1, \ldots, d_k - d_{k-1}, n-d_k),
\end{equation}
where $n = |\lambda|$.

\subsection{Crystals}
\label{section.crystal}

We briefly review the crystal of type $A_{n-1}$ on tableaux. More details can be found in~\cite[Chapters 3, 8]{BumpSchilling.2017}. 
The \defn{crystal} $B(\lambda)_n$ is the set $\mathsf{SSYT}(\lambda)_n$ together with the maps
\begin{equation}
\begin{split}
	\mathsf{wt} &\colon B(\lambda)_n \to \mathbb{Z}^n_{\geqslant 0},\\
	e_i, f_i &\colon B(\lambda)_n \to  B(\lambda)_n \cup \{ \emptyset\} \qquad \text{for $i\in \{1,2,\ldots,n-1\}$.}
\end{split}
\end{equation}

The \defn{crystal raising} and \defn{crystal lowering} operators $e_i$ and $f_i$ are defined as follows.
The operators $e_i$ and $f_i$ act on the subword of $w=\mathsf{row}(b)$ containing only the letters $i$ and $i+1$, denoted by $w^{(i)}$.
Associate to each letter $i+1$ an open bracket $($ and to each letter $i$ a closed bracket $)$. Then pair letters with matched open and closed brackets, called
bracketed letters. This procedure is also referred to as the signature rule (see for example~\cite[Chapter 2.4]{BumpSchilling.2017})
The subword of unmatched or unbracketed letters is of the form $i^r (i+1)^s$. On this subword
\begin{equation*}
	e_i(i^r (i+1)^s) = \begin{cases} i^{r+1} (i+1)^{s-1} & \text{if $s>0$,}\\
	\emptyset & \text{else,}
	\end{cases}
	\quad
	f_i(i^r (i+1)^s) = \begin{cases} i^{r-1} (i+1)^{s+1} & \text{if $r>0$,}\\
	\emptyset & \text{else.}
	\end{cases}
\end{equation*}
All other letters in $w$ remain unchanged.

Crystal operators are well-behaved with respect to Knuth equivalence (see for example~\cite[Theorem 8.4]{BumpSchilling.2017}) as stated
in the next proposition.
\begin{proposition}
\label{proposition.fi and Knuth}
Let $w\equiv_K v$ be two Knuth equivalent words. Then $f_i w \equiv_K f_i v$ as long as $f_iw \neq \emptyset$.
\end{proposition}
Proposition~\ref{proposition.fi and Knuth} ensures that crystal operators are well-defined on tableaux as well. 

\begin{example}
\label{example.f3}
Consider how $f_3$ acts on $b \in \SSYT(10,9,3,1)_4$:
\[
b = \TIKZ[scale=.4]{
\filldraw[ssyt3] (6,0) rectangle (9,1) 
		(3,1) rectangle (5,2)
		(0,2) rectangle (2,3);
\filldraw[ssyt4] (9,0) rectangle (10,1) 
		(5,1) rectangle (9,2)
		(2,2) rectangle (3,3)
		(0,3) rectangle (1,4);
\Tableau{{1,1,1,1,2,2,3,3,3,4},{2,2,2,3,3,4,4,4,4}, {3,3,4},{4}}}
\qquad \text{with} \qquad
f_3b=
\TIKZ[scale=.4]{
\filldraw[ssyt3] (6,0) rectangle (9,1) 
		(3,1) rectangle (4,2)
		(0,2) rectangle (2,3);
\filldraw[ssyt4] (9,0) rectangle (10,1) 
		(4,1) rectangle (9,2)
		(2,2) rectangle (3,3)
		(0,3) rectangle (1,4);
\Tableau{{1,1,1,1,2,2,3,3,3,4},{2,2,2,3,4,4,4,4,4}, {3,3,4},{4}}}.
\]
\end{example}

A crystal is often encoded in a \defn{crystal graph} with vertices in $\SSYT(\lambda)$ and an edge $b\stackrel{i}{\longrightarrow} b'$ if 
$b'=f_i(b)$. Since the crystal raising and lowering operators are partial inverses, that is, $b'=f_i(b)$ if and only if $e_i(b')=b$, there is no loss of
information just recording the action of $f_i$ as directed edges. The advantage of viewing the crystal as an edge-labeled, vertex-weighted graph is that
one can study subgraphs and graph symmetries.

The \defn{string lengths} for $b\in B(\lambda)_n$ and $1\leqslant i<n$, which play an important role for Stembridge axioms~\cite{Stembridge.2003}, 
are defined as follows
\begin{equation}
\label{equation.string lengths}
	\varphi_i(b) = \max\{k \in \mathbb{Z}_{\geqslant 0} \mid f_i^k(b) \neq \emptyset\} \quad \text{and} \quad
	\varepsilon_i(b) = \max\{k \in \mathbb{Z}_{\geqslant 0} \mid e_i^k(b) \neq \emptyset\}.
\end{equation}

The crystal $B(\lambda)_n$ of type $A_{n-1}$ enjoys a symmetry under the Sch\"utzenberger or Lusztig involution 
(see for example \cite[p. 79]{BumpSchilling.2017}). To state it, we first need the definition of highest and lowest weight vectors. A vector $u$ in a crystal
is called \defn{highest weight} if $e_i u = \emptyset$ for all $1\leqslant i<n$. Similarly, a vector $v$ in a crystal is called \defn{lowest weight}
if $f_i v = \emptyset$ for all $1\leqslant i<n$. The crystal $B(\lambda)_n$ has a unique highest weight vector, denoted $u_\lambda$, of weight 
$\lambda=(\lambda_1,\ldots,\lambda_n)$. As a tableau, it is given by the superstandard tableau with letters $i$ in row $i$. Similarly, the crystal 
$B(\lambda)_n$ has a unique lowest weight vector, denoted $v_\lambda$, of weight $(\lambda_n,\ldots,\lambda_1)$.

\begin{definition}
\label{definition.lusztig}
Let $B(\lambda)_n$ be the crystal of type $A_{n-1}$ with highest weight element $u_\lambda$ of highest weight $\lambda$.
The \defn{Lusztig involution}
\[
	\LL \colon B(\lambda)_n \to B(\lambda)_n
\]
is defined as follows. Any $b\in B(\lambda)_n$ can be obtained from $u_\lambda$ by applying a sequence of lowering operators, that is, 
$b=f_{i_1}\cdots f_{i_k} u_\lambda$ for some sequence $1\leqslant i_j<n$. We define
\[
	\LL(b) = e_{n-i_1} \cdots e_{n-i_k} v_\lambda.
\]
In particular, $\LL(u_\lambda)=v_\lambda$ and $\LL(v_\lambda)=u_\lambda$.
\end{definition}

\begin{remark}
\label{remark.evacuation}
For crystals of tableaux in type $A$, the Lusztig involution is the same as \defn{evacuation} (also known as the \defn{Sch\"utzenberger involution}) on tableaux. 
Evacuation can be defined using jeu de taquin (see for example~\cite[p. 425]{Stanley.EC2}). Alternatively, let $T \in \mathsf{SSYT}(\lambda)$ and 
$w=w_1\ldots w_\ell = \mathsf{row}(T)$. Then $\mathsf{evac}(T)$ is the RSK insertion tableau of the word $w^\# = (n+1-w_\ell) \ldots (n+1-w_1)$.
Equivalently, Butler~\cite{Butler.1994} and Lenart~\cite[\S2.4 and Fig. 1]{Lenart.Lusztig-involution} used the following description of evacuation on $T$:
rotate $T$ by $180^\circ$, take the complement of each entry by replacing $i$ by $n+1-i$ and finally apply jeu de taquin. 
\end{remark}

The \defn{character} of a crystal is the sum over all elements in the crystal weighted by the weight function. 
For the crystal $B(\lambda)_n$, the character is the \defn{Schur polynomial} in $n$ variables
\begin{equation}
	s_\lambda(x_1,\ldots,x_n) = \sum_{b \in \mathsf{SSYT}(\lambda)_n} x^{\wt(b)}.
\end{equation}

\subsection{Quasi-crystals}
\label{section.QS}

We now define \defn{quasi-crystals} inside of the crystal $B(\lambda)_n$, which were first defined by Maas-Gari\'epy~\cite{MG.2023}. 
Her definition arises from studying ~\eqref{equation.s in F} from the perspective of crystal theory. 

\begin{definition}
For a given crystal $B(\lambda)_n$, the \defn{quasi-crystal} $Q_T$ associated to a standard tableau $T \in \mathsf{SYT}(\lambda)$  is the subgraph  
of $B(\lambda)_n$ by restricting to the elements $b\in B(\lambda)_n$ such that $\std(b) = T$. We refer to the distinct quasi-crystals in
$B(\lambda)_n$ as its quasi-crystal components.
\end{definition}

\begin{example}
Suppose $T \in \SYT(3,2,1)$ has $\row(T)= 645123$. Then the tableau with reading word $w = 322111$ is in the same quasi-crystal as $T$ in 
$B(3,2,1)_6$, since $\std(w) = \row(T)$.
\end{example}

\begin{example}\label{example.crystalandquasicrystals}
Figure~\ref{figure.B21} shows the crystal $B(2,1)_3$ and its two quasi-crystal components.
\end{example}

It was shown in~\cite[Thm.\ 1]{MG.2023} that the components $Q_T$ are connected in $B(\lambda)_n$\footnote{Note that in~\cite{MG.2023} the 
quasi-crystal components are defined in a slightly different way by fixing the descent composition. This is not quite accurate since the descent 
composition does not uniquely specify a quasi-crystal component.}. 
A local characterization of quasi-crystals was given in~\cite{CMRR.2025}.
Each quasi-crystal $Q_T$ contains exactly one standard tableau, namely $T$. 
In~\cite[Section 2.5.2]{CMRR.2026}, the crystal operators staying within a quasi-crystal component are characterized and denoted by
$\ddot{f}_i$ and $\ddot{e}_i$.

\begin{proposition}[\cite{CMRR.2026}, Section 2.5.2]
\label{proposition.quasi edges}
For $b\in Q_T \subseteq B(\lambda)_n$ with $f_i(b) \in B(\lambda)_n$, we have $f_i(b) \in Q_T$ 
if and only if no $i+1$ and $i$ are bracketed in $\mathsf{row}(b)$ (or equivalently, no $i+1$ appears to the left of an $i$ in $\mathsf{row}(b)$).
\end{proposition}

Proposition~\ref{proposition.quasi edges} also follows from the observation that the crystal operator $f_i$ under standardization acts 
as a cycle. To explain this, we start with a couple of definitions. This observation will be important later.

Let $b\in B(\lambda)_n$ and $1\leqslant i<n$ such that $f_i(b)\neq \emptyset$. Set $w = w_1 w_2 \ldots w_N = \mathsf{row}(b)$. Let 
$p$ be the position of the letter $i$ in $w$ on which $f_i$ acts, which is the rightmost unbracketed letter $i$. Furthermore, let 
$\nl$ (resp. $\nr$) be the number of bracketed  $(i+1, i)$ pairs to the left (resp. right) of $w_p$ in $w$. Note that $w_p$ is mapped to 
$\pi_p$ under standardization $\pi=\mathsf{std}(w)$. We define the \defn{cycle} of size $\nl+\nr+1$ as
\begin{equation}
\label{equation.sigma}
	\cycle(b,i) = (\pi_p+\nl+\nr, \pi_p+\nl+\nr-1,\ldots,\pi_p).
\end{equation}

\begin{remark}
\label{rmk:lengthofsigma}
The permutation $\cycle(b,i) \in S_N$ is a $(\nl + \nr + 1)$-cycle. In particular, if $w=\mathsf{row}(b)$ has no bracketed pair $(i+1, i)$
(or equivalently $\nl=\nr=0$), then $\cycle(b,i)$ is the identity. If $w$ has one bracketed pair $(i+1, i)$, then $\cycle(b,i)$ is a transposition. 
\end{remark}

\begin{example}
\label{example.sigma b}
Let
\[
	b=
	\raisebox{0.3cm}{\begin{ytableau}
	3\\
	2&3\\
	1&2&2
	\end{ytableau}}
	\qquad \text{with} \qquad
	f_2 b=
	\raisebox{0.3cm}{\begin{ytableau}
	3\\
	2&3\\
	1&2&3
	\end{ytableau}}.
\]
Then $\mathsf{row}(b) = 323122$, $\pi=526134$, $p=6$, $\pi_p=4$, $\nl=2$, and $\nr=0$. Hence $\cycle(b,2)=(6,5,4)$ written in cycle 
notation for permutations.
\end{example}

\begin{example}
Consider $b$ and $f_3b$ of Example~\ref{example.f3}. Then 
\[
	\TIKZ[scale=.3, xscale =1.1]{
	\node[left] at (.5,0) {$\mathsf{row}(b) = $};
	\foreach \w [count=\j from 1] in {4,3,3,4,2,2,2,3,3,4,4,4,4,1,1,1,1,2,2,3,3,3,4}{
		\ifthenelse
			{\w=3}
			{\node[bssyt3] (\j) at (\j,0) {\bf \w};}
			{\ifthenelse
				{\w=4}
				{\node[bssyt4] (\j) at (\j,0) {\bf \w};}
				{\node (\j) at (\j,0) {\w};}
				}
		}
	\node[ssyt3, rounded corners, draw, thick, inner sep=2pt] (j) at (9) {$\phantom{3}$};
	\draw[bssyt3, thick, -] (j) to +(.5,1.5) node[above, inner sep=1pt]{\footnotesize$p=9$};
		\Bracketed{1}{2}{0}	
		\Bracketed{4}{8}{0}	
		\Bracketed{11}{22}{2}
		\Bracketed{12}{21}{1}
		\Bracketed{13}{20}{0}
	\draw[rounded corners] (1.5,-1.1) node[v]{} to (1.5,-1.5)
	 	to node[below]{\footnotesize $\nl = 2$} 
		 (6, -1.5) to (6,-1.1) node[v]{};
	 \draw (16.5,-1.1) node[v]{} to (16.5,-1.1 - .4) node[v]{} to (16.5,-1.1-2*.4) node[v]{} node[below] {\footnotesize $\nr = 3$} ;
	}
	\]
	so that $p=9$, $\pi_p=13$, $\nl=2$, and $\nr=3$. Hence $\cycle(b,3)=(18,17,16,15,14,13)$.
\end{example}

\begin{lemma}
\label{lemma.sigma}
Suppose that $f_i(b) = b'$ for $b,b'  \in B(\lambda)_n$. Let 
\[ w = \mathsf{row}(b), \quad w' = \mathsf{row}(b'), \quad \pi=\mathsf{std}(w), \quad \pi'=\mathsf{std}(w'), \quad T=\std(b), \quad T'=\std(b'). \] Then
\[
	\cycle(b,i) \cdot \pi = \pi' \qquad \text{and} \qquad \cycle(b,i) \cdot T = T'.
\]
\end{lemma}

\begin{proof}
We use the same notation as in the definition of $\cycle(b,i)$. We compare the standardization $\pi$ and $\pi'$.
First note that all letters $1,2,\ldots,i-1$ and all letters $i$ to the left of position $p$ in $w$ and $w'$ standardize in the same way.
Similarly, the letters $i+1$ to the right of position $p$ and all letters bigger than $i+1$ in $w$ and $w'$ standardize in the same
way. Hence $\pi$ and $\pi'$ agree in these letters. It remains to compare the standardization of the remaining letters:
\begin{itemize}
\item the element $w_p=i$ standardizes to $\pi_p$; the element $w'_p=i+1$ standardizes to $\pi_p+\nr+\nl$ in $\pi'$;
\item the (bracketed) elements $i$ appearing to the right of $w_p$ are labeled $\pi_p+1,\pi_p+2,\ldots,\pi_p+\nr$ in $\pi$, whereas in $\pi'$ they 
are labeled $\pi_p,\pi_p+1,\ldots,\pi_p+\nr-1$;
\item the (bracketed) elements $i+1$ appearing to the left of $w_p$ are labeled $\pi_p+\nr+1,\pi_p+\nr+2,\ldots,\pi_p+\nr+\nl$ in $\pi$, whereas in 
$\pi'$ they are labeled $\pi_p+\nr,\pi_p+\nr+1,\ldots,\pi_p+\nr+\nl-1$.
\end{itemize}
Note that this difference is precisely explained by the action of $\cycle(b,i)$ as defined in \eqref{equation.sigma} on $\pi$ to obtain $\pi'$. 

Since $\pi = \row(T)$ and $\pi'=\row(T')$, it follows that also $\cycle(b,i) \cdot T = T'$.
\end{proof}

\begin{example}
Continuing Example~\ref{example.sigma b}, we have $w=323122$, $w'=323123$, $\pi=526134$ and $\pi'=425136$. Recall that 
$\cycle(b,2)=(654)$ and indeed $(654) \cdot 526134 = 425136$.
\end{example}

Since the standardization $\mathsf{std}(b)$ for $b \in B(\lambda)_n$ uniquely determines the quasi-crystal component of $b$, 
Proposition~\ref{proposition.quasi edges} follows immediately from Lemma~\ref{lemma.sigma} since by Remark~\ref{rmk:lengthofsigma}
the cycle $\cycle(b,i)$ is trivial if and only if $\nl=\nr=0$.

There are other tableaux which can index quasi-crystals. Instead of taking standard tableaux, one can index a quasi-crystal by its highest weight 
element.

\begin{definition}[\cite{AS.2017,Wang.2019}]
\label{lemma.bijection SYT QYT}
A semistandard Young tableau $T$ is a \defn{quasi-Yamanouchi tableau} if, for every $i>1$ that appears in the tableau, some instance
of $i$ is in a higher row than some instance of $i - 1$. Denote the set of all quasi-Yamanouchi tableaux of shape $\lambda$
by $\mathsf{QYT}(\lambda)$.
\end{definition}

\begin{lemma}[\cite{AS.2017}]
\label{lemma.QY}
The standardization map $\mathsf{std} \colon \mathsf{QYT}(\lambda) \to \mathsf{SYT}(\lambda)$ is a bijection such that
$\mathsf{wt}(T) = \mathsf{Des}(\mathsf{std}(T))$ for $T\in \mathsf{QYT}(\lambda)$.
\end{lemma}

The connection between quasi-crystals and quasisymmetric functions is as follows.
\defn{Gessel's fundamental quasisymmetric function}~\cite{Gessel.1984} is indexed by compositions $\alpha$
\[
	F_\alpha = \sum_{\substack{\beta \preccurlyeq \alpha \\ \text{refinement}}} M_\beta
	\qquad \text{with} \qquad 
	M_\beta = \sum_{i_1<i_2<\cdots<i_\ell} x_{i_1}^{\beta_1} x_{i_2}^{\beta_2} \cdots x_{i_\ell}^{\beta_\ell}
\]
and $\beta \preccurlyeq \alpha$ indicates that $\beta$ is a refinement of $\alpha$, that is, adjacent parts of $\beta$ can be summed to
obtain $\alpha$.

The character of the quasi-crystal $Q_T$ with $T\in\SYT(\lambda)$ inside $B(\lambda)_n$ is Gessel's quasisymmetric function $F_\alpha$
in $n$ variables, where $\alpha=\mathsf{Des}(T)$. Since $B(\lambda)_n$ is the union over all $Q_T$ with $T \in \SYT(\lambda)$ (as long as
$T$ appears in $B(\lambda)_n$, which is true for $n$ large enough), this yields Gessel's formula~\cite{Gessel.1984}
\[
    s_\lambda = \sum_{T \in \mathsf{SYT}(\lambda)} F_{\mathsf{Des}(T)} = \sum_{T \in \mathsf{QYT}(\lambda)} F_{\wt(T)}.
\]

\subsection{Dual equivalence graphs}
\label{section.DE}

Dual equivalence graphs were first introduced by Haiman~\cite{Haiman.1992}. The vertices of the \defn{dual equivalence graph} $\mathsf{DE}(\lambda)$
indexed by the partition $\lambda$ are all standard Young tableaux $\mathsf{SYT}(\lambda)$ of shape $\lambda$. 
The edges in the dual equivalence graph $\mathsf{DE}(\lambda)$ are given by the elementary dual equivalence relations $D_i$ ($1< i <
|\lambda|=:N$) defined on permutations as follows:
\begin{equation}
\label{equation.dual equivalence}
\begin{split}
	\ldots i \ldots i+1 \ldots i-1 \ldots &\stackrel{i}{\longleftrightarrow} \ldots i-1 \ldots i+1 \ldots i \ldots  \\
	\ldots i \ldots i-1 \ldots i+1 \ldots &\stackrel{i}{\longleftrightarrow} \ldots i+1 \ldots i-1 \ldots i \ldots  
\end{split}
\end{equation}
The operator $D_i$ is not defined for other configurations of the letters $i-1,i,i+1$ in the permutation.
Note that descents in the permutation do not change under $D_i$. Hence $D_i$ is defined on a standard Young tableau $T$ as well using the 
reading word $\mathsf{row}(T)$.

Assaf~\cite{Assaf.2008} proved that $D_i(T)$ on $T\in \mathsf{SYT}(\lambda)$ can be expressed in terms of crystal operators:
\begin{equation}
\label{equation.dual equivalence crystal}
	D_i(T) = \begin{cases}
	f_{i-1} f_i e_{i-1} e_i(T) & \text{if $e_i(T) \neq \emptyset$,}\\
	f_i f_{i-1} e_i e_{i-1}(T) & \text{if $e_{i-1}(T) \neq \emptyset$.}
	\end{cases}
\end{equation}
Note that if $D_i(T)$ is defined, then either $e_i(T) \neq \emptyset$ or $e_{i-1}(T) \neq \emptyset$, but not both. This can be seen by 
inspecting~\eqref{equation.dual equivalence}.

\section{Crystal skeletons}
\label{section.CS}

In~\cite{MG.2023}, Maas-Gari\'epy introduced the \defn{crystal skeleton} by contracting the quasi-crystals in $B(\lambda)_n$ to a vertex, assuming 
that $n$ is sufficiently large. Since there is a unique standard tableau $T\in \mathsf{SYT}(\lambda)$ in each quasi-crystal in $B(\lambda)_n$,
it is natural to label the vertices of the crystal skeleton by standard tableaux.

\begin{definition}
\label{definition.CS}
Let $\lambda$ be a partition and consider the ambient crystal $B(\lambda)_n$ for $n \geqslant |\lambda|$.

The \defn{crystal skeleton} $\mathsf{CS}(\lambda)$ is an edge-labeled, directed graph whose vertices are elements in $\mathsf{SYT}(\lambda)$. 
For $T,T' \in \mathsf{SYT}(\lambda)$, there is an edge from $T$ to $T'$ in $\mathsf{CS}(\lambda)$ if there exist $b\in Q_T$ and $b' \in Q_{T'}$ 
such that $f_i(b) = b'$ for some $1\leqslant i < n$.
\end{definition}

\noindent
This definition is rather indirect. One goal of our work is to give a combinatorial description of $\CS(\lambda)$ as a graph
without reference to the original crystal $B(\lambda)_n$.
In addition, we will give a natural way to label the edges with intervals (described below) different from the edge-labels in~\cite{MG.2023}.
See Figure~\ref{figure.B21} for an example of a crystal skeleton.

\subsection{Vertices of the crystal skeleton}
\label{section.vertices}
As discussed in Definition \ref{definition.CS}, the vertices of $\CS(\lambda)$ are labeled by standard tableaux in $\SYT(\lambda)$. 
For our axiomatic description of crystal skeletons in Section~\ref{section.axioms}, it will be important to also associate to each vertex $T \in \SYT(\lambda)$
its descent composition $\Des(T)$ as defined in~\eqref{equation.descent composition}.

\begin{example}
\label{example.descentcomposition}
We color $T \in \SYT(3,2,1)$ by its descent composition below: 
\[ T = 
	\TIKZ[scale=.4]{
		\SSYTcolors{{1,2,4},{2,3},{4}}
		\Tableau{{1,3,6},{2,4},{5}}
		}
. \]
In particular, $\Des(T) = (1,2,1,2)$. 
\end{example}

\subsection{Edges of the crystal skeleton}
\label{section.edges}
The edges of the crystal skeleton are more subtle. In~\cite{MG.2023}, the edge from $T$ to $T'$ in $\mathsf{CS}(\lambda)$ is indexed by the minimal 
index $j$ such that $f_j(b)=b'$ for $b\in Q_T$ and $b'\in Q_{T'}$. We give two alternative edge labels in $\CS(\lambda)$.

\subsubsection{Dyck pattern intervals}
\label{section.Dyck pattern}

The next definition is crucial for our labeling of the edges in the crystal skeleton.

\begin{definition}
\label{definition.Dyck pattern}
Let $I=[i,i+2m] \subseteq [n]$ be an interval of length $2m+1\geqslant 3$.
\begin{enumerate}
\item Let $\pi \in S_n$ be a permutation and $\pi|_I$ be the subword of $\pi$ restricted to the letters in $I$. 
We call $I$ a \defn{Dyck pattern interval} of $\pi$ if
\begin{equation}
\label{equation.P pi I}
	P(\pi|_I) = \scalebox{0.6}{\TIKZ[scale=1.5]{
	\Tableau{{i,i+1,\ldots,i+m-1,i+m},{i+m+1,i+m+2,\ldots,i+2m}}
	}}\, ,
\end{equation}
that is, the letters in $[i,i+m]$ occur in the bottom row and the letters in $[i+m+1,i+2m]$ occur in the top row of $P(\pi|_I)$.
We call the subword $\pi|_I$ a \defn{Dyck pattern} of $\pi$. \medskip
\item Let $T \in \mathsf{SYT}(\lambda)$ and $n=|\lambda|$. Then $I$ is a \defn{Dyck pattern interval} of $T$ if it is a Dyck pattern interval of 
$\pi=\row(T)$. Similarly, $\pi|_I$ is a \defn{Dyck pattern} of $T$ if it is a Dyck pattern of $\pi$.
\end{enumerate}
\end{definition}

\begin{example}
Consider $T$ of Example~\ref{example.descentcomposition}, so that $\pi=\row(T)=524136$. The interval $I=[2,4]$ is a Dyck pattern interval of $T$ since 
$\pi|_{[2,4]}=243$ and 
\[
	P(\pi|_{[2,4]}) = 
		\TIKZ[scale=.4]{
		\Tableau{{2,3},{4}}
		}.
\]
\end{example}

\begin{example}
\label{example.pi}
Let $\pi=10 \, 783142596 \in S_{10}$. Then $I=[3,9]$ is a Dyck pattern interval of $\pi$ since $\pi|_I = 7834596$ and
\[
	P(\pi|_I) = 
		\TIKZ[scale=.4]{
		\Tableau{{3,4,5,6},{7,8,9}}
		}\, .
\]
\end{example}

We may think about a Dyck pattern interval $I=[i,i+2m]$ also in terms of crystal operators. To this end, we first define the \defn{destandardization} 
of $\pi|_I$ as follows. Note that the letters in $[i,i+m]$ (resp. $[i+m+1,i+2m]$) need to appear in increasing order in $\pi|_I$ due to the form of 
$P(\pi|_I)$ in Definition \ref{definition.Dyck pattern}. Replacing the letters in $[i,i+m]$ by $i$ and the letters in $[i+m+1,i+2m]$ by $i+1$ in $\pi|_I$, 
we obtain a word in the letters $i$ and $i+1$ such that all letters $i+1$ are bracketed with some $i$. Since there is one more letter $i$ than $i+1$,
there is one unbracketed letter $i$. We call this word the \defn{destandardization} of $\pi|_I$ denoted
\begin{equation}
\label{equation.destd}
	\destd(\pi|_I).
\end{equation}

\begin{remark}
\label{remark.Dyck pattern}
An alternate definition of a Dyck pattern interval $I=[i,i+2m]$ of $\pi$ is as follows. Suppose $m\geqslant 1$ and the letters in $[i,i+m]$ 
(resp. $[i+m+1,i+2m]$) appear in increasing order in $\pi|_I$. As before, let $w=\destd(\pi|_I)$ be the word by replacing the letters in $[i,i+m]$ by 
$i$ and the letters in $[i+m+1,i+2m]$ by $i+1$. Then $I$ is a Dyck pattern interval if 
\[
	e_i(w)=\emptyset \quad \text{and} \quad f_i(i+1 \, w)=\emptyset.
\]
The condition $e_i(w)=\emptyset$ ensures that all letters $i+1$ are bracketed with some $i$. 
Viewing the letters $i+1$ as down steps and $i$ as up steps, the word $i+1\, w$ associated to a Dyck pattern is a Dyck path. This motivates the
name \defn{Dyck pattern}.
\end{remark}

\begin{remark}
By the RSK correspondence, the Dyck pattern $\pi|_I$ is uniquely determined by the insertion and recording tableaux $P(\pi|_I)$ and
$Q(\pi|_I)$. Since $P(\pi|_I)$ is fixed by definition, a Dyck pattern is hence determined by the recording tableau, which is a standard tableau
of shape $(m+1,m)$ or equivalently a standard tableau of shape $(m+1,m+1)$ since the entry in the unique corner cell $(2,m+1)$ is fixed to be $2m+2$.
It is well-known that standard tableaux of shape $(m+1,m+1)$ are in bijection with Dyck paths of length $2m+2$ (see~\cite[pg. 243]{Stanley.EC2}).
\end{remark}

\begin{example}
Consider $\pi$ and $I$ from Example~\ref{example.pi}. Then 
\[
	\TIKZ[scale=.3, xscale =.9]{
	\node[left] at (.5,0.1) {$\pi|_I = \viridian{78}\sapphire{345}\viridian{9}\sapphire{6}$
	\quad and \quad $w=\destd(\pi|_I)=$};
	\foreach \w [count=\j from 1] in {4,4,3,3,3,4,3}{
		\ifthenelse
			{\w=4}
			{\node[viridian] (\j) at (\j,0) {\w\strut};
			\node[viridian] (\j) at (\j,1.25) {$\scriptstyle ($\strut};}
			{\ifthenelse
				{\w=3}
				{\node[sapphire] (\j) at (\j,0) {\w\strut};
				\node[sapphire] (\j) at (\j,1.25) {$\scriptstyle )$\strut};}
				{\node (\j) at (\j,0) {\w\strut};}
				}
		}
		\node[right] at (7,0.1){\strut.};
		\draw [rounded corners, thick](1, -.6) node[v,viridian]{} to (1,-1.2) to  (4, -1.2) to (4,-.6) node[v,sapphire]{};
		\draw [bend right=60, thick](2, -.6) node[v,viridian]{} to  (3,-.6) node[v,sapphire]{};
		\draw [bend right=60, thick](6, -.6) node[v,viridian]{} to  (7,-.6) node[v,sapphire]{};
	}
\]
Note that $e_3(w)=\emptyset$ and $f_3(4w)=\emptyset$. The Dyck path corresponding to $4w$ is
\[
\TIKZ[scale=.3]{
\draw[very thick, viridian] (0,0) to ++(3,-3);
\draw[very thick, sapphire] (3,-3) to ++(3,3);
\draw[very thick, viridian] (6,0) to ++(1,-1);
\draw[very thick, sapphire] (7,-1) to ++(1,1);
\draw (0,0) to ++(8,0);
}
\]
\end{example}

We show that the edges of $\CS(\lambda)$ are labeled by Dyck pattern intervals.

\begin{theorem}
\label{theorem.Dyck pattern}
There is a bijection between the edges in $\mathsf{CS}(\lambda)$ and Dyck pattern intervals that occur in $\row(T)$ for  $T \in \mathsf{SYT}(\lambda)$. 
\end{theorem}

\begin{example}
\label{example.shortpattern}
We describe in detail the edge below found in $\CS(3,2,1)$ in Figure \ref{figure.CS321}:
\[ T=   
\TIKZ[scale=.4]{
		\SSYTcolors{{1,1,1},{2,2},{3}}
		\Tableau{{1,2,3},{4,5},{6}}
		}
    \quad \xrightarrow[]{~I = [1,5]~} \quad T' =  
\TIKZ[scale=.4]{
		\SSYTcolors{{1,1,2},{2,2},{3}}
		\Tableau{{1,2,5},{3,4},{6}}
		}
. \]
In this case $\pi = \row(T) = 645123$, $\pi|_{[1,5]} = 45123$, and
\[
	P(\pi|_I) = 
	\TIKZ[scale=.4]{
		\SSYTcolors{{1,1,1},{2,2}}
		\Tableau{{1,2,3},{4,5}}
		}.
\]
Thus $I$ is a Dyck pattern interval on $\pi$. 

To see why there is an edge from $T$ to $T'$ in $\CS(3,2,1)$, note that 
the word $w = 622111$ is in the same quasi-crystal component of $B(3,2,1)_6$ as $\pi$, since $\std(w) = \pi$. Similarly, $w' = 622112$ is in the same 
quasi-crystal component as $\pi' := \row(T') = 634125$ in $B(3,2,1)_6$, since $\std(w') = \pi'$. Next, observe that in $B(3,2,1)_6$, we have
\[ 
 	f_{1}(w) = f_1(62211\textcolor{viridian}{1})  = 62211\textcolor{viridian}{2} = w'. 
\]
The idea behind Theorem \ref{theorem.Dyck pattern} is that the Dyck pattern interval detects this edge in $B(\lambda)_n$. Note that $\Des(T) = (3,2,1)$ 
and $\Des(T') = (2,3,1)$. The way descent compositions change between edges in $\CS(\lambda)$ will be described in 
Section~\ref{sec:descent compositions}.
\end{example}

\begin{proof}[Proof of Theorem~\ref{theorem.Dyck pattern}]
Suppose that $I=[i,i+2m]$ with $m\geqslant 1$ is a Dyck pattern interval in $\pi = \row(T)$. Replace the letters $i,i+1,\ldots,i+m$ by $i$ and the letters
$i+m+1,i+m+2,\ldots,i+2m$ by $i+1$ in $T$ to obtain a new semistandard tableau $b$. Note that since $I$ is a Dyck pattern interval,
the letters $i,i+1,\ldots,i+m$ (resp. $i+m+1,i+m+2,\ldots,i+2m$) have to appear left to right in $T$. Hence we have $\std(b)=T$. Furthermore, 
since $b$ contains one more $i$ than $i+1$ not all $i$ are bracketed and hence $f_ib=b'$ exists.
In addition, since $I$ is a Dyck pattern interval, by Remark~\ref{remark.Dyck pattern} all $i+1$ are 
bracketed with an $i$. Hence by Proposition~\ref{proposition.quasi edges} since $m\geqslant 1$, the standard tableau $T'=\std(b')$ is different
from $T$ and there is an edge from $T$ to $T'$ in $\CS(\lambda)$. This shows that every Dyck pattern interval gives rise to an edge
in $\CS(\lambda)$.

Now suppose that there is an edge from $T$ to $T'$ in $\CS(\lambda)$. This means that there exist $b\in Q_T$ and $b'\in Q_{T'}$ such that
$f_ib=b'$ for some $i$. By the explanations in Section~\ref{section.QS} before~\eqref{equation.sigma} and using the same notation,
there are $\nr$ bracketed pairs $(i+1,i)$ after $w_p$ and $\nl$ bracketed pairs $(i+1,i)$ before $w_p$. Restricting $\pi=\row(T)$ to the 
interval $I=[\pi_p-\nl,\pi_p+\nl+2\nr]$ yields a Dyck pattern of size $2m+1$ with $m=\nl+\nr>0$ by Remark~\ref{remark.Dyck pattern}.
Hence every edge in $\CS(\lambda)$ leaving $T$ is associated to a Dyck pattern interval $I$ of $\pi$.
\end{proof}

\subsubsection{Cycles}\label{section:cycles}

By Theorem~\ref{theorem.Dyck pattern}, every edge leaving $T$ in $\mathsf{CS}(\lambda)$ is associated with a Dyck pattern
interval $I=[i,i+2m]$ with $m\geqslant 1$ of $\pi=\row(T)$. As in the proof of Theorem~\ref{theorem.Dyck pattern}, let $b$ be the 
semistandard tableau obtained from $T$ by replacing the letter $i,i+1,\ldots,i+m$ by $i$ and the letters
$i+m+1,i+m+2,\ldots,i+2m$ by $i+1$. Note that $b|_I$ is $\destd(\pi|_I)$.
By Lemma~\ref{lemma.sigma}, the edge leaving $T$ associated to $I$ goes to
\begin{equation}
\label{equation.cycle on T}
	T' = (m+\pi_p,m+\pi_p-1,\ldots,\pi_p) \cdot T,
\end{equation}
where $\pi_p$ is the letter in $T$ (or $\pi$) which corresponds to the $i$ in $b$ on which $f_i$ acts. Hence the edge from $T$ to
$T'$ in $\CS(\lambda)$ can also be labeled by a \defn{cycle} which we shall denote by 
\[ 
	\cycle(\pi|_I):=(m+\pi_p,m+\pi_p-1,\ldots,\pi_p). 
\]

\begin{example}
Continuing Example~\ref{example.shortpattern}, the semistandard tableau $b$ associated to $T$ is
\[ 
b=   
\TIKZ[scale=.4]{
		\SSYTcolors{{1,1,1},{2,2},{3}}
		\Tableau{{1,1,1},{2,2},{6}}
		}.
\]
The crystal operator $f_1$ acts on the rightmost 1 and hence $\pi_p=3$. Since $|I|=5$, we have $m=2$, and $\cycle(\pi|_I) = (5,4,3)$, and indeed 
\[ (5,4,3)\cdot T = T'. \]
Figure~\ref{figure.CS321} shows the crystal skeleton $\CS(3,2,1)$ with edges labeled by both Dyck pattern intervals and cycles.
\end{example}

As we will see later, due to self-similarity properties which respect the Dyck pattern intervals, but not the cycles, it is more natural
to use the Dyck pattern intervals as edge labels.

\section{Properties of the crystal skeleton}
\label{section.CS properties}

In this section, we study various properties of the crystal skeleton $\CS(\lambda)$.

\subsection{Dual equivalence graphs as subgraphs}

Recall the dual equivalence graphs of Section~\ref{section.DE}. The next theorem proves~\cite[Conjecture~5.3]{MG.2023}.

\begin{theorem}\label{thm:DEinsideCS}
The dual equivalence graph $\mathsf{DE}(\lambda)$ is a subgraph of the crystal skeleton $\mathsf{CS}(\lambda)$ (disregarding edge labels and edge 
directions).
\end{theorem}

\begin{proof}
Suppose that there is an edge labeled $i$ in $\mathsf{DE}(\lambda)$ between $T$ and $T'$, that is, $D_i(T)=T'$. 
By~\eqref{equation.dual equivalence crystal}, either $T' = f_{i-1} f_i e_{i-1} e_i(T)$ or $T' = f_i f_{i-1} e_i e_{i-1}(T)$.
Suppose first that $T' = f_{i-1} f_i e_{i-1} e_i(T)$. In this case, $i$ is to the left of $i+1$ in $\mathsf{row}(T)$, so that
\begin{equation}
\label{equation.DE1}
\begin{split}
	\ldots i \ldots i+1 \ldots i-1 \ldots \stackrel{e_i}{\longrightarrow} \ldots i \ldots i \ldots i-1 \ldots 
	&\stackrel{e_{i-1}}{\longrightarrow} \ldots i-1 \ldots i \ldots i-1 \ldots \\
	\stackrel{f_i}{\longrightarrow} \ldots i-1 \ldots i+1 \ldots i-1 \ldots
	&\stackrel{f_{i-1}}{\longrightarrow} \ldots i-1 \ldots i+1 \ldots i \ldots
\end{split}
\end{equation}
or
\begin{equation}
\label{equation.DE2}
\begin{split}
	\ldots i \ldots i-1 \ldots i+1 \ldots \stackrel{e_i}{\longrightarrow} \ldots i \ldots i-1 \ldots i \ldots 
	&\stackrel{e_{i-1}}{\longrightarrow} \ldots i \ldots i-1 \ldots i-1 \ldots \\
	\stackrel{f_i}{\longrightarrow} \ldots i+1 \ldots i-1 \ldots i-1 \ldots
	&\stackrel{f_{i-1}}{\longrightarrow} \ldots i+1 \ldots i-1 \ldots i \ldots
\end{split}
\end{equation}
In both cases, the crystal operator $e_{i-1}$ changes quasi-crystal components by Proposition~\ref{proposition.quasi edges}
since there is a bracketing between $i$ and $i-1$, whereas all other crystal operators in this sequence do not change
quasi-crystal components. Hence there is an edge in $\mathsf{CS}(\lambda)$ between $T$ and $T'$.

If, on the other hand, $T' = f_i f_{i-1} e_i e_{i-1}(T)$, the letter $i-1$ must be to the left of the letter $i$. In this case, we start
with the last elements in the sequences of~\eqref{equation.DE1} and~\eqref{equation.DE2}. Reversing the arrows (changing $e$ to $f$ and vice versa),
the dual equivalence operator $D_i$ takes the last element to the first element. The quasi-crystal component again changes in the same arrow.
Hence there is an edge in $\mathsf{CS}(\lambda)$ between $T$ and $T'$. This proves that $\mathsf{DE}(\lambda)$ is a subgraph of $\mathsf{CS}(\lambda)$.
\end{proof}

We can characterize which edges in $\CS(\lambda)$ belong to $\DE(\lambda)$. Namely, they are the edges labeled by Dyck pattern intervals $I$ of size 3, or
equivalently, by transpositions when labeling by cycles.

\begin{proposition}
\label{proposition.transposition}
Let $T \stackrel{I}{\longrightarrow} T'$ be an edge in $\CS(\lambda)$.
Then there is an edge between $T$ and $T'$ in the dual equivalence graph $\DE(\lambda)$ if and only if $|I|=3$.
\end{proposition}

\begin{proof}
All elementary dual equivalence relations $D_i$ are simple transpositions on the corresponding standard tableaux $T$ and $T'$. Hence in the
notation of~\eqref{equation.cycle on T} we have $m=1$ or equivalently $|I|=3$.

Conversely, assume that $|I|=3$, so that $I=[i-1,i+1]$ for some $i$. Since $I$ is a Dyck pattern interval by Theorem~\ref{theorem.Dyck pattern}, 
we must have for $\pi=\row(T)$ that $\pi|_I = i-1 \ i+1 \ i$ with transposition $(i,i-1)$ or $\pi|_I=i+1\ i-1 \ i$ with transposition $(i+1,i)$. 
These correspond to the dual equivalence transitions on the top and bottom (right to left) in~\eqref{equation.dual equivalence}. This proves the claim.
\end{proof}

\subsection{Properties of Dyck pattern intervals}
\label{section.properties}

We collect several properties of Dyck pattern intervals which will be useful later. 

First, there is an inclusion property of Dyck pattern intervals.

\begin{proposition}
\label{proposition.inside}
Let $T \in \SYT(\lambda)$ and suppose that $I=[i,i+2m]$ is a Dyck pattern interval of $T$ of size $2m+1>3$ (or $m>1$). Then
$J=[i+1,i+2m-1]$ of size $2m-1$ is also a Dyck pattern interval of $T$.
\end{proposition}

\begin{proof}
Let $\pi=\row(T)$. Since $I$ is a Dyck pattern interval, the insertion tableau $P(\pi|_I)$ is of the form~\eqref{equation.P pi I} of
shape $(m+1,m)$. The insertion tableau $P(\pi|_J)$ can be obtained from $P(\pi|_I)$ by removing the cells containing $i$ and $i+2m$ and 
performing jeu de taquin to straighten the tableau. This yields the tableau of shape $(m,m-1)$ with the letters $i+1,i+2,\ldots,i+m$ 
(resp. $i+m+1,i+m+2,\ldots,i+2m-1$) in the bottom (resp. top) row. Hence $J$ is a Dyck pattern interval.
\end{proof}

\begin{example}
Consider the edge in Example~\ref{example.shortpattern} labeled by the Dyck pattern interval $I=[1,5]$. As can be seen from Figure~\ref{figure.CS321},
there is also an edge starting at $T$ labeled by the interval $[2,4]$.
\end{example}

Second, recall from Proposition~\ref{proposition.fi and Knuth} that crystal operators respect Knuth equivalence. The same is true for edges in the crystal
skeleton.

\begin{proposition}
Suppose $\pi \equiv_K \sigma$ are two Knuth equivalent permutations. Let $I$ be a Dyck pattern interval in $\pi$. Then $I$ is also a Dyck pattern interval
in $\sigma$. Furthermore, if 
\[ \pi \stackrel{I}{\longrightarrow} \pi' \quad \textrm{ and } \quad \sigma \stackrel{I}{\longrightarrow} \sigma' \] 
are edges in the corresponding crystal
skeletons, then $\pi' \equiv_K \sigma'$.
\end{proposition}

\begin{proof}
Since $I=[i,i+2m]$ is a Dyck pattern interval of $\pi$, $P(\pi|_I)$ is of the form~\eqref{equation.P pi I}.
It is not hard to check from~\eqref{equation.Knuth} that if $\pi \equiv_K \sigma$ then also $\pi|_I \equiv_K \sigma|_I$. Hence,
$P(\sigma|_I)=P(\pi|_I)$, so that $I$ is also a Dyck pattern interval in $\sigma$, proving the first claim.

By~\cite[Theorem 8.4]{BumpSchilling.2017}, the crystals containing $\pi$ and $\sigma$ are isomorphic since $\pi \equiv_K \sigma$.
Since the crystal skeleton is fully determined by the underlying crystal and its edges, the second claim follows from Proposition~\ref{proposition.fi and Knuth}.
\end{proof}

\subsection{Self-similarity and branching property}
\label{section.self similar}

In this section, we study self-similarity properties of the crystal skeleton. For a standard tableau $T\in \mathsf{SYT}(\lambda)$ and an interval 
$[a,b]$, we define $T_{[a,b]}$ to be the skew tableau $T$ restricted to the interval $[a,b]$. Furthermore, for a skew tableau $T$, let
$\mathsf{jdt}(T)$ be the jeu de taquin straightening of $T$ (see for example~\cite[Chapter A1.2]{Stanley.EC2} and~\cite[Chapter 3.7]{Sagan.2001}).

\begin{example}
\[
\text{Let} \quad 
	T=
	\raisebox{0.3cm}{
	\begin{ytableau}
	6\\
	4&5\\
	1&2&3
	\end{ytableau}}.
	\quad
	\text{Then} \quad
	T_{[2,6]} = 
	\raisebox{0.3cm}{
	\begin{ytableau}
	6\\
	4&5\\
	*(gray) &2&3
	\end{ytableau}}
	\quad
	\text{and}
	\quad
	\mathsf{jdt}(T_{[2,6]})=
	\raisebox{0.3cm}{
	\begin{ytableau}
	6\\
	4&5\\
	2&3
	\end{ytableau}}\, .
\]
\end{example}

The self-similarity property for crystal skeletons can be stated as follows. 

\begin{theorem}
\label{theorem.self similar}
Let $T\in \mathsf{SYT}(\lambda)$, $[a,b]\subseteq [1,|\lambda|]$ be an interval, and $\mu=\mathsf{shape}(\mathsf{jdt}(T_{[a,b]}))$.
Then $\mathsf{CS}(\mu)$ is a subgraph of $\mathsf{CS}(\lambda)$ (up to relabelling of the edges). In particular, labeling the edges in the crystal skeleton by
the Dyck pattern intervals, the edges in $\mathsf{CS}(\mu)$ differ from those in $\mathsf{CS}(\lambda)$ by a shift of $-a+1$.
\end{theorem}

Before proving Theorem \ref{theorem.self similar}, we give an example.

\begin{example}
Recall that Figure~\ref{figure.CS321} shows $\CS(3,2,1)$. Take the subgraph of $\CS(3,2,1)$ with the vertices
\begin{equation}
\label{equation.tableaux221}
	\raisebox{0.3cm}{
	\begin{ytableau}
	5\\
	3&4\\
	1&2&6
	\end{ytableau}},
	\quad
	\raisebox{0.3cm}{
	\begin{ytableau}
	4\\
	3&5\\
	1&2&6
	\end{ytableau}},
	\quad
	\raisebox{0.3cm}{
	\begin{ytableau}
	5\\
	2&4\\
	1&3&6
	\end{ytableau}},
	\quad
	\raisebox{0.3cm}{
	\begin{ytableau}
	4\\
	2&5\\
	1&3&6
	\end{ytableau}},
	\quad
	\raisebox{0.3cm}{
	\begin{ytableau}
	3\\
	2&5\\
	1&4&6
	\end{ytableau}}.
\end{equation}
It is isomorphic to the crystal skeleton $\mathsf{CS}(2,2,1)$ in Figure~\ref{figure.subgraphs}. Note that all tableaux in~\eqref{equation.tableaux221}
restricted to the interval $[1,5]$ have shape $(2,2,1)$.
\begin{figure}
\begin{center}
\definecolor{permutation}{HTML}{2660A4}
\scalebox{0.8}{
\begin{tikzpicture}[>=latex,line join=bevel,]
\node (node_0) at (88.361bp,341.0bp) [draw,draw=none] {${\def\lr#1{\multicolumn{1}{|@{\hspace{.6ex}}c@{\hspace{.6ex}}|}{\raisebox{-.3ex}{$#1$}}}\raisebox{-.6ex}{$\begin{array}[t]{*{2}c}\cline{1-1}\lr{5}\\\cline{1-2}\lr{3}&\lr{4}\\\cline{1-2}\lr{1}&\lr{2}\\\cline{1-2}\end{array}$}}$};
  \node (node_1) at (30.361bp,235.0bp) [draw,draw=none] {${\def\lr#1{\multicolumn{1}{|@{\hspace{.6ex}}c@{\hspace{.6ex}}|}{\raisebox{-.3ex}{$#1$}}}\raisebox{-.6ex}{$\begin{array}[t]{*{2}c}\cline{1-1}\lr{4}\\\cline{1-2}\lr{3}&\lr{5}\\\cline{1-2}\lr{1}&\lr{2}\\\cline{1-2}\end{array}$}}$};
  \node (node_2) at (117.36bp,235.0bp) [draw,draw=none] {${\def\lr#1{\multicolumn{1}{|@{\hspace{.6ex}}c@{\hspace{.6ex}}|}{\raisebox{-.3ex}{$#1$}}}\raisebox{-.6ex}{$\begin{array}[t]{*{2}c}\cline{1-1}\lr{5}\\\cline{1-2}\lr{2}&\lr{4}\\\cline{1-2}\lr{1}&\lr{3}\\\cline{1-2}\end{array}$}}$};
  \node (node_3) at (30.361bp,129.0bp) [draw,draw=none] {${\def\lr#1{\multicolumn{1}{|@{\hspace{.6ex}}c@{\hspace{.6ex}}|}{\raisebox{-.3ex}{$#1$}}}\raisebox{-.6ex}{$\begin{array}[t]{*{2}c}\cline{1-1}\lr{4}\\\cline{1-2}\lr{2}&\lr{5}\\\cline{1-2}\lr{1}&\lr{3}\\\cline{1-2}\end{array}$}}$};
  \node (node_4) at (30.361bp,23.0bp) [draw,draw=none] {${\def\lr#1{\multicolumn{1}{|@{\hspace{.6ex}}c@{\hspace{.6ex}}|}{\raisebox{-.3ex}{$#1$}}}\raisebox{-.6ex}{$\begin{array}[t]{*{2}c}\cline{1-1}\lr{3}\\\cline{1-2}\lr{2}&\lr{5}\\\cline{1-2}\lr{1}&\lr{4}\\\cline{1-2}\end{array}$}}$};
  
     \draw[->] (node_0)  to[bend right=20] node[opLabel]{$[1,3]$}  (node_2);
     \draw[->] (node_0)  to[bend right=20] node[opLabel]{$[3,5]$}  node[pos=.6, above left, permutation]{$(54)$} (node_1);
     \draw[->] (node_2)  to[bend right=30] node[opLabel]{$[2,4]$}  node[pos=.5, above left, permutation]{$(32)$} (node_0);
     \draw[->] (node_1)  to node[opLabel]{$[1,3]$}  node[pos=.8, above left, permutation]{$(32)$} (node_3);
     \draw[->] (node_3)  to[bend right=30] node[opLabel]{$[2,4]$}  node[pos=.5, below right, permutation]{$(43)$} (node_4);
     \draw[->] (node_4)  to[bend right=20] node[opLabel]{$[3,5]$}  (node_3);
\end{tikzpicture}
}
\hspace{3cm}
\scalebox{0.8}{
\begin{tikzpicture}[>=latex,line join=bevel,]
\node (node_0) at (64.929bp,299.0bp) [draw,draw=none] {${\def\lr#1{\multicolumn{1}{|@{\hspace{.6ex}}c@{\hspace{.6ex}}|}{\raisebox{-.3ex}{$#1$}}}\raisebox{-.6ex}{$\begin{array}[t]{*{3}c}\cline{1-2}\lr{4}&\lr{5}\\\cline{1-3}\lr{1}&\lr{2}&\lr{3}\\\cline{1-3}\end{array}$}}$};
  \node (node_1) at (31.929bp,111.0bp) [draw,draw=none] {${\def\lr#1{\multicolumn{1}{|@{\hspace{.6ex}}c@{\hspace{.6ex}}|}{\raisebox{-.3ex}{$#1$}}}\raisebox{-.6ex}{$\begin{array}[t]{*{3}c}\cline{1-2}\lr{3}&\lr{4}\\\cline{1-3}\lr{1}&\lr{2}&\lr{5}\\\cline{1-3}\end{array}$}}$};
  \node (node_2) at (117.93bp,111.0bp) [draw,draw=none] {${\def\lr#1{\multicolumn{1}{|@{\hspace{.6ex}}c@{\hspace{.6ex}}|}{\raisebox{-.3ex}{$#1$}}}\raisebox{-.6ex}{$\begin{array}[t]{*{3}c}\cline{1-2}\lr{2}&\lr{5}\\\cline{1-3}\lr{1}&\lr{3}&\lr{4}\\\cline{1-3}\end{array}$}}$};
  \node (node_3) at (59.929bp,17.0bp) [draw,draw=none] {${\def\lr#1{\multicolumn{1}{|@{\hspace{.6ex}}c@{\hspace{.6ex}}|}{\raisebox{-.3ex}{$#1$}}}\raisebox{-.6ex}{$\begin{array}[t]{*{3}c}\cline{1-2}\lr{2}&\lr{4}\\\cline{1-3}\lr{1}&\lr{3}&\lr{5}\\\cline{1-3}\end{array}$}}$};
  \node (node_4) at (117.93bp,205.0bp) [draw,draw=none] {${\def\lr#1{\multicolumn{1}{|@{\hspace{.6ex}}c@{\hspace{.6ex}}|}{\raisebox{-.3ex}{$#1$}}}\raisebox{-.6ex}{$\begin{array}[t]{*{3}c}\cline{1-2}\lr{3}&\lr{5}\\\cline{1-3}\lr{1}&\lr{2}&\lr{4}\\\cline{1-3}\end{array}$}}$};
  
    \draw[->] (node_0)  to[bend right=20] node[opLabel]{$[1,5]$}  node[pos=.6, above left, permutation]{$(543)$} (node_1);
    \draw[->] (node_0)  to[bend right=30] node[opLabel]{$[2,4]$}  node[pos=.4, below right, permutation]{$(43)$} (node_4);
    \draw[->] (node_4)  to[bend right=20] node[opLabel]{$[3,5]$}  (node_0);
    \draw[->] (node_1)  to[bend right=20] node[opLabel]{$[1,3]$}  (node_3);
    \draw[->] (node_3)  to[bend right=30] node[opLabel]{$[2,4]$}  node[pos=.6, above left, permutation]{$(32)$} (node_1);
    \draw[->] (node_2)  to[bend left=20] node[opLabel]{$[3,5]$}  node[pos=.6, below right, permutation]{$(54)$} (node_3);
    \draw[->] (node_4)  to node[opLabel]{$[1,3]$}  node[pos=.6, right, permutation]{$(32)$} (node_2);
\end{tikzpicture}
}
\end{center}
\caption{Left: The crystal skeleton $\mathsf{CS}(2,2,1)$. Right: The crystal skeleton $\mathsf{CS}(3,2)$.
\label{figure.subgraphs}
}
\end{figure}
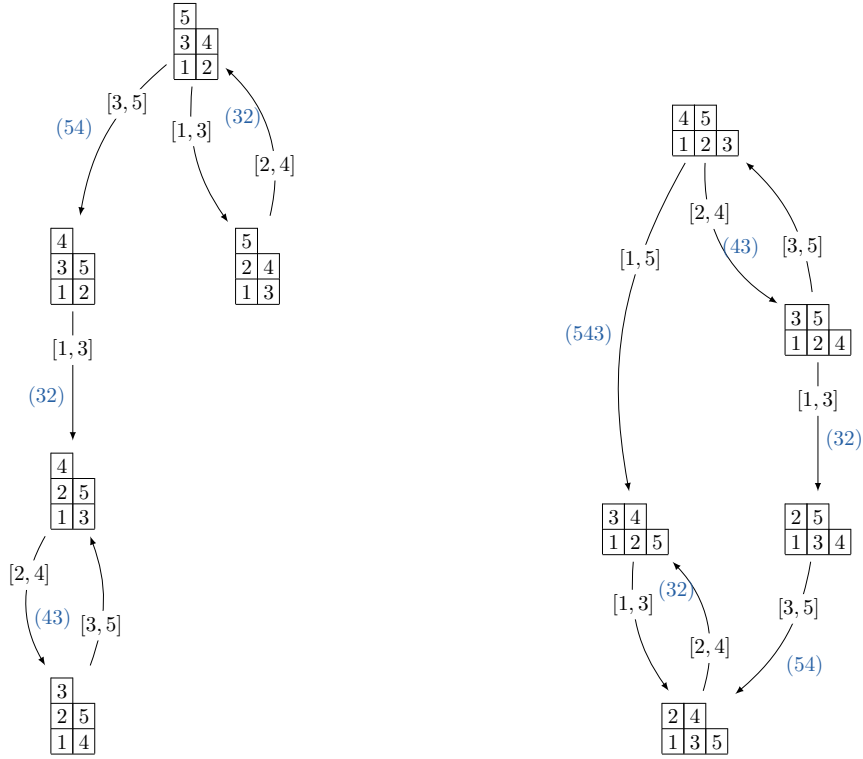

Take the subgraph of Figure~\ref{figure.CS321} with the vertices
\begin{equation}
\label{equation.tableaux32}
	\raisebox{0.3cm}{
	\begin{ytableau}
	5\\
	2&6\\
	1&3&4
	\end{ytableau}},
	\quad
	\raisebox{0.3cm}{
	\begin{ytableau}
	4\\
	3&6\\
	1&2&5
	\end{ytableau}},
	\quad
	\raisebox{0.3cm}{
	\begin{ytableau}
	4\\
	2&5\\
	1&3&6
	\end{ytableau}},
	\quad
	\raisebox{0.3cm}{
	\begin{ytableau}
	3\\
	2&6\\
	1&4&5
	\end{ytableau}},
	\quad
	\raisebox{0.3cm}{
	\begin{ytableau}
	3\\
	2&5\\
	1&4&6
	\end{ytableau}}.
\end{equation}
It is isomorphic to the crystal skeleton $\mathsf{CS}(3,2)$ in Figure~\ref{figure.subgraphs} up to shifting all intervals by $-1$.
Note that applying jeu de taquin to all tableaux in~\eqref{equation.tableaux32} restricted to the interval $[2,6]$ yields tableaux of shape~$(3,2)$.
\end{example}

\begin{proof}[Proof of Theorem~\ref{theorem.self similar}]
The row reading word of a skew semistandard tableau $P$ and the row reading word of $\mathsf{jdt}(P)$ are Knuth equivalent (see for 
example~\cite{Stanley.EC2}). Knuth equivalent words lie in isomorphic crystals (see for example~\cite[Chapter 8]{BumpSchilling.2017}).

Restricting the standard tableau $T \in \mathsf{CS}(\lambda)$ to the subinterval $[a,b]$, we can embed $T$ into the subcrystal
generated by the crystal operators $e_i,f_i$ for $a\leqslant i<b$. Shifting all letters in $T_{[a,b]}$ by $-a+1$ and operating with the 
crystal operators $f_i,e_i$ for $1\leqslant i \leqslant b-a$, $T_{[a,b]}-a+1$ sits inside a crystal isomorphic to $B(\mu)$ by the above arguments.
Considering the crystal skeleton of $B(\mu)$ proves the claim.
\end{proof}

\begin{remark}
Note that Knuth relations and jeu de taquin do not preserve cycles, even before shifting the labels of the corresponding tableaux. For example
\[
	\raisebox{0.2cm}{
	\begin{ytableau}
	3&4&6\\
	1&2&5
	\end{ytableau}}
	\quad
	\stackrel{(54)}{\longrightarrow}
	\quad
	\raisebox{0.2cm}{
	\begin{ytableau}
	3&5&6\\
	1&2&4
	\end{ytableau}}
\]
is an edge in the crystal skeleton $\mathsf{CS}(3,3)$. However, restricting to the interval $[4,6]$ and applying jeu de taquin we obtain the edge
\[
	\raisebox{0.2cm}{
	\begin{ytableau}
	6\\
	4&5
	\end{ytableau}}
	\quad
	\stackrel{(65)}{\longrightarrow}
	\quad
	\raisebox{0.2cm}{
	\begin{ytableau}
	5\\
	4&6
	\end{ytableau}}\, ,
\]
where the cycle edge label has changed from $(54)$ to $(65)$. If we label the edges of the crystal skeleton by the interval of the Dyck pattern
instead of the cycle, then the label of the corresponding edge is $[4,6]$ in both cases. This is one reason we will prefer Dyck pattern labeling in what follows.
\end{remark}

One corollary of Theorem \ref{theorem.self similar} is that one can always find a subgraph of $\CS(\lambda)$ obtained from products of intervals.
\begin{corollary}
\label{corollary.direct product}
Let $[a_k,b_k]\subseteq [1,|\lambda|]$ be intervals for $1\leqslant k \leqslant N$ for some $N>0$ such that $[a_j,b_j] \cap [a_k,b_k] = \emptyset$ 
whenever $j\neq k$. Let $\mu_k=\mathsf{shape}(\mathsf{jdt}(T_{[a_k,b_k]}))$ for all $1\leqslant k \leqslant N$.
Then
\[
	\mathsf{CS}(\mu_1) \times \cdots \times \mathsf{CS}(\mu_N)
\]
is a subgraph of $\mathsf{CS}(\lambda)$ (up to relabelling of the edges).
\end{corollary}

\begin{proof}
This follows directly from Theorem~\ref{theorem.self similar} since the intervals are pairwise disjoint.
\end{proof}

We can also use Theorem \ref{theorem.self similar} to prove the crystal skeleton has branching properties similar to the symmetric group. We first briefly review the branching properties of the 
symmetric group $S_n$. Let $S^\lambda$ be the Specht module, which is the finite-dimensional, irreducible $S_n$-representation indexed by the partition 
$\lambda \vdash n$. Restricting $S^\lambda$ to $S_{n-1}$ gives the branching
\[
	\operatorname{Res}^{S_n}_{S_{n-1}} S^\lambda  = \bigoplus_{\lambda^-} S^{\lambda^-},
\]
where the sum is over all $\lambda^-$ such that $\lambda/\lambda^-$ is a skew shape with a single box. Theorem \ref{theorem.branching} shows that $\CS(\lambda)$ has analogous restriction properties.

\begin{theorem}
\label{theorem.branching}
Let $\lambda \vdash n$ be a partition of $n$. Denote by $\CS(\lambda)_{[1,n-1]}$ the restriction of the crystal skeleton graph $\CS(\lambda)$ 
by removing all edges labeled by $I$ such that $n\in I$. Then we have a graph isomorphism
\[
	\CS(\lambda)_{[1,n-1]} \cong \bigoplus_{\lambda^-} \CS(\lambda^-),
\]
where the sum is over all $\lambda^-$ such that $\lambda/\lambda^-$ is a skew shape with a single box.
\end{theorem}

\begin{proof}
By Theorem~\ref{theorem.self similar}, each $\CS(\lambda^-)$ with $\lambda/\lambda^-$ a skew shape with $|\lambda/\lambda^-|=1$ is a subgraph of 
$\CS(\lambda)_{[1,n-1]}$. Since for $T \in \SYT(\lambda)$ the letter $n$ can sit in any corner cell of $\lambda$, each vertex $T$ in $\CS(\lambda)$ 
corresponds to precisely one vertex in $\bigoplus_{\lambda^-} \CS(\lambda^-)$. Note also that $I \subseteq [1,n-1]$ fixes entry $n$ and hence fixes the corner
of $T$. This proves the claim.
\end{proof}

\begin{corollary}
Let $\lambda \vdash n$ be a partition of $n$. Denote by $\CS(\lambda)_{[2,n]}$ the restriction of the crystal skeleton graph $\CS(\lambda)$ 
by removing all edges labeled by $I$ such that $1\in I$. Then we have a graph isomorphism
\[
	\CS(\lambda)_{[1,n-1]} \cong \CS(\lambda)_{[2,n]}.
\]
\end{corollary}

\begin{proof}
This follows from Theorems~\ref{theorem.branching} and~\ref{theorem.self similar}, which imply that both $\CS(\lambda)_{[1,n-1]}$ and
$\CS(\lambda)_{[2,n]}$ are isomorphic to $\bigoplus_{\lambda^-} \CS(\lambda^-)$ (the second via the jeu de taquin map).
\end{proof}

\subsection{Lusztig involution}\label{section.lusztig}

Recall from Section~\ref{section.crystal} that the crystal $B(\lambda)_n$ of type $A_{n-1}$ enjoys a symmetry under the Sch\"utzenberger or Lusztig 
involution. We show below that this translates into a symmetry of the crystal skeleton as well.

If $T \in \mathsf{CS}(\lambda)$ has descent composition $\Des(T) = (\alpha_1,\ldots,\alpha_\ell)$, then it follows from the definition
of evacuation in Remark~\ref{remark.evacuation} that 
\[ \Des(\mathsf{evac}(T)) = (\alpha_\ell,\ldots,\alpha_1) 
=: \rev(\Des(T)).\] 
Furthermore, if there is an edge $T \stackrel{I}{\longrightarrow} T'$ in $\mathsf{CS}(\lambda)$, then under Lusztig involution
this edge becomes $\mathsf{evac}(T') \stackrel{I^\LL}{\longrightarrow} \mathsf{evac}(T)$, where $I^\LL = [n+1-b,n+1-a]$ if $I=[a,b]$.

\begin{definition}\label{definition.lusztig-on-CS}
Let $\lambda \vdash n$. We define the \defn{Lusztig involution} on the crystal skeleton 
\[
	\LL_n \colon \CS(\lambda) \to \CS(\lambda),
\]
by mapping the vertex $T\in \SYT(\lambda)$ to $\mathsf{evac}(T)$ and an edge $T \stackrel{I}{\longrightarrow} T'$ to
$\mathsf{evac}(T') \stackrel{I^\LL}{\longrightarrow} \mathsf{evac}(T)$, where $I^\LL = [n+1-b,n+1-a]$ if $I=[a,b]$.
\end{definition}

\begin{corollary}
\label{corollary.lusztig}
Let $\lambda \vdash n$. The crystal skeleton $\mathsf{CS}(\lambda)$ is invariant under the Lusztig involution $\LL_n$, that is,
\[
	\CS(\lambda) \cong \LL_n(\CS(\lambda)).
\]
\end{corollary}

\begin{proof}
This follows directly from the fact that the crystal graph $B(\lambda)_n$ underlying $\CS(\lambda)$ is invariant under the Lusztig involution. 
Furthermore, by Proposition~\ref{proposition.quasi edges} an edge $f_i$ in the crystal gives rise to an edge in the crystal skeleton if and only
if there is some $i+1$ that is bracketed with an $i$. Under evacuation, since a word is reversed and complemented,  an $(i+1,i)$  bracketed
pair will transform to an $(n+1-i,n-i)$ pair. Equivalently, whenever an edge $f_i$ is in a quasi-crystal component, it will be mapped to another 
quasi-crystal component under Lusztig involution.
\end{proof}

\subsection{Strongly-connected components}\label{sec:strongly-connected}

We now turn to the strongly-connected components of the crystal skeleton.

\begin{definition}
A directed graph $G=(V,E)$ with vertex set $V$ and edge set $E$ is \defn{strongly-connected} if for any two vertices $u,v\in V$ there is a sequence
of directed edges $e_1,e_2,\ldots, e_k$ with $e_i=(u_{i-1},u_i) \in E$, $u_i\in V$, $u=u_0$ and $v=u_k$.
\end{definition}

\begin{example}
The subgraph
\[
\scalebox{0.7}{
\begin{tikzpicture}[>=latex,line join=bevel,]
\definecolor{permutation}{HTML}{2660A4}
\node (N0) at (0,0) [draw,draw=none] {${\def\lr#1{\multicolumn{1}{|@{\hspace{.6ex}}c@{\hspace{.6ex}}|}{\raisebox{-.3ex}{$#1$}}}
\raisebox{-.6ex}{$\begin{array}[t]{*{3}c}\cline{1-1}\lr{5}\\\cline{1-2}\lr{2}&\lr{4}\\\cline{1-3}\lr{1}&\lr{3}&\lr{6}\\\cline{1-3}\end{array}$}}$};
  \node (N1) at (4,0) [draw,draw=none] {${\def\lr#1{\multicolumn{1}{|@{\hspace{.6ex}}c@{\hspace{.6ex}}|}{\raisebox{-.3ex}{$#1$}}}
  \raisebox{-.6ex}{$\begin{array}[t]{*{3}c}\cline{1-1}\lr{5}\\\cline{1-2}\lr{3}&\lr{4}\\\cline{1-3}\lr{1}&\lr{2}&\lr{6}\\\cline{1-3}\end{array}$}}$};
  
   \draw[->] (N0)  to[bend right=20] node[opLabel]{$[2,4]$}  node[pos=.6, above, permutation]{$(32)$} (N1);
   \draw[->] (N1)  to[bend right=20] node[opLabel]{$[1,3]$}  (N0);
\end{tikzpicture}
}
\]
of the graph in Figure~\ref{figure.CS321} forms a strongly-connected component.
\end{example}

We characterize the crystal skeletons that are strongly-connected.

\begin{theorem}
\label{theorem.rectangle}
A crystal skeleton $\mathsf{CS}(\lambda)$ is strongly-connected if and only if $\lambda=(w^\ell)$ is a rectangle.
\end{theorem}

\begin{proof}
First we prove that $\mathsf{CS}(\lambda)$ is strongly-connected if $\lambda=(w^\ell)$. Let $u$ be the highest weight element in the crystal
$B(\lambda)_n$ for $n\geqslant |\lambda|$. It is given by the Yamanouchi tableau with letters $r$ in row $r$. Let $T_u =\mathsf{std}(u)$, which has 
the letters $(r-1)w+1,(r-1)w+2,\ldots, rw$ in row $r$. Similarly, let $v$ be the lowest weight element in $B(\lambda)_n$. The quasi-crystal $Q_{T_u}$ 
contains both $u$ and $v$ since $\mathsf{std}(u)=\mathsf{std}(v)$. For any $b\in B(\lambda)_n$, there exists a sequence $j_1,j_2,\ldots, j_k$ 
for some $1\leqslant j_i<n$ such that $b=f_{j_1}\cdots f_{j_k} u$. Hence there exists a path from $T_u$ to $T$ in the crystal 
skeleton $\mathsf{CS}(\lambda)$, where $T = \mathsf{std}(b)$.
Conversely, for any $b\in B(\lambda)_n$, there exists a sequence $j_1,j_2,\ldots, j_k$ for some $1\leqslant j_i<n$ such that
$v=f_{j_1}\cdots f_{j_k} b$. Since $u$ and $v$ are in the same quasi-crystal component, this means that there is a path in $\mathsf{CS}(\lambda)$
from $T$ to $T_u$. This proves that $\mathsf{CS}(\lambda)$ is strongly-connected.

Next we show that $\mathsf{CS}(\lambda)$ is not strongly-connected if $\lambda$ is not rectangular. As before, let $u$ be the highest weight element
and $v$ be the lowest weight element in $B(\lambda)_n$. Furthermore, let $T_u = \mathsf{std}(u)$ and $T_v = \mathsf{std}(v)$.
Since $\lambda$ is not rectangular, $T_u\neq T_v$. Namely, $T_u$ has letters $\lambda_1+\cdots+\lambda_{r-1} + 1,\ldots,
\lambda_1+\cdots+\lambda_r$ in row $r$ and $T_v$ has the letters $\lambda_{r+1}+\cdots+\lambda_\ell+1,\ldots,\lambda_r+\cdots+\lambda_\ell$ 
in row $r$ from the top, where $\ell=\ell(\lambda)$. We now demonstrate that it is not possible to find a directed path in $\mathsf{CS}(\lambda)$
from $T_v$ to $T_u$.

Since $\lambda$ is not a rectangle, $\lambda$ has at least two corner cells. Let $\alpha$ (resp. $\beta$) be the topmost (resp. second topmost) 
corner cell in $\lambda$. Let $a_u,b_u$ be the entries in cells $\alpha,\beta$ in $T_u$ and similarly, $a_v,b_v$ be the entries in cells $\alpha,\beta$ 
in $T_v$. We have $a_v<b_v$ and $a_u>b_u$. To find a path from $T_v$ to $T_u$ in $\mathsf{CS}(\lambda)$ we can apply any crystal operator 
$f_i$, but only $e_i$ if there is no bracketing between letters $i+1$ and $i$ (see Proposition~\ref{proposition.quasi edges}). In $\mathsf{row}(T_v)$, 
$a_v$ is to the left of $b_v$ and all letters to the left of $b_v$ are weakly smaller than $b_v$ since $a_v<b_v$. The crystal operators $f_i$ act on the 
rightmost unbracketed letter $i$. Hence applying $f_i$ operators to $T_v$, the entry in cell $\alpha$ is always weakly smaller than the entry in cell 
$\beta$. Similarly, applying $e_i$ only in the case when no $i+1$ is bracketed with $i$, the relative order of the entries in cells $\alpha$ and $\beta$
cannot change. This proves that there is no path from $T_v$ to $T_u$ in $\mathsf{CS}(\lambda)$.
 \end{proof}

\begin{example}
The crystal skeleton $\mathsf{CS}(3,2,1)$ in Figure~\ref{figure.CS321} is not strongly-connected. The crystal skeleton $\mathsf{CS}(3,3)$ in
Figure~\ref{figure.CS33} is strongly-connected.
\begin{figure}
\definecolor{permutation}{HTML}{2660A4}
\begin{tikzpicture}[>=latex,line join=bevel,]
\node (N2) at (0,-3) [draw,draw=none] {${\def\lr#1{\multicolumn{1}{|@{\hspace{.6ex}}c@{\hspace{.6ex}}|}{\raisebox{-.3ex}{$#1$}}}\raisebox{-.6ex}
{$\begin{array}[t]{*{3}c}\cline{1-3}\lr{3}&\lr{5}&\lr{6}\\\cline{1-3}\lr{1}&\lr{2}&\lr{4}\\\cline{1-3}\end{array}$}}$};
  \node (N4) at (0,-6) [draw,draw=none] {${\def\lr#1{\multicolumn{1}{|@{\hspace{.6ex}}c@{\hspace{.6ex}}|}{\raisebox{-.3ex}{$#1$}}}\raisebox{-.6ex}
  {$\begin{array}[t]{*{3}c}\cline{1-3}\lr{2}&\lr{4}&\lr{6}\\\cline{1-3}\lr{1}&\lr{3}&\lr{5}\\\cline{1-3}\end{array}$}}$};
  \node (N1) at (-4,-3) [draw,draw=none] {${\def\lr#1{\multicolumn{1}{|@{\hspace{.6ex}}c@{\hspace{.6ex}}|}{\raisebox{-.3ex}{$#1$}}}\raisebox{-.6ex}
  {$\begin{array}[t]{*{3}c}\cline{1-3}\lr{3}&\lr{4}&\lr{6}\\\cline{1-3}\lr{1}&\lr{2}&\lr{5}\\\cline{1-3}\end{array}$}}$};
  \node (N0) at (0,0) [draw,draw=none] {${\def\lr#1{\multicolumn{1}{|@{\hspace{.6ex}}c@{\hspace{.6ex}}|}{\raisebox{-.3ex}{$#1$}}}\raisebox{-.6ex}
  {$\begin{array}[t]{*{3}c}\cline{1-3}\lr{4}&\lr{5}&\lr{6}\\\cline{1-3}\lr{1}&\lr{2}&\lr{3}\\\cline{1-3}\end{array}$}}$};
  \node (N3) at (4,-3) [draw,draw=none] {${\def\lr#1{\multicolumn{1}{|@{\hspace{.6ex}}c@{\hspace{.6ex}}|}{\raisebox{-.3ex}{$#1$}}}\raisebox{-.6ex}
  {$\begin{array}[t]{*{3}c}\cline{1-3}\lr{2}&\lr{5}&\lr{6}\\\cline{1-3}\lr{1}&\lr{3}&\lr{4}\\\cline{1-3}\end{array}$}}$};
  
  \draw[->] (N0)  to[bend right=20] node[opLabel]{$[1,5]$}  node[pos=.6, above left, permutation]{$(543)$} (N1);
  \draw[->] (N3)  to[bend right=20] node[opLabel]{$[2,6]$}  node[pos=.6, above right, permutation]{$(432)$} (N0);
  \draw[->] (N1)  to node[opLabel]{$[4,6]$}  node[pos=.6, above, permutation]{$(54)$} (N2);
  \draw[->] (N2)  to node[opLabel]{$[1,3]$}  node[pos=.6, above, permutation]{$(32)$} (N3);
  \draw[->] (N0)  to[bend right=30] node[opLabel]{$[2,4]$}  node[pos=.6, right, permutation]{$(43)$} (N2);
  \draw[->] (N2)  to[bend right=30] node[opLabel]{$[3,5]$}  (N0);
  \draw[->] (N1)  to[bend right=18] node[opLabel]{$[1,3]$}  node[pos=.7, above right, permutation]{$(32)$} (N4);
  \draw[->] (N4)  to[bend right=18] node[opLabel]{$[2,4]$}  (N1);
  \draw[->] (N3)  to[bend right=18] node[opLabel]{$[3,5]$}  node[pos=.6, below right, permutation]{$(54)$} (N4);
  \draw[->] (N4)  to[bend right=18] node[opLabel]{$[4,6]$}  (N3);
\end{tikzpicture}
\caption{The strongly-connected crystal skeleton $\mathsf{CS}(3,3)$.
\label{figure.CS33}}
\end{figure}
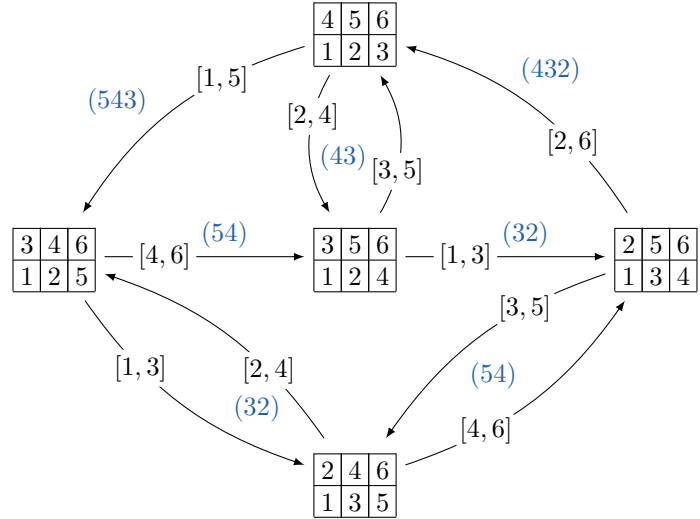
\end{example}

While only crystal skeletons of rectangular shape are strongly-connected, we have the following description of the strongly-connected components 
of $\CS(\lambda)$ in general.

\begin{theorem}\label{theorem.rectangle-covers}
The strongly-connected components of $\mathsf{CS}(\lambda)$ can be covered by unions of direct products
\[
	\mathsf{CS}(w_1^{\ell_1}) \times \cdots \times \mathsf{CS}(w_N^{\ell_N})
\]
of crystal skeletons of rectangular shape.
\end{theorem}

\begin{proof}
By Theorem~\ref{theorem.rectangle}, only crystal skeletons of rectangular shape are strongly-connected. By Corollary~\ref{corollary.direct product},
whenever restricting to disjoint subintervals, by jeu de taquin the subgraphs are isomorphic to products of crystal skeletons. Hence determining the 
strongly-connected component containing a tableau $T \in \mathsf{SYT}(\lambda)$, one can determine all sets of disjoint unions of intervals 
$I_1 \cup \cdots \cup I_k$ of $[|\lambda|]$ such that $\mathsf{jdt}(T_{I_j})$ has rectangular shape $(w_j^{\ell_j})$ for all $1\leqslant j \leqslant k$. 
One then repeats this process for all tableaux in the component isomorphic to $\mathsf{CS}(w_1^{\ell_1}) \times \cdots \times \mathsf{CS}(w_k^{\ell_k})$
and obtains a union of direct products of rectangular crystal skeletons.
\end{proof}

More can be said about the strongly-connected components of $\CS(\lambda)$ when $\lambda$ has at most two parts.
\begin{corollary}
\label{corollary.rectangle-covers}
Let $\lambda$ be a partition with at most two parts. 
The strongly-connected components of $\mathsf{CS}(\lambda)$ are isomorphic to direct products
\[
	\mathsf{CS}(w_1^2) \times \cdots \times \mathsf{CS}(w_N^2)
\]
of crystal skeletons of rectangular shape.
\end{corollary}

\begin{proof}
The crystal skeletons $\mathsf{CS}(\mu)$ for $\mu$ a single column or single row are trivial. Since $\lambda$ has at most two rows,
$\mathsf{shape}(\mathsf{std}(T|_I))$ can have at most two rows.
\end{proof}

A stronger version of Corollary~\ref{corollary.rectangle-covers} is stated in Corollary~\ref{thm:two-row connected components}.

\subsection{Descent compositions}
\label{sec:descent compositions}

We have established that the edges of a crystal skeleton are labeled by odd length intervals coming from Dyck patterns and that 
the vertices are labeled by compositions of $n$ coming from $\Des(T)$ for $T \in \SYT(\lambda)$. 

Our goal is to next describe in Theorem \ref{theorem.descent composition} the relationship between the compositions $\alpha$ and $\beta$ 
in an edge in $\CS(\lambda)$:
\[ 
	(v,\alpha) \xrightarrow{~I~} (w,\beta). 
\]

We first set some notation. Given a composition $\alpha \models n$, write 
\[ 
	\alpha = (\alpha_1, \ldots, \alpha_\ell). 
\]
Sometimes we will identify $\alpha$ with the partition 
$(\alpha^{(1)},\ldots,\alpha^{(\ell)})$ of $[n]$, where 
\begin{equation}
\label{equation.alpha}
	\alpha^{(i)} = \{\alpha_1+\cdots + \alpha_{i-1}+1,\alpha_1+\cdots + \alpha_{i-1}+2,\ldots,\alpha_1+\cdots+\alpha_i\}.
\end{equation}
For instance, if $\alpha = (3,2,2)$, we will also write it as $\alpha = (123|45|67)$.

Since the edges in $\CS(\lambda)$ are labeled by odd length intervals $I = [i,i+2m]$ by Theorem~\ref{theorem.Dyck pattern}, they can be 
written 
\[ I=I^- \cup \{ i+m\} \cup I^+ \quad \textrm{ where } \quad |I^-|=|I^+|>0.\] 
In other words, $I^- =[i,i+m-1]$ and $I^+ = [i+m+1,i+2m]$.
The decomposition of $I$ interacts with the descent composition as follows:

\begin{lemma}
\label{lemma:I.inside.alpha}
Let $T \in \SYT(\lambda)$. Suppose $\alpha = \Des(T)$ and $I = [i,i+2m] = I^- \cup \{i+m\} \cup I^+$ is a Dyck pattern interval in $T$. Then 
$I^- \cup \{i+m\} \subseteq \alpha^{(j)}$ and $I^+ \subseteq \alpha^{(j+1)}$ for some $1\leqslant j < \ell$, 
where $\ell$ is the length of $\alpha$:
\begin{center}
\TIKZ[scale=.75, thick]{
\draw[thin, black!20] (-1.5, 0) to (8.5,0);
\draw[Bracket-Bracket] (0,0) to (3,0);
\draw[Bracket-Bracket] (4,0) to (7,0);
\begin{scope}[every node/.style={draw, thin, inner sep=1.75pt, fill=white, rounded corners}]
\node at (1.5,0) {\small $I^-\vphantom{I_-}$};
\node at (3.5,0) {\small $i+m\vphantom{I_-^-}$};
\node at (5.5,0) {\small $I^+\vphantom{I_-}$};
\end{scope}
\draw[rounded corners=2, sapphire] (-1+.2, .25) to (-1+.2, .5) to 
	node[pos=.15, inner sep=.5pt, fill=white]{$\ \!\cdots$} 
	node[above] {$\alpha^{(j)}$} (4-.15-.05, .5) to (4-.15-.05, .25);
\draw[rounded corners=2, sapphire] (8+.2, .25) to (8+.2, .5) to 
	node[pos=.15, inner sep=.5pt, fill=white]{$\ \!\cdots$} 
	node[above] {$\alpha^{(j+1)}$} (4-.15+.05, .5) to (4-.15+.05, .25);
}
\end{center}
\end{lemma}

\begin{proof}
If $I$ is a Dyck pattern interval on $T$, then (1) the elements of $I^{-} \cup \{ i+m \}$ are in increasing order in $\pi = \row(T)$, and similarly 
for $I^+$, and (2) there is a descent at $i+m$ in $\pi|_I$, i.e. $i+m$ appears after the smallest element of $I^+$ in $\pi$. It follows from (1) that 
$I^- \cup \{ i+m \} \subseteq \alpha^{(j)}$ for some $j$ and $I^+ \subseteq \alpha^{(j')}$ for some $1 \leqslant j, j'\leqslant \ell$, and from 
(2) that $I^+ \not \subseteq \alpha^{(j)}.$ Thus $I^+ \subseteq \alpha^{(j+1)}$, since $I$ is a consecutive interval.
\end{proof}

\begin{example}\label{example.alphacontainsI}
Suppose 
\[ T=
	\raisebox{0.3cm}{
	\begin{ytableau}
	5&6&7&8\\
	1&2&3&4
	\end{ytableau}}, \]
so $\alpha = (1234|5678)$. Then $T$ contains a Dyck pattern interval on $I = [3,5]$, where $I^{-} = \{ 3 \},$ $i+m = 4$ and $I^{+} = \{ 5 \}$. 
We see that $\{ 3, 4 \} \subseteq \alpha^{(1)}$ and $\{ 5 \} \subseteq \alpha^{(2)}$.
\end{example}

Our characterization of how descent compositions change between edges in the crystal skeleton in Theorem~\ref{theorem.descent composition} 
will involve three cases. Recall that for an interval $[a,b] \subseteq [n]$, the restriction $T_{[a,b]}$ is the skew tableau restricted to $[a,b]$ and 
$\jdt(T_{[a,b]})$ is the jeu de taquin straightening of $T_{[a,b]}$ with shape $\shape(\jdt(T_{[a,b]}))$. 

\begin{definition}
Given an interval $[a,b] \subseteq [n]$ of length $2m$, the tableau $T \in \SYT(\lambda)$ has a \defn{rectangle on $[a,b]$} if 
\[ \shape(\jdt(T_{[a,b]})) = (m,m). \]  
Otherwise, we say that $T$ is \defn{rectangle-free on $[a,b]$}.
\end{definition}

\begin{example}\label{example.rectangle}
Consider $T$ as in Example \ref{example.alphacontainsI}. Then $T$ contains a rectangle on $[3,6]$ since 
\[ T_{[3,6]} = 
	\raisebox{0.3cm}{
	\begin{ytableau}
	5&6 & *(gray) & *(gray)\\
	*(gray) & *(gray)&3&4
	\end{ytableau}}  \]
    and 
   \[ \jdt(T_{[3,6]}) = \raisebox{0.3cm}{
	\begin{ytableau}
	5&6 \\
	3&4
	\end{ytableau}}.  \] 
\end{example}

We will be particularly interested in the case, where $[a,b]$ comes from a Dyck pattern interval as in Example~\ref{example.rectangle}.
Recall from Section~\ref{section:cycles} that every Dyck pattern interval on $\pi = \row(T)$ bijects with a cycle $\cycle(\pi|_I)$  and that given $\pi|_I$, 
we have that $\destd(\pi|_I)$ is the destandardization of $\pi|_I$ in the alphabet $\{ i \}^{m+1}, \{ i+1\}^{m}$. The presence or absence of a rectangle is 
closely related to the elements that appear in $\cycle(\pi|_I)$. 

We first prove the following preliminary lemma.

\begin{lemma}
\label{lemma:firstlastbracket}
Given $\pi = \row(T)$ for $T \in \SYT(\lambda)$, let $I=[i,i+2m]$ be a Dyck pattern interval on $\pi$. Then
\begin{enumerate}
    \item 
    \[ \cycle(\pi|_I) = (i+2m, i+2m-1,\ldots, i+m) \] if and only if  the unbracketed $i$ in $\destd(\pi|_I)$ is the rightmost $i$ appearing 
in $\destd(\pi|_I)$, and \medskip
\item 
\[ \cycle(\pi|_I) = (i+m, i+m-1, \ldots, i) \] 
if and only if the unbracketed $i$ in $\destd(\pi|_I)$ is 
the leftmost $i$ in $\destd(\pi|_I)$.
\end{enumerate}
\end{lemma}

\begin{proof}
Suppose there is an edge from $T$ to $T'$ in $\CS(\lambda)$ labeled by the Dyck pattern interval $I$. Let $\pi = \row(T)$ and $\pi' = \row(T')$.
Recall from the discussion in Section \ref{section:cycles} that 
\[ 
	T' = (m+\pi_p,m+\pi_p-1,\ldots,\pi_p) \cdot T, 
\]
where $\pi_p$ is the letter in $T$ (or $\pi$) which corresponds to the $i$ in $\destd(\pi|_I)$ on which $f_i$ acts. It follows that 
$\cycle(\pi|_I) = (i+2m, i+2m-1,\ldots, i+m)$ if and only if $\pi_p = i+m$, which corresponds to $f_i$ acting on the rightmost $i$ in $\destd(\pi|_I)$. 
Similarly, $\cycle(\pi|_I) = (i+m, i+m-1, \ldots, i)$ if and only $\pi_p = i$, which corresponds to $f_i$ acting on the leftmost $i$ in $\destd(\pi|_I)$. 
Since $f_i$ acts on the unique unbracketed $i$ in $\destd(\pi|_I)$, this proves the claim. 
\end{proof}

We will now describe when $T$ has a rectangle coming from a Dyck pattern interval.
\begin{lemma}
\label{lemma:squares.equivalent.shortpattern}
    Suppose $I = [i, i + 2m]$ is a Dyck pattern interval on $T \in \SYT(\lambda)$. Let $\pi = \row(T)$. Then the following is true:
    \begin{enumerate}
            \item $T$ has a rectangle on $[i, i+2m+1]$ if and only if 
        \[ \pi|_{ \{ i+m, \ i+2m, \ i+2m+1 \}}  = i+2m \quad  i+2m+1 \quad i+m  \]
        and $\cycle(\pi|_I) = (i+2m, i+2m-1, \ldots, i+m)$. \medskip
        \item $T$ has a rectangle on $[i-1, i+2m]$ if and only if 
        \[ \pi|_{\{ i-1, \ i, \  i+1 \}} = i \quad i-1 \quad i+1 \]
        and $\cycle(\pi|_I) = (i+m, i+m-1, \ldots, i)$. \medskip
        \item $T$ cannot have a rectangle on both $[i-1,i+2m]$ and $[i, i+2m+1]$.
    \end{enumerate}
\end{lemma}

\begin{proof}
Note that (3) follows immediately from (1) and (2). Since $I$ is a Dyck pattern interval on $T$ and $T$ has a rectangle on both 
$[i-1,i+2m]$ and $[i,i+2m+1]$, then by (1) and (2) we have
\[ 
	(i+2m, i+2m-1, \ldots, i+m) =\cycle(\pi|_I) = (i+m, i+m-1, \ldots, i),
\]
which can occur only if $i+m = i$, a contradiction.

Thus it suffices to prove (1) and (2). We will prove (1); the proof of (2) is analogous and therefore omitted. 

First suppose that $T$ has a rectangle on $J:= [i,i+2m+1]$ where $I = [i, i+2m]$ is a Dyck pattern interval. Note that this is equivalent to the condition that 
\begin{equation}
\label{eq:Pinrectangle} 
	P(\pi|_{J}) =  \scalebox{0.6}{\TIKZ[scale=1.5]{
	\Tableau{{i,i+1,\ldots,i+m-1,i+m},{i+m+1,i+m+2,\ldots,i+2m, i+2m+1}}
	}}\ . 
\end{equation}
Since $i+2m+1$ appears in the top row, it cannot be the last letter in $\pi|_{J}$, and thus $i+m$ must appear to the right of $i+2m+1$. We must 
also have that $i+2m$ appears to the left of $i+2m+1$, otherwise $\shape(\jdt(\pi|_J)) $ would have three rows. Thus we have the following relative 
order in $\pi|_J:$
\begin{equation}
\label{eq:relorder1} 
	i+ 2m \ldots i+2m+1 \ldots i+m.
\end{equation}
Moreover, we must have the relative order
\begin{equation}
\label{eq:relorder2} 
	i+2m \ldots i+m-1 \ldots i+m,
\end{equation}
in $\pi|_I$, since otherwise $\shape(\jdt(\pi|_J))$ would have first row of size larger than $m$. Thus we have that $i+m$ is the last letter in both
$\pi|_I$ and $\pi|_J$. This implies that the unbracketed $i$ in $\destd(\pi|_I)$ is the rightmost $i$. By Lemma~\ref{lemma:firstlastbracket},
this means $\cycle(\pi|_I) = (i+2m, \ldots, i+m)$. 

In the other direction, note that because $I$ is a Dyck pattern interval, 
\[ 
	P(\pi|_{I}) =  \scalebox{0.6}{\TIKZ[scale=1.5]{
	\Tableau{{i,i+1,\ldots,i+m-1,i+m},{i+m+1,i+m+2,\ldots,i+2m}}
	}}\ . 
\]
By assumption, $\pi|_J$ has the same relative order as \eqref{eq:relorder1}. By Lemma~\ref{lemma:firstlastbracket}, if 
$\cycle(\pi|_I) = (i+2m, i+2m-1, \ldots, i+m)$, then the unbracketed $i$ in $\destd(\pi|_I)$ is its rightmost $i$, which implies that $\pi|_J$ 
has the same relative order as \eqref{eq:relorder2}. It follows that $P(\pi|_J)$ has the form in \eqref{eq:Pinrectangle},
and so $T$ has a rectangle on $J = [i, i+2m+1]$. 
\end{proof}

We will use Lemma \ref{lemma:squares.equivalent.shortpattern} to characterize how descent compositions change in $\CS(\lambda)$. 

\begin{theorem}
\label{theorem.descent composition}
Let $T, T' \in \SYT(\lambda)$ with $\Des(T) = \alpha$ and $\Des(T') = \beta$. Suppose there is an edge \[ I = [i,i+2m]= I^- \cup \{ i+m \} \cup I^+ \] from 
$T$ to $T'$ in $\CS(\lambda)$ with $I^- \cup \{ i+m \} \subseteq \alpha^{(j)}$. Then $\beta$ must be of the form   
\[
	\beta = (\alpha^{(1)}, \dots, \alpha^{(j-2)}, 
	\TIKZ[yscale=.45]{\draw[thin] (0,.5) to (0,0) to (2,0) to +(0,.5); \node[circle, draw, inner sep=1pt] at (1,.5) {$*$};}\ ,
	\alpha^{(j+2)}, \dots, \alpha^{(\ell)}),
\]
where \circasterisk\ is determined as follows:
\begin{enumerate}
    \item If $T$ is rectangle-free on $[i-1,i+2m]$ and $[i,i+2m+1]$, then $I$ is a \defn{length-preserving edge}, and 
        \[ \circasterisk = \bigg(\alpha^{(j-1)}, 			\hspace{5pt}
				\alpha^{(j)} \setminus \{i+m\}, 	\hspace{5pt}
				\alpha^{(j+1)} \cup \{i+m\} \bigg). \]
\item If $T$ contains a rectangle on $[i,i+2m+1]$, then $I$ is a \defn{length-increasing edge} and 
\[ \circasterisk = \bigg(\alpha^{(j-1)}, 			\hspace{5pt}
				\alpha^{(j)} \setminus \{ i+m\}, 	\hspace{5pt}
				I^+ \cup \{i+m\}, 				\hspace{5pt}
				\alpha^{(j+1)} \setminus I^+\bigg) .\]
\item If $T$ contains a rectangle on $[i-1, i+2m]$, then $\alpha^{(j)}=I^-\cup \{i+m\}$, and $I$ is a \defn{length-decreasing edge}:
\[ \circasterisk = \bigg(\alpha^{(j-1)} \cup I^-,  		\hspace{5pt}
				\alpha^{(j+1)} \cup \{i+m\} \bigg).\]
\end{enumerate}
\end{theorem}

\begin{proof}
By Lemma \ref{lemma:I.inside.alpha}, we may assume that $I^- \cup \{ i+m \} \subseteq \alpha^{(j)}$ for some $1 \leqslant j \leqslant \ell$. 
Recall that every Dyck pattern interval defines a cycle $\cycle(\pi|_I)$ such that 
\[ 
	\cycle(\pi|_I) \cdot T = T'.
\]

Note that the elements in $\alpha^{(1)}, \ldots, \alpha^{(j-2)}$ and $\alpha^{(j+2)}, \ldots, \alpha^{(\ell)}$ are unchanged by applying $\cycle(\pi|_I)$. Thus 
\[
	\beta = (\alpha^{(1)}, \dots, \alpha^{(j-2)}, 
	\TIKZ[yscale=.45]{\draw[thin] (0,.5) to (0,0) to (2,0) to +(0,.5); \node[circle, draw, inner sep=1pt] at (1,.5) {$*$};}\ ,
	\alpha^{(j+2)}, \dots, \alpha^{(\ell)}),
\]
where it remains to determine the makeup of \circasterisk\ . Furthermore, by definition of the Dyck pattern interval, in $\beta$ we know that $I^-$ 
is in a distinct composition part from $\{ i+m \} \cup I^+$. We will go through each case in turn:
\begin{enumerate}
\item Note that 
      \[ \circasterisk = \bigg(\alpha^{(j-1)}, 			\hspace{5pt}
				\alpha^{(j)} \setminus \{i+m\}, 	\hspace{5pt}
				\alpha^{(j+1)} \cup \{i+m\} \bigg) \]
if the relative positions of $i-1, i$ and $i+2m, i+2m+1$  are preserved under $\cycle(\pi|_I)$. This follows because the presence or absence 
of a descent between $i-1$ and $i$ in $\pi$ will be unchanged after applying $\cycle(\pi|_I)$, and similarly for a descent between $i+2m$ and 
$i+2m+1$. Thus the only descent that is created between $T$ and $T'$ is between $i+m-1$ and $i+m$, which follows from the definition of a 
Dyck pattern interval. 
      
The relative positions of $i-1, i$ and $i+2m, i+2m+1$ are preserved if $i, i+2m \not \in \cycle(\pi|_I)$, since then $\cycle(\pi|_I)$ fixes $i-1,i, i+2m$ 
and $i+2m+1$.  Inspection shows it is also true when $i \in \cycle(\pi|_I)$ but $i-1$ does not appear between $i$ and $i+1$, as well as if 
$i + 2m \in \cycle(\pi|_I)$ but $i+2m+1$ does not occur between $i+2m$ and $i+m+1$. By Lemma~\ref{lemma:squares.equivalent.shortpattern}, 
these conditions are met precisely when $T$ is rectangle-free on $[i-1,i+2m]$ and $[i, i+2m+1]$. \medskip

\item Suppose $T$ contains a rectangle on $[i,i+2m+1]$. Then by Lemma~\ref{lemma:squares.equivalent.shortpattern},  
\[ 
	\pi|_{ \{ i+m, \ i+2m, \ i+2m+1 \}}  = i+2m \quad i+2m+1 \quad i+m  
\]
and $\cycle(\pi|_I) = (i+2m, i+2m-1, \ldots, i+m)$. Note that since $i+2m$ occurs before $i+2m+1$ in $\pi$, we have that 
$\{ i+m, \ldots, i+2m, i+2m+1 \} = I^+ \cup \{ i+2m+1 \} \subseteq \alpha^{(j+1)}.$
Under $\cycle(\pi|_I)$, $i+m+1$ is sent to $i+2m$, and $i+2m+1$ is fixed. Thus in $\cycle(\pi|_I) \cdot \pi$, we have that $i+2m$ 
occurs after $i+2m+1$, and thus these elements occur in different blocks in $\beta$. It follows that
\[ \circasterisk = \bigg(\alpha^{(j-1)}, 			\hspace{5pt}
				\alpha^{(j)} \setminus \{i+m\}, 	\hspace{5pt}
				I^+ \cup \{i+m\}, 				\hspace{5pt}
				\alpha^{(j+1)} \setminus I^+\bigg) .\]
Note that $i+2m+1 \in \alpha^{(j+1)} \setminus I^+$, so this part in $\beta$ is non-empty. \medskip

\item Suppose $T$ contains a rectangle on $[i-1,i+2m]$. Then by Lemma~\ref{lemma:squares.equivalent.shortpattern}, we have 
\[ 
	 \pi|_{\{ i-1, \ i, \  i+1 \}} = i \quad i-1 \quad i+1 
\]
and $\cycle(\pi|_I) = (i+m, i+m-1, \ldots, i)$. Since $I^- \subseteq \alpha^{(j)}$ and $i-1$ appears to the right of $i$ in $\pi$, we must have that 
$i-1 \in \alpha^{(j-1)}$. It follows that $\alpha^{(j)} = I^- \cup \{ i+m \}$. Since $\cycle(\pi|_I)$ maps $i+1$ to $i$, in $\cycle(\pi|_I) \cdot \pi$, 
we have $i-1$ occurring before $i$. Thus in $\beta$, the blocks $\alpha^{(j-1)}$ and $\alpha^{(j)} \setminus \{ i+m \} = I^-$ merge, giving
\[ \circasterisk = \bigg(\alpha^{(j-1)} \cup I^-,  		\hspace{5pt}
				\alpha^{(j+1)} \cup \{i+m\} \bigg).\]
\end{enumerate}
\end{proof}

\begin{example}
\label{ex:compositionsshapes2}
Let 
\[
	T=
	\raisebox{0.3cm}{
	\begin{ytableau}
	6\\
	4&5\\
	1&2&3
	\end{ytableau}} 
	\qquad \text{and} \qquad T'=
	\raisebox{0.3cm}{
	\begin{ytableau}
	6\\
	3&5\\
	1&2&4
	\end{ytableau}}.
\]
Then $\pi = 645123$ and $\pi' = 635124$, and there is an edge from $T$ to $T'$ with Dyck pattern $\pi=423$, so that  $I=[2,4]$ with 
$I^- = \{2 \}$ and $I^+ = \{4 \}$. We have $\alpha = ( 123| 45|6 )$. Note that $T$ contains a rectangle on $[2,5]$. Hence $I$ is a length-increasing 
edge, and $\beta = (12|34|5|6)$.
\end{example}

\begin{example}
Suppose 
\[
	T=
	\raisebox{0.3cm}{
	\begin{ytableau}
	2&5&6\\
	1&3&4
	\end{ytableau}} 
	\qquad \text{and} \qquad T'=
	\raisebox{0.3cm}{
	\begin{ytableau}
	4&5&6\\
	1&2&3
	\end{ytableau}}.
\]
Then $\pi = 256134$ and $\pi'= 456123$, and $T$ contains a Dyck pattern interval $I=[2,6]$ with $I^- = \{ 2,3 \}$, $i+m = 4$, and 
$I^+ = \{ 5, 6 \}$. Note that $T$ contains a rectangle on $[1,6]$. Thus $I$ is a length-decreasing edge, and $\beta =(123|456)$.
\end{example}

By Proposition~\ref{proposition.inside} we know that if $I=[i,i+2m]$ is a Dyck pattern interval of $T$, then so is $J=[i+1,i+2m-1] \subsetneq I$.
We show that $J$ is always length-increasing.

\begin{corollary}
\label{corollary.increase}
Let $T \in \CS(\lambda)$ and $m>1$.
\begin{enumerate}
\item Suppose there is an edge $I=[i,i+2m]$ from $T$ to $T'$ in $\CS(\lambda)$. Then the edge  $J=[i+1,i+2m-1]$ out of $T$ is
length-increasing.
\item Suppose there is an edge $I=[i,i+2m]$ from $T'$ to $T$ in $\CS(\lambda)$. Then the edge $J=[i+1,i+2m-1]$ into $T$ is
length-decreasing.
\end{enumerate}
\end{corollary}

\begin{proof}
We prove part 1 as part 2 is analogous. Let $\pi = \row(T)$. Since $I$ is a Dyck pattern interval of $T$, we have by Definition~\ref{definition.Dyck pattern}
that $P(\pi|_I)$ is as in~\eqref{equation.P pi I}. Removing the letter $i$ from $P(\pi|_I)$ and straightening, gives a tableau of shape $(m,m)$.
Hence by Theorem~\ref{theorem.descent composition}, the interval $J$ is length-increasing.
\end{proof}

\subsection{Fans}
\label{section.fans}
We analyze the local properties of the crystal skeleton near its length increasing and decreasing edges. We call the structures that occur \defn{fans}.

\begin{proposition}
\label{proposition.shorter}
Suppose $T\stackrel{I}{\longrightarrow} T'$ is a length-decreasing edge in $\CS(\lambda)$ (as defined in Theorem~\ref{theorem.descent composition}).
Then one of the two cases holds:
\begin{enumerate}
\item We either have
\begin{center}
\raisebox{-0.8cm}{
\scalebox{0.8}{
\begin{tikzpicture}[>=latex,line join=bevel,]
\node (N0) at (0,0) [draw,draw=none] {$T'$};
\node (N1) at (-2,2) [draw,draw=none] {$T$};
\node (N2) at (2,2) [draw,draw=none] {$T''$};
  \draw[->] (N1)  to node[opLabel]{$I$} (N0);
  \draw[->] (N2)  to node[opLabel]{$I'$} (N0);
  \draw[->] (N2)  to node[opLabel]{$J$} (N1);
\end{tikzpicture}
}}
\end{center}
\[ \textrm{with} \quad |I|>3, \quad \quad  I = [i,i+2m], \quad \quad J = [i-1,i+1], \quad \quad I' = [i+1,i+2m-1] \quad \quad \textrm{ or}\]
\item we have
\begin{center}
\raisebox{-1cm}{
\scalebox{0.8}{
\begin{tikzpicture}[>=latex,line join=bevel,]
\node (N0) at (0,0) [draw,draw=none] {$T'$};
\node (N1) at (0,2) [draw,draw=none] {$T$};
  \draw[->] (N1) to[bend right=30] node[opLabel]{$I$} (N0);
  \draw[->] (N0) to[bend right=30] node[opLabel]{$J$} (N1);
\end{tikzpicture}
}}
\end{center}
\[ \textrm{with} \quad \quad |I| = |J| = 3, \quad \quad I = [i,i+2], \quad \quad J = [i-1,i+1].\]
\end{enumerate}
\end{proposition}

\begin{proof}
Since the edge labeled $I$ is length-decreasing, by Theorem~\ref{theorem.descent composition} and
Lemma~\ref{lemma:squares.equivalent.shortpattern} the letters $i-1,i,i+1$ in the reading word of $T$ appear in the order $\ldots i \ldots i-1 \ldots i+1 \ldots$.
Hence there is an incoming arrow into $T$ labeled $J=[i-1,i+1]$ from $T''$ which is obtained from $T$ by the application of the cycle $(i+1,i)$.
Furthermore, $T'$ is obtained from $T$ by the application of the cycle $(i+m,i+m-1,\ldots,i)$. 

First suppose that $|I|=2m+1>3$. In this case the letter $i+2$ has to be to the right of the letter $i+1$, that is, 
$\ldots i \ldots i-1 \ldots i+1 \ldots i+2 \ldots$. In $T''$ this reads $\ldots i+1 \ldots i-1 \ldots i \ldots i+2 \ldots$ and hence there is a Dyck
pattern comprised of the letters in the interval $I'=[i+1,i+2m-1]$. Combining Lemma~\ref{lemma:squares.equivalent.shortpattern} and
Theorem~\ref{theorem.descent composition}, it follows that 
the new tableau is obtained by the application of the cycle $(i+m,i+m-1,\ldots,i+1)$ on $T''$. Note that this is the same as $T'$, proving the first
claim.

Next assume that $|I|=3$. In this case $T'$ is obtained from $T$ by the arrow indexed by $I=[i,i+2]$ by the application of the cycle $(i+1,i)$.
Recall that the incoming arrow labeled $J=[i-1,i+1]$ also interchanges $i$ and $i+1$. Hence in this case $T''=T'$, which proves the second claim.
\end{proof}

\begin{corollary}
\label{corollary.fans shorter}
Suppose $T\stackrel{I}{\longrightarrow} S$ is a length-decreasing edge in $\CS(\lambda)$ (as defined in Theorem~\ref{theorem.descent composition}).
Then all edges entering $S$ from tableaux with longer descent compositions look as follows (with $I=I_j$ for some $j$):
\begin{center}
\raisebox{-0.8cm}{
\scalebox{0.8}{
\begin{tikzpicture}[>=latex,line join=bevel,]
\node (N0) at (0,0) [draw,draw=none] {$S$};
\node (N1) at (-4,3) [draw,draw=none] {$T_1$};
\node (N2) at (-2,3) [draw,draw=none] {$T_2$};
\node (N3) at (0,3) [draw,draw=none] {$T_3$};
\node (N4) at (2,3) [draw,draw=none] {$T_{m-1}$};
\node (M) at (1,1.5) [draw,draw=none] {$\cdots$};
\node (M1) at (1,3) [draw,draw=none] {$\cdots$};
\node (N5) at (4,3) [draw,draw=none] {$T_m$};
  \draw[->] (N1)  to node[opLabel]{$I_1$} (N0);
  \draw[->] (N2)  to node[opLabel]{$I_2$} (N0);
  \draw[->] (N3)  to node[opLabel]{$I_3$} (N0);
  \draw[->] (N5)  to node[opLabel]{$I_m$} (N0);
  \draw[->] (N2)  to node[opLabel]{$J_1$} (N1);
  \draw[->] (N3)  to node[opLabel]{$J_2$} (N2);
  \draw[->] (N5)  to node[opLabel]{$J_{m-1}$} (N4);
  \draw[->] (N0) to[bend right=30] node[opLabel]{$J_m$} (N5);
\end{tikzpicture}
}}
\end{center}
where 
\[
\begin{aligned}
	I_1&=[i,i+2m], & I_2&=[i+1,i+2m-1], &&\ldots, & I_m&=[i+m-1,i+m+1],\\
	J_1&=[i-1,i+1], & J_2&=[i,i+2], &&\ldots, & J_{m}&=[i+m-2,i+m].
\end{aligned}
\]
\end{corollary}

\begin{example}
In the crystal skeleton $\mathsf{CS}(4,3)$, we have the following local behavior:
\begin{center}
\definecolor{permutation}{HTML}{2660A4}
\scalebox{0.8}{
\def\YScale{1.5}
\TIKZ[yscale=-1, scale=2]{
\node (P134-26) at (1.5,3*\YScale) {
	\TIKZ[scale=.25]{
		\SSYTcolors{{1,1,1,2},{2,3,3}}
		\Tableau{{1,2,3,5},{4,6,7}}
		}
	};
\node (P126-35) at (-1.5,3*\YScale) {
	\TIKZ[scale=.25]{
		\SSYTcolors{{1,1,2,2},{2,3,3}}
		\Tableau{{1,2,4,5},{3,6,7}}
		}
	};
\node (P136-25) at (0,4*\YScale) {
	\TIKZ[scale=.25]{
		\SSYTcolors{{1,1,1,1},{2,2,2}}
		\Tableau{{1,2,3,4},{5,6,7}}
		}
	};
\begin{scope}[-latex, font=\scriptsize, inner sep=2pt]
\draw (P126-35) to node[opLabel]{$[3,7]$} (P136-25);
\draw(P134-26) to node[opLabel]{$[4,6]$}  (P136-25);
\draw(P134-26) to node[opLabel]{$[2,4]$}  (P126-35);
\draw(P136-25) to[bend left=30] node[opLabel]{$[3,5]$}  (P134-26);
\end{scope}
}}
\end{center}
\end{example}
\noindent More examples in the two-row case can be found in Section \ref{two-row: fans}.

There are similar descriptions when $I$ is a length-increasing edge.
\begin{proposition}
\label{proposition.longer}
Suppose $T\stackrel{I}{\longrightarrow} T'$ is a length-increasing edge in $\CS(\lambda)$ (as defined in Theorem~\ref{theorem.descent composition}).
Then one of the two cases holds:
\begin{enumerate}
\item We either have
\begin{center}
\raisebox{-0.8cm}{
\scalebox{0.8}{
\begin{tikzpicture}[>=latex,line join=bevel,]
\node (N0) at (0,0) [draw,draw=none] {$T$};
\node (N1) at (-2,2) [draw,draw=none] {$T'$};
\node (N2) at (2,2) [draw,draw=none] {$T''$};
  \draw[->] (N0)  to node[opLabel]{$I$} (N1);
  \draw[->] (N0)  to node[opLabel]{$I'$} (N2);
  \draw[->] (N1)  to node[opLabel]{$J$} (N2);
\end{tikzpicture}
}}
\end{center}
\[ \textrm{with} \quad \quad |I|>3, \quad \quad I = [i,i+2m], \quad \quad J = [i+2m-1, i+2m+1], \quad \quad I' = [i+1,i+2m-1] \quad \quad \textrm{ or}\]
\item we have
\begin{center}
\raisebox{-1cm}{
\scalebox{0.8}{
\begin{tikzpicture}[>=latex,line join=bevel,]
\node (N0) at (0,0) [draw,draw=none] {$T$};
\node (N1) at (0,2) [draw,draw=none] {$T'$};
  \draw[->] (N0) to[bend right=30] node[opLabel]{$I$} (N1);
  \draw[->] (N1) to[bend right=30] node[opLabel]{$J$} (N0);
\end{tikzpicture}
}}
\end{center}
\[ \textrm{ with } \quad \quad |I| = |J| = 3, \quad \quad I = [i,i+2], \quad \quad J = [i+1,i+3]. \]
\end{enumerate}
\end{proposition}

\begin{proof}
The proof is analogous to the proof of Proposition~\ref{proposition.shorter}.
\end{proof}

\begin{corollary}
\label{corollary.fans longer}
Suppose $T\stackrel{I}{\longrightarrow} S$ is a length-increasing edge in $\CS(\lambda)$ (as defined in Theorem~\ref{theorem.descent composition}).
Then all edges leaving $T$ to tableaux with longer descent compositions look as follows (with $I=I_j$ for some $j$):
\begin{center}
\raisebox{-0.8cm}{
\scalebox{0.8}{
\begin{tikzpicture}[>=latex,line join=bevel,]
\node (N0) at (0,0) [draw,draw=none] {$T$};
\node (N1) at (-4,3) [draw,draw=none] {$S_1$};
\node (N2) at (-2,3) [draw,draw=none] {$S_2$};
\node (N3) at (0,3) [draw,draw=none] {$S_3$};
\node (N4) at (2,3) [draw,draw=none] {$S_{m-1}$};
\node (M) at (1,1.5) [draw,draw=none] {$\cdots$};
\node (M1) at (1,3) [draw,draw=none] {$\cdots$};
\node (N5) at (5,3) [draw,draw=none] {$S_m$};
  \draw[->] (N0)  to node[opLabel]{$I_1$} (N1);
  \draw[->] (N0)  to node[opLabel]{$I_2$} (N2);
  \draw[->] (N0)  to node[opLabel]{$I_3$} (N3);
  \draw[->] (N0)  to node[opLabel]{$I_m$} (N5);
  \draw[->] (N1)  to node[opLabel]{$J_1$} (N2);
  \draw[->] (N2)  to node[opLabel]{$J_2$} (N3);
  \draw[->] (N4)  to node[opLabel]{$J_{m-1}$} (N5);
  \draw[->] (N5) to[bend left=30] node[opLabel]{$J_m$} (N0);
\end{tikzpicture}
}}
\end{center}
where
\[
\begin{aligned}
	I_1&=[i,i+2m], &I_2&=[i+1,i+2m-1], &&\ldots, &I_m&=[i+m-1,i+m+1],\\
	J_1&=[i+2m-1,i+2m+1], &J_2&=[i+2m-2,i+2m], &&\ldots, &J_{m}&=[i+m,i+m+2].
\end{aligned}
\]
\end{corollary}

\subsection{Commutation relations}
\label{section.commutation}

The fans of Section~\ref{section.fans} are part of a larger set of local commutation relations within the crystal
skeleton, which are inherited from local relations in a crystal. Stembridge~\cite{Stembridge.2003} showed that if
$f_i(b) \neq \emptyset$ and $f_j(b) \neq \emptyset$ for $b\in B(\lambda)_n$ with $1\leqslant i<j<n$, then either
\begin{equation}
\label{equation.stembridge square}
	f_i f_j (b) = f_j f_i (b) \quad \text{if $\varphi_j(f_ib) = \varphi_j(b)$ or $\varphi_i(f_jb) = \varphi_i(b)$,}
\end{equation}
or
\begin{equation}
\label{equation.stembridge octagon}
	f_i f_j^2 f_i (b) = f_j f_i^2 f_j (b) \quad \text{if $\varphi_j(f_ib) = \varphi_j(b)+1$ and $\varphi_i(f_jb) = \varphi_i(b)+1$.}
\end{equation}
The case~\eqref{equation.stembridge octagon} can only happen when $j=i+1$.
The relation~\eqref{equation.stembridge square} gives rise to a \defn{commuting square}, 
whereas~\eqref{equation.stembridge octagon} gives rise to a \defn{commuting octagon}.
Similar (dual) relations hold for the raising operators $e_i$ and $e_j$.
See~\cite[Chapter~4]{BumpSchilling.2017} for more details.

To set up the commutation relations in the crystal skeleton, we assume that there is a vertex $T\in \CS(\lambda)$ such that
\[T\stackrel{I}{\longrightarrow} T^I \quad \textrm{ and } \quad T\stackrel{J}{\longrightarrow} T^J \quad \textrm{ with } \min I < \min J.  \]
Note that if $I\neq J$, then $\min I = \min J$ is not possible by the definition of Dyck pattern intervals. Recall that $I$ can either be length-preserving, 
increasing or decreasing according to Theorem~\ref{theorem.descent composition}; we call this the \defn{type} of $I$.

\colorlet{f-one}{sapphire}
\colorlet{f-two}{darkred}
\tikzstyle{preservingEdge}=[thick, |->]
\tikzstyle{pres-decrEdge}=[thick, {Tee Barb[inset=0pt]}->]
\tikzstyle{decrEdge}=[thick, {Arc Barb[reversed]}->]
\tikzstyle{pres-incrEdge}=[thick, {Tee Barb[inset=0pt, reversed]}->]
\tikzstyle{incrEdge}=[thick, {Arc Barb}->]
In what follows, we will describe commutation relations by color coding edges; these colors refer to applications of $f_i$, made specific in the proof of 
Theorem \ref{theorem.commutation}. Edge decorations indicate the following edge types:
\[\TIKZ[scale=.75]{
	\draw[preservingEdge] (0,0) to +(1,0) 
                node[right]{~length-preserving};
	\draw[decrEdge] (0,-1) to +(1,0) 
                node[right]{~length-decreasing};
	\draw[pres-decrEdge] (0,-2) to +(1,0) 
                node[right]{~length-preserving or decreasing};
	\draw[->, thick] (10,0) to +(1,0) 
                node[right]{~any};
	\draw[incrEdge] (10,-1) to +(1,0) 
                node[right]{~length-increasing};
	\draw[pres-incrEdge] (10,-2) to +(1,0) 
                node[right]{~length-preserving or increasing};
}\]

\begin{theorem}
\label{theorem.commutation}
Let $T$ be a vertex in $\CS(\lambda)$ such that $T\stackrel{I}{\longrightarrow} T^I$ and $T\stackrel{J}{\longrightarrow} T^J$ with
$I=[i,i+2m]$ and $J=[j,j+2\ell]$ with $i<j$. Then we have the following local commutation relations:
\begin{enumerate}[label=\rm {Case }\arabic{enumi}.,ref=\arabic{enumi}]
\item \label{case1}
 \TIKZ[scale=.75, thick, every node/.style={draw, thin, inner sep=2pt, fill=white, rounded corners}]{
\draw[thin, black!20] (-.5, 0) to (6.5,0);
\draw[intEdge] (0,0) to node{\small $I$} +(2,0);
\draw[intEdge] (3,0) to node{\small $J$} +(3,0);
}

\noindent If $I\cap J = \emptyset$, then we have a \defn{square} 
\[\TIKZ[line join=bevel, scale=1.25]{
\node (N0) at (0,0)  {$T$};
\node (N1) at (-1,-1) {$T^I$};
\node (N2) at (1,-1)  {$T^J$};
\node (N3) at (0,-2) {$S$};
  \draw[->,f-one,thick] (N0)  to node[edgeLabel]{$I$} (N1);
  \draw[->,f-two,thick] (N0)  to node[edgeLabel]{$J$} (N2);
  \draw[->,f-two,thick] (N1)  to node[edgeLabel]{$J'$} (N3);
  \draw[->,f-one,thick] (N2) to node[edgeLabel]{$I'$} (N3);
}\]
with $I'=I$ and $J'=J$, such that: \smallskip
\begin{enumerate}
    \item \label{case1a}
    The edges labeled $I$ and $I'$ (resp. $J$ and $J'$) are of same type if 
    \begin{itemize}
    \item $\max I + 1 \neq \min J \quad \quad$ or
    \item $\max I + 1 = \min J$ and there is no edge labeled $[i,i+2m+2]$ into $T$ or no edge labeled $[j-2,j+2\ell]$ out of $S$.
    \end{itemize}
\item \label{case1b}
If $\max I + 1 = \min J$, there is an edge labeled $[i,i+2m+2]$ into $T$ and an edge labeled $[j-2,j+2\ell]$ out of $S$,
edge $I$ is length-increasing, edge $J'$ is length-decreasing, and the edges $J$ and $I'$ are length-preserving.
\end{enumerate} \medskip

\item \label{case2} \TIKZ[scale=.75, thick, every node/.style={draw, thin, inner sep=3pt, fill=white, rounded corners}]{
\draw[intEdge] (0,0) to node{\small $I$} +(2,0);
\draw[intEdge] (1.5,-.15) to node{\small $J$} +(3,0);
}

\noindent If $I \cap J \neq \emptyset$, $J \not \subseteq I$, and $I$ is not length-increasing, consider $T^I\stackrel{J'}{\longrightarrow} S^{J'}$
and $T^J\stackrel{I'}{\longrightarrow} S^{I'}$. One of the following must hold: \smallskip
\begin{enumerate}
\item \label{case2a}
We have $J \subsetneq J'$ and one of the following cases:\smallskip
\begin{enumerate}
\item \label{case2ai}
If $|I|>3$, we have a \defn{square} 
\[\TIKZ[line join=bevel, scale=1.25]{
\node (N0) at (0,0) {$T$};
\node (N1) at (-1,-1) {$T^I$};
\node (N2) at (1,-1) {$T^J$};
\node (N3) at (0,-2) {$S$};
  \draw[preservingEdge,f-one] (N0)  to node[edgeLabel]{$I$} (N1);
  \draw[preservingEdge,f-two] (N0)  to node[edgeLabel]{$J$} (N2);
  \draw[preservingEdge,f-two] (N1)  to node[edgeLabel]{$J'$} (N3);
  \draw[preservingEdge,f-one] (N2) to node[edgeLabel]{$I'$} (N3);
}
\qquad
\TIKZ[scale=.75, thick, every node/.style={draw, thin, inner sep=3pt, fill=white, rounded corners}]{
\coordinate (I1) at (0,0); \coordinate (I2) at (3,0); 
\coordinate (I'1) at (.25,.5); \coordinate (I'2) at (.25+2.5,.5); 
\coordinate (J1) at (2,-.25); \coordinate (J2) at (2+4,-.25); 
\coordinate (J'1) at (2-.25,-.75); \coordinate (J'2) at (2-.25+4.5,-.75); 
\begin{scope}[densely dotted, thin] 
	\draw (I'1) to +(0, -.5);
	\draw (I'2) to +(0, -.5);
	\draw (J1) to +(0, -.5);
	\draw (J2) to +(0, -.5);
\end{scope}
\draw[intEdge, f-one] (I'1) to node[pos=.6]{\small $I'$} (I'2);
\draw[intEdge, f-one] (I1) to node[pos=.3]{\small $I$} (I2);
	\foreach \x in {.25, 2.75}{\draw[f-one, thin] (\x,-.1) to +(0,.2);}
\draw[intEdge, f-two] (J1) to node[pos=.6]{\small $J$} (J2);
\draw[intEdge, f-two] (J'1) to node[pos=.4]{\small $J'$} (J'2);
	\foreach \x in {2, 6}{\draw[f-two, thin] (\x,-.75-.1) to +(0,.2);}
}\]
with $J'=[j-1,j+2\ell+1]$ and $I'=[i+1,i+2m-1]$, where all edges are length-preserving. \smallskip
\item \label{case2aii}
If $|I|=3$, we have a \defn{triangle}
\[
\TIKZ{
\node (N0) at (0,0) {$T$};
\node (N1) at (-1,-1.5) {$T^I$};
\node (N2) at (1,-1.5) {$T^J$};
  \draw[preservingEdge,f-one] (N0)  to node[edgeLabel]{$I$} (N1);
  \draw[decrEdge,f-two] (N0)  to node[edgeLabel]{$J$} (N2);
  \draw[decrEdge,f-two] (N1)  to node[edgeLabel]{$J'$} (N2);}
\quad\text{with}\quad
\begin{array}{r@{\ }l}
I&=[i,i+2]\\
J&=[i+2,i+2+2\ell]\\
J'&=[i+1,i+3+2\ell]
\end{array}
\quad
\TIKZ[scale=.75, thick, every node/.style={draw, thin, inner sep=3pt, fill=white, rounded corners}]{
\node[draw=none] at (0,-1.75){};
\coordinate (I1) at (0,0); \coordinate (I2) at (3*.25,0); 
\coordinate (J1) at (.5,-.5); \coordinate (J2) at (.5+3,-.5); 
\coordinate (J'1) at (.25,-1); \coordinate (J'2) at (.25+3+.5,-1); 
\begin{scope}[densely dotted, thin] 
	\draw (J1) to +(0, .5);
	\draw (J1) to +(0, -.5);
	\draw (J2) to +(0,- .5); 
	\draw (J'1) to +(0,1);
	\draw (I2) to +(0,-1);
\end{scope}
\draw[intEdge, f-one] (I1) to node[above, draw=none, outer sep=1pt] {\small $I\strut$} (I2); 
	\foreach \x in {.25, .5}{\draw[f-one, thin] (\x,-.1) to +(0,.2);}
\draw[intEdge, f-two] (J1) to node[pos=.6]{\small $J$} (J2);
\draw[intEdge, f-two] (J'1) to node[pos=.3]{\small $J'$} (J'2);
	\foreach \x/\y in {.5/-1, .75/-1, 3.5/-1, .75/-.5}{\draw[f-two, thin] (\x,\y-.1) to +(0,.2);}
}
\]
where $I$ is length-preserving and $J$ and $J'$ are length-decreasing. 
\end{enumerate} \smallskip
\item \label{case2b}
We have $J=J'$ and one of the following cases must hold:
\begin{enumerate}
\item \label{case2bi}
We have $I' \subsetneq I$ and an \defn{octagon}
\[\TIKZ[scale=1.2]{
\node (N0) at (0,0) {$T$};
\node (N1) at (-1.5,-1) {$T^I$};
\node (N2) at (1.5,-1) {$T^J$};
\node (N3) at (-1.5,-2.5) {$S^{J'}$};
\node (N4) at (1.5,-2.5) {$S^{I'}$};
\node (N5) at (-1.5,-4) {$\bullet$};
\node (N6) at (1.5,-4) {$\bullet$};
\node (N7) at (0,-5) {$\bullet$};
  \draw[pres-decrEdge,f-one] (N0)  to node[edgeLabel]{$I$} (N1);
  \draw[preservingEdge, f-two] (N0)  to node[edgeLabel]{$J$} (N2);
  \draw[pres-incrEdge,f-two] (N1)  to node[edgeLabel]{$J'$} (N3);
  \draw[preservingEdge,f-one] (N2) to node[edgeLabel]{$I'$} (N4);
  \draw[preservingEdge,f-two] (N3) to node[edgeLabel]{$J''$} (N5);
  \draw[pres-decrEdge,f-one] (N4) to node[edgeLabel]{$I''$} (N6);
  \draw[preservingEdge,f-one] (N5) to node[edgeLabel]{$I'''$} (N7);
  \draw[pres-incrEdge,f-two] (N6) to node[edgeLabel]{$J'''$} (N7);
}\quad \text{with} \quad
\begin{matrix}
\TIKZ[scale=.75, thick, every node/.style={draw, thin, inner sep=3pt, fill=white, rounded corners}]{
\coordinate (I1) at (0,0); \coordinate (I2) at (3,0); 
\coordinate (I'1) at (.25,.5); \coordinate (I'2) at (.25+2.5,.5); 
\coordinate (I''1) at (0,1); \coordinate (I''2) at (2.5,1); 
\coordinate (J1) at (2,-.25); \coordinate (J2) at (2+4,-.25); 
\coordinate (J''1) at (2-.25,-.75); \coordinate (J''2) at (2-.25+4,-.75); 
\coordinate (J'''1) at (2-.5,-1.25); \coordinate (J'''2) at (2-.5+4.5,-1.25); 
\begin{scope}[densely dotted, thin] 
	\draw (I'1) to +(0, -.5); 
	\draw (I'1) to +(0, .5); 
	\draw (I'2) to +(0, -.5); 
	\draw (I''1) to +(0, -1); 
	\draw (I''2) to +(0, -1); 
	\draw (J''1) to +(0, -.5); 
	\draw (J''2) to +(0, -.5); 
	\draw (J''2) to +(0, .5); 
	\draw (J1) to +(0, -1);  
	\draw (J'''2) to +(0, 1);
\end{scope}
\draw[intEdge, f-one] (I1) to node[pos=.3]{\small $I$} (I2);
\draw[intEdge, f-one] (I'1) to node[pos=.5]{\small $I'$} (I'2);
\draw[intEdge, f-one] (I''1) to node[pos=.3]{\small $I''$} (I''2);
	\foreach \x/\y in {.25/0, 2.75/0, 2.5/.5, .25/1, 2.5/0}{\draw[f-one, thin] (\x,\y-.1) to +(0,.2);}
\draw[intEdge, f-two] (J1) to node[pos=.7]{\small $J$} (J2);
\draw[intEdge, f-two] (J''1) to node[pos=.5]{\small $J''$} (J''2);
\draw[intEdge, f-two] (J'''1) to node[pos=.3]{\small $J'''$} (J'''2);
	\foreach \x/\y in {2/-1.25, 1.75/-1.25, 2/-.75, 5.75/-.25, 5.75/-1.25}{\draw[f-two, thin] (\x,\y-.1) to +(0,.2);}
}\\~\\
\begin{array}{r@{\ }l}
I&=[i, i+2m]\\
I'&=[i+1,i+2m-1]\\
I''=I'''&=[i,i+2m-2]\\
J=J'&=[j, j+2\ell]\\
J''&=[j-1,j+2\ell-1]\\
J'''&=[j-2,j+2\ell]\\
\end{array}
\end{matrix}
\]
where edges $I$ and $I''$ are length-decreasing or preserving, edges $J'$ and $J'''$ are length-increasing or preserving, and all other 
edges are length-preserving. When $J'$ is length-increasing, there is also an edge $[j-1,j+2\ell+1]$ out of $T^I$ and case~\ref{case2ai} applies
as well.
\smallskip
\item \label{case2bii}
We have $I \subsetneq I'$ with $|I|=3$ and a \defn{pentagon}
\[\TIKZ[scale=1.5]{
\node (N0) at (0,0) {$T$};
\node (N1) at (-1,-.8) {$T^I$};
\node (N2) at (1,-1.5) {$T^J$};
\node (N3) at (-1,-2.2) {$S^{J'}$};
\node (N5) at (0,-3) {$\bullet$};
  \draw[pres-decrEdge,f-one,thick] (N0)  to node[edgeLabel]{$I$} (N1);
  \draw[decrEdge,f-two] (N0)  to node[edgeLabel]{$J$} (N2);
  \draw[->,f-two,thick] (N1)  to node[edgeLabel]{$J'$} (N3);
  \draw[->,thick,f-two] (N2) to node[edgeLabel]{$I'$} (N5);
  \draw[decrEdge,f-two] (N3) to node[edgeLabel]{$J''$} (N5);
}
\quad \text{with} \quad
\begin{matrix}
\TIKZ[scale=.75, thick, every node/.style={draw, thin, inner sep=3pt, fill=white, rounded corners}]{
\coordinate (I1) at (0,0); \coordinate (I2) at (3*.25,0); 
\coordinate (I'1) at (0,.5); \coordinate (I'2) at (.75+3,.5); 
\coordinate (J1) at (.5,-.5); \coordinate (J2) at (.25+.5+3,-.5); 
\coordinate (J'1) at (.25,-1); \coordinate (J'2) at (3+.5,-1); 
\begin{scope}[densely dotted, thin] 
	\draw (I2) to +(0,-1);
	\draw (I2) to +(0,.5);
	\draw (I'1) to +(0, -.5);
	\draw (J1) to +(0, 1); 
	\draw (J1) to +(0, -.5); 
	\draw (J2) to +(0,1);
	\draw (J'2) to +(0, 1.5); 
	\draw (J'1) to +(0,1.5);
\end{scope}
\draw[intEdge, f-one] (I1)  node[left, fill=none, draw=none] {\small $I\strut$}to (I2); 
\draw[intEdge, f-two] (I'1) to node{$I'$} (I'2);
	\foreach \x in {.25, .5}{\draw[f-one, thin] (\x,-.1) to +(0,.2);}
\draw[intEdge, f-two] (J1) to node[pos=.6]{\small $J$} (J2);
\draw[intEdge, f-two] (J'1) to node[pos=.4]{\small $J''$} (J'2);
	\foreach \x/\y in {.5/-1, .75/-.5,.75/-1, 3.5/-.5, 3.5/.5, .25/.5, .5/.5, .75/.5}{\draw[f-two, thin] (\x,\y-.1) to +(0,.2);}
}\\~\\
\begin{array}{r@{\ }l}
I &= [i, i+2]\\
I' &= [i, i+2+2\ell]\\
J =J' &= [i+2, i+2+2\ell]\\
J''&= [i+1, i+1+2\ell]
\end{array}
\end{matrix}
\]
where $J$ and $J''$ are length-decreasing.
\smallskip
\item \label{case2biii}
We have $I'=I$ and a \defn{square} with edges labeled $I$ (resp. $J$) of the same type and $J$ is not length-decreasing
\[\TIKZ[line join=bevel, scale=1.25]{
\node (N0) at (0,0)  {$T$};
\node (N1) at (-1,-1) {$T^I$};
\node (N2) at (1,-1)  {$T^J$};
\node (N3) at (0,-2) {$S$};
  \draw[pres-decrEdge,f-one,thick] (N0)  to node[edgeLabel]{$I$} (N1);
  \draw[pres-incrEdge,f-two,thick] (N0)  to node[edgeLabel]{$J$} (N2);
  \draw[pres-incrEdge,f-two,thick] (N1)  to node[edgeLabel]{$J$} (N3);
  \draw[pres-decrEdge,f-one,thick] (N2) to node[edgeLabel]{$I$} (N3);
}\]
\end{enumerate}
\end{enumerate}
\item \label{case3} \TIKZ[scale=.75, thick, every node/.style={draw, thin, inner sep=2pt, fill=white, rounded corners}]{
\draw[intEdge] (0,0) to node[pos=.3]{\small $I$} +(4,0);
\draw[intEdge] (.5,-.5) to node[pos=.6]{\small $J$} +(3,0);
}

\smallskip

\noindent 
If $J \subsetneq I=[i,i+2m]$ with $|J|=|I|-2$, we have the following: \smallskip
\begin{enumerate}
\item \label{case3a} If $I$ is length-increasing, we have a \defn{triangle} 
\[
\TIKZ{
\node (N0) at (0,0) {$T$};
\node (N1) at (-1,-1.5) {$T^I$};
\node (N2) at (1,-1.5) {$T^J$};
  \draw[incrEdge,f-one] (N0)  to node[edgeLabel]{$I$} (N1);
  \draw[incrEdge,f-two] (N0)  to node[edgeLabel]{$J$} (N2);
  \draw[preservingEdge,f-two] (N1)  to node[edgeLabel]{$J'$} (N2);}
\quad\text{with}\quad
\begin{array}{r@{\ }l}
I&=[i,i+2m]\\
J&=[i+1,i+2m-1]\\
J'&=[i+2m-1,i+2m+1]
\end{array}
\quad
\TIKZ[scale=.75, thick, every node/.style={draw, thin, inner sep=3pt, fill=white, rounded corners}]{
\node[draw=none] at (0,.75){};
\coordinate (I1) at (0,0); \coordinate (I2) at (3,0); 
\coordinate (J1) at (.25,-.5); \coordinate (J2) at (3-.25,-.5); 
\coordinate (J'1) at (3-.5,-1); \coordinate (J'2) at (3+.25,-1); 
\begin{scope}[densely dotted, thin] 
	\draw (J1) to +(0, .5);
	\draw (J2) to +(0,- .5);
	\draw (J2) to +(0,.5); 
	\draw (J'1) to +(0,1);
	\draw (I2) to +(0,-1);
\end{scope}
\draw[intEdge, f-one] (I1) to node[pos=.35]{\small$I$} (I2); 
\draw[intEdge, f-two] (J1) to node[pos=.6]{\small $J$} (J2);
\draw[intEdge, f-two] (J'1) to node[below, draw=none, outer sep=1pt] {\small $J'\strut$} (J'2);
	\foreach \x in {.25, 2.5, 2.75}{\draw[f-one, thin] (\x,-.1) to +(0,.2);}
	\foreach \x in {3-.25, 3}{\draw[f-two, thin] (\x,-1-.1) to +(0,.2);}
	\draw[f-two, thin] (2.5,-.5-.1) to +(0,.2);
}
\]
where $J$ is length-increasing and $J'$ is length-preserving. \smallskip
\item \label{case3b}
We have a \defn{pentagon}
\[\TIKZ[scale=1.5, yscale=.9]{
\node (N0) at (0,0) {$T$};
\node (N1) at (-1,-1.5) {$T^I$};
\node (N2) at (1,-.8) {$T^J$};
\node (N4) at (1,-2.2) {$S^{I'}$};
\node (N5) at (0,-3) {$\bullet$};
  \draw[->,thick,f-one] (N0)  to node[edgeLabel]{$I$} (N1);
  \draw[incrEdge,f-two] (N0)  to node[edgeLabel]{$J$} (N2);
  \draw[incrEdge,f-two] (N1)  to node[edgeLabel]{$J'$} (N5);
  \draw[->,thick,f-one] (N2) to node[edgeLabel]{$I'$} (N4);
  \draw[pres-incrEdge,f-two] (N4) to node[edgeLabel]{$J''$} (N5);
}
\quad \text{with} \quad 
\begin{matrix}
\TIKZ[scale=.75, thick, every node/.style={draw, thin, inner sep=3pt, fill=white, rounded corners}]{
\coordinate (I1) at (0,0); \coordinate (I2) at (3,0); 
\coordinate (J1) at (.25,-.5); \coordinate (J2) at (3-.25,-.5); 
\coordinate (J'1) at (0,-1.25); \coordinate (J'2) at (2.5,-1.25); 
\coordinate (J''1) at (3-.75,-1.75); \coordinate (J''2) at (3,-1.75); 
\begin{scope}[densely dotted, thin] 
	\draw (I1) to +(0,-1.25);
	\draw (J1) to +(0, .5);
	\draw (J1) to +(0, -.75); 
	\draw (J2) to +(0,-1.25); 
	\draw (J2) to +(0,.5);
	\draw (J''1) to +(0,1.75); 
	\draw (3-.5,-1.75) to +(0,1.75);
	\draw (I2) to +(0,-1.75);
\end{scope}
\draw[intEdge, f-one] (I1) to node[pos=.35]{\small$I$} (I2); 
\draw[intEdge, f-two] (J1) to node[pos=.6]{\small $J$} (J2);
\draw[intEdge] (J'1) to node[pos=.5]{\small {\color{f-one}$I'$}, {\color{f-two}$J'$}} (J'2);
\draw[intEdge, f-two] (J''1) to node[below, draw=none, fill=none, outer sep=1pt] {\small $J''\strut$} (J''2);
	\foreach \x in {.25, 2.25, 2.5, 2.75}{\draw[f-one, thin] (\x,-.1) to +(0,.2);}
	\foreach \x in {3-.5, 3-.25}{\draw[f-two, thin] (\x,-1.75-.1) to +(0,.2);}
	\foreach \x in {.25, 2.25}{\draw[thin] (\x,-1.25-.1) to +(0,.2);}
	\foreach \x in {2.25, 2.5}{\draw[f-two, thin] (\x,-.5-.1) to +(0,.2);}
}
~\\
\begin{array}{r@{\ }l}
I&=[i,i+2m]\\
J&=[i+1,i+2m-1]\\
I' = J'&=[i,i+2m-2]\\
J''&= [i+2m-2,i+2m]
\end{array}\end{matrix}\]
where $J$ and $J'$ are length-increasing and $J''$ is not length-decreasing.
\end{enumerate}
\end{enumerate}
\end{theorem}

\begin{proof}
Recall from Remark~\ref{remark.Dyck pattern} that an edge labeled by the Dyck pattern interval $I=[i,i+2m]$ (resp. $J=[j,j+2\ell]$) in the crystal skeleton 
$\CS(\lambda)$ corresponds to an edge $f_i$ (resp. $f_j$) in the crystal $B(\lambda)$ by destandardizing the letters in $I$ (resp. $J$) to letters in 
$\{i,i+1\}$ (resp. $\{j,j+1\}$). 

In the crystal $f_if_j = f_j f_i$ if $i+1<j$. In this case $I \cap J = \emptyset$ with $\max I +1 \neq \min J$ and we obtain a square as in case~\ref{case1}.
By Theorem~\ref{theorem.descent composition}, the type of the edges is determined by the presence or absence of certain rectangles. These
do not change when $I\cap J=\emptyset$ unless $\max I +1 = \min J$, hence the edges $I$ (resp. $J$) have the same type.

We can have $I \cap J = \emptyset$ also when $j=i+1$. By self-similarity (see Section~\ref{section.self similar}), we may assume without 
loss of generality that $I=[1,2m+1]$, that is $i=1$. If $I$ is not length-increasing, we have
\begin{equation}
\label{equation.t jdt}
	t:=\jdt(\destd(T_{I\cup J})) =
	\TIKZ[scale=.35]{
	\fill[ssyt2] (0,0) to (5,0) to (5,1) to (4,1) to (4,2) to (0,2) to (0,0);
	\Tableau{{1,1,1,1,1,2,2,2,2},{2,2,2,2,3,3,3},{3,3,3}}
	\draw[thick, saffron] (0,0) to (5,0) to (5,1) to (4,1) to (4,2) to (0,2) to (0,0);
	\draw[|-|] (0,-.5) to node[below]{$m$} +(4,0);
	\draw[|-|] (5,-.5) to node[below]{$a$} +(4,0);
	\draw[|-|] (0,3.5) to node[above]{$b$} +(3,0);
	\draw[|-|] (4,2.5) to node[above]{$c$} +(3,0);
	},
\end{equation}
where $0\leqslant a$, $0\leqslant b \leqslant m$ and $0\leqslant c \leqslant a+1$. 
The boxes shaded in yellow correspond to the Dyck pattern interval $I$.
In this setup, the edge $I$ (resp. $J$) in the crystal skeleton corresponds to $f_1t$ (resp. $f_2t$).
In the commutation relations, we have colored the edges according to the  corresponding $f_i$:
\[\TIKZ{
\draw[->, thick, f-one] (0,0) to 
    node[below]{(blue)} 
    node[above]{$f_1$} +(2,0);
\draw[->, thick, f-two] (3,0) to 
    node[below]{(red)} 
    node[above]{$f_2$} +(2,0);
}\]
Since $I \cap J = \emptyset$, we have $b+c < a$.
The crystal operator $f_1$ acts on the rightmost 1 and $f_2$ acts on the rightmost 2, so that $f_1 f_2 = f_2 f_1$. Furthermore the types of the edges $I$
(resp. $J$) are the same. Hence case~\ref{case1} holds. 

If $I$ is length-increasing, we have
\begin{equation}
\label{equation.t jdt1}
	t:=\jdt(\destd(T_{I\cup J})) =
	\TIKZ[scale=.35]{
	\fill[ssyt2] (0,0) to (5,0) to (5,1) to (4,1) to (4,2) to (0,2) to (0,0);
	\Tableau{{1,1,1,1,1,2,2,2,2},{2,2,2,2,2,3,3,3},{3,3,3}}
	\draw[thick, saffron] (0,0) to (5,0) to (5,1) to (4,1) to (4,2) to (0,2) to (0,0);
	\draw[|-|] (0,-.5) to node[below]{$m$} +(4,0);
	\draw[|-|] (5,-.5) to node[below]{$a$} +(4,0);
	\draw[|-|] (0,3.5) to node[above]{$b$} +(3,0);
	\draw[|-|] (5,2.5) to node[above]{$c$} +(3,0);
	},
\end{equation}
where $0\leqslant a$, $0\leqslant b \leqslant m+1$ and $0\leqslant c \leqslant a$. 
Since $I \cap J = \emptyset$, we have $b+c \leqslant a$.
When $b>0$ or when $b+c<a$, it can be checked explicitly that the square holds and that the edges $I$ (resp. $J$) have the same type.
When $b=0$ and $c=a$, we can have $\max I +1 = \min J$. In this case, there is an incoming edge $[i,i+2m+2]$ into $T$ and an outgoing edge
$[j-2,j+2\ell]$ out of $S$. Again, it can be checked explicitly that case~\ref{case1} holds.

Let us now assume that the conditions of case~\ref{case2} hold. Again by self-similarity (see Section~\ref{section.self similar}), we may assume without 
loss of generality that $I=[1,2m+1]$, that is $i=1$ and that $t$ is as in~\eqref{equation.t jdt}.
Since $I \cap J \neq \emptyset$, we have $b+c \geqslant a$.

If $c=a+1$, $b<m$, and $|I|>3$, it can be checked explicitly using the combinatorial definitions of $f_1$ and $f_2$ that case~\ref{case2ai} holds. 
If $c=a+1$, $b<m$, and $|I|=3$, the edge $f_2t \stackrel{f_1}{\longrightarrow} f_1 f_2 t$ is contracted in the crystal skeleton by 
Proposition~\ref{proposition.quasi edges} since there is no bracketed pair for the letters $\{1,2\}$. Edge $J$ is length-decreasing by 
Corollary~\ref{corollary.increase}. Hence case~\ref{case2aii} holds.

If $c=a$, $b<m$, and $|I|>3$, it can be checked explicitly that the octagon in case~\ref{case2bi} holds.
Note that $J'$ is length-increasing if there is a 4 next to the rightmost 3 in~\eqref{equation.t jdt}. This is in the same quasi-crystal component
as the same tableau with $c=a+1$, so that case~\ref{case2ai} applies.
If $c=a$, $b<m$, and $|I|=3$, the octagon relation $f_1 f_2^2 f_2 t = f_2 f_1^2 f_2 t$ holds in the crystal, but the edges
\[ 
	f_2t \stackrel{f_1}{\longrightarrow} f_1 f_2 t, \quad \quad f_1 f_2t \stackrel{f_1}{\longrightarrow} f_1^2 f_2 t,  \quad \textrm{ and } 
	\quad f_2^2 f_1t \stackrel{f_1}{\longrightarrow} f_1 f_2^2 f_1 t\]
are contracted (i.e. in the same quasi-crystal) since there are no bracketed pairs in the letters $\{1,2\}$.
Edge $J''$ is length-decreasing by Corollary~\ref{corollary.increase}. Since $b=0$ in this case, edge $J$ is length-decreasing.
Hence the pentagon of case~\ref{case2bii} holds.

Finally, if $c<a$ or if $b=m$, it can be checked explicitly that case~\ref{case2biii} holds.

Case~\ref{case3a} follows from Proposition~\ref{proposition.longer}. Alternatively, it can be derived from the (contracted) square
relation in the crystal on the tableau
\[t:=\jdt(\destd(T_{I}))  =
	\TIKZ[scale=.35]{
        \node[above] at (1,2.5){$\phantom{m}$};
	\fill[ssyt2] (0,0) to (5,0) to (5,1) to (4,1) to (4,2) to (0,2) to (0,0);
	\Tableau{{1,1,1,1,1},{2,2,2,2,3}}
	\draw[thick, saffron] (0,0) to (5,0) to (5,1) to (4,1) to (4,2) to (0,2) to (0,0);
	\draw[|-|] (0,-.5) to node[below]{$m$} +(4,0);
	}.
\]

Case~\ref{case3b} can be derived from the (contracted) octagon
relation in the crystal on the previous tableau with $3$ replaced by $4$ when $I$ is length-increasing or 
on the following tableau when $I$ is length-preserving or length-decreasing
\[t:=\jdt(\destd(T_{I}))  =
	\TIKZ[scale=.35]{
        \node[above] at (1,2.5){$\phantom{m}$};
	\fill[ssyt2] (0,0) to (5,0) to (5,1) to (4,1) to (4,2) to (0,2) to (0,0);
	\Tableau{{1,1,1,1,1},{2,2,2,2}}
	\draw[thick, saffron] (0,0) to (5,0) to (5,1) to (4,1) to (4,2) to (0,2) to (0,0);
	\draw[|-|] (0,-.5) to node[below]{$m$} +(4,0);
	}.
\]
In this case all but one of the edges given by $f_2$ are contracted in the crystal skeleton to give a pentagon.
The edge $J$ in case~\ref{case3} is length-increasing by Corollary~\ref{corollary.increase}. It can be checked explicitly that the edge $J'$ is
length-increasing.
\end{proof}

\begin{remark}
\label{remark.commutations duality}
We note that the commutation relations stated in Theorem~\ref{theorem.commutation} satisfy a duality. Reversing all edges and
mapping edge label $I=[i,i+2m]$ to $[n+1-i-2m,n+1-i]$ as in the Lusztig involution, we have that:
\begin{itemize}
\item cases~\ref{case1}, \ref{case2ai}, \ref{case2bi} and~\ref{case2biii} are self-dual;
\item case~\ref{case2aii} is dual to case~\ref{case3a};
\item case~\ref{case2bii} is dual to case~\ref{case3b}.
\end{itemize}
\end{remark}

\subsection{Sub-skeleton with shortest descent composition}\label{section.subcrystal}

Let $\CS(\lambda)$ be a crystal skeleton. The standardization of the highest weight element in $B(\lambda)_n$ has a descent 
composition with $\ell(\lambda)$ parts; note that this is the minimal length possible for compositions appearing in $\CS(\lambda)$.
This shows that 
\[ 
	\ell(\lambda) =\min\{\len(\Des(T)) \mid T \in \CS(\lambda)\}. 
\] 
Let $\CS(\lambda)_{\ell(\lambda)}$ be the induced subgraph
of $\CS(\lambda)$ consisting of the vertices $\{T \in \CS(\lambda) \mid \len(\Des(T)) =\ell(\lambda)\}$.

\begin{theorem}
\label{theorem.B short}
$\CS(\lambda)_{\ell(\lambda)}$ is isomorphic to the crystal graph $B(\lambda)_{\ell(\lambda)}$. Under this isomorphism:
\begin{itemize}
\item The vertex $T\in \CS(\lambda)_{\ell(\lambda)}$ is mapped to $b\in B(\lambda)_{\ell(\lambda)}$ by replacing all letters in 
$\alpha^{(j)}$ in $T$ in~\eqref{equation.alpha} by $j$, where $\alpha=\Des(T)$. Under this map, $\Des(T)$ becomes $\wt(b)$. \smallskip
\item The edge label $I= [i,i+2m] = I^-\cup \{i+m\} \cup I^+$ with $|I^+|=|I^-|$ in $T\stackrel{I}{\longrightarrow} T'$ in $\CS(\lambda)_{\ell(\lambda)}$ becomes 
$b \stackrel{j}{\longrightarrow} b'$ in $B(\lambda)_{\ell(\lambda)}$ if $I^-\cup \{i+m\} \subseteq \alpha^{(j)}$ and $I^+ \subseteq \alpha^{(j+1)}$.
\end{itemize}
\end{theorem}

\begin{proof}
Elements in $B(\lambda)_{\ell(\lambda)}$ are semistandard tableaux of shape $\lambda$ in the alphabet $\{1,2,\ldots, \ell(\lambda)\}$. Since 
$\lambda$ has $\ell(\lambda)$ parts, the first column of every tableau in $B(\lambda)_{\ell(\lambda)}$ is filled with the numbers
$1,2,\ldots,\ell(\lambda)$. In particular, this means that for every $1\leqslant j <\ell(\lambda)$, there is a bracketed $j+1,j$ pair. Hence by
Proposition~\ref{proposition.quasi edges} each crystal operator $f_j$ in $B(\lambda)_{\ell(\lambda)}$ moves between quasi-crystal components and 
hence is associated with an arrow in $\mathsf{CS}(\lambda)$. Furthermore, the weight $\mathsf{wt}(b)$ for $b\in B(\lambda)_{\ell(\lambda)}$ is
equal to the descent composition since the first column contains the letters $1,2,\ldots,\ell(\lambda)$ and has length $\ell(\lambda)$ in this case. 
In a crystal, $\mathsf{wt}(f_jb) = \mathsf{wt}(b)-\alpha_j$, where $\alpha_j$ is the $j$-th simple root with 1 in position $j$ and $-1$ in position $j+1$. 
These correspond to the length-preserving edges in Theorem~\ref{theorem.descent composition}.
Hence the crystal operators $f_j$ in $B(\lambda)_{\ell(\lambda)}$ correspond to an edge in $\mathsf{CS}(\lambda)_{\ell(\lambda)}$.

Conversely, suppose $T \stackrel{I}{\longrightarrow} T'$ is an edge in $\CS(\lambda)_{\ell(\lambda)}$ with 
$I \subseteq \alpha^{(j)} \cup \alpha^{(j+1)}$. By Corollary~\ref{corollary.increase}, $I$ is maximal under containment among all Dyck pattern 
intervals $J$ for $T$ with $J \subseteq \alpha^{(j)} \cup \alpha^{(j+1)}$. Under the map replacing all letters in $\alpha^{(k)}$ by $k$, this edge
becomes $f_j$.
\end{proof}

\begin{example}
The subgraph $\CS(3,2,1)_3$ in $\CS(3,2,1)$ is indicated by bold edges in Figure~\ref{figure.CS321}.
Figure~\ref{figure.CS321-rest} also shows $\CS(3,2,1)_3$, which is isomorphic to the crystal $B(3,2,1)$ of type $A_2$.
\begin{figure}
\begin{center}
\definecolor{permutation}{HTML}{2660A4}
\includegraphics{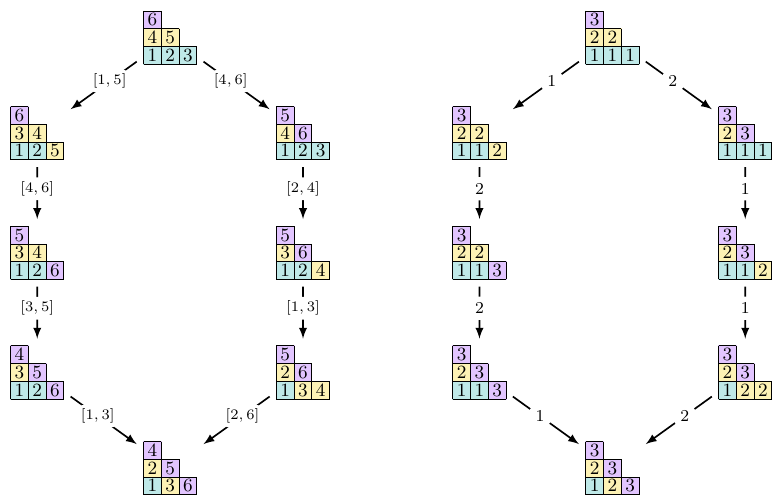}
\end{center}
\caption{Left: Crystal skeleton subgraph $\mathsf{CS}(3,2,1)_3$ of Figure~\ref{figure.CS321}. Right: Crystal $B(3,2,1)_3$ of type $A_2$.
\label{figure.CS321-rest}}
\end{figure}
\end{example}

Recall the string lengths $\varphi_i$ and $\varepsilon_i$ in a crystal defined in~\eqref{equation.string lengths}.
We can use Theorem \ref{theorem.B short} to determine string lengths on $\CS(\lambda)_{\ell(\lambda)}$ as well. 

\begin{corollary}
\label{corollary.string}
Let $T \in \CS(\lambda)_{\ell(\lambda)}$ which maps to $b \in B(\lambda)_{\ell(\lambda)}$ under the isomorphism in Theorem~\ref{theorem.B short}
and set $\alpha = \Des(T)$.
Suppose $T \stackrel{I}{\longrightarrow} T'$ is an edge in $\CS(\lambda)_{\ell(\lambda)}$ with $I \subseteq \alpha^{(i)} \cup \alpha^{(i+1)}$.
Then
\[
	\varphi_i(b) = \min I -\min \alpha^{(i)} +1 \quad \text{and} \quad \varepsilon_i(b) = \max \alpha^{(i+1)} - \max I.
\]
\end{corollary}

\begin{proof}
By self-similarity (see Section~\ref{section.self similar}), this can be checked explicitly on a semistandard Young tableau
with entries $i$ and $i+1$. If this tableau has $m$ letters $i+1$ in the second row and $\ell$ letters $i$, then
$\varphi_i(b) = \ell-m = \min I -\min \alpha^{(i)}+1$. Similarly, if the tableau has $\ell'$ letters $i+1$ in the first row, then
$\varepsilon_i(b) = \ell' = \max \alpha^{(i+1)} - \max I$.
\end{proof}

\begin{proposition}
\label{proposition.incoming and outgoing}
Let $T \in \CS(\lambda)$. Then the following is true:
\begin{enumerate}
\item If $\Des(T) \neq \lambda$, then $T$ has at least one incoming edge which is length-preserving or length-increasing.
\item If $\rev(\Des(T)) \neq \lambda$, then $T$ has at least one outgoing edge which is length-preserving or length-decreasing.
\end{enumerate}
\end{proposition}

\begin{proof}
We prove (2) as (1) follows by applying Lusztig involution.

Assume that all outgoing edges from $T$ are length-increasing. Then for every Dyck pattern interval $I=[i,i+2m]$, the tableau $T$ contains
a rectangle on $[i,i+2m+1]$ by Theorem~\ref{theorem.descent composition}. In particular, this means that if
$\Des(T) = (\alpha_1,\ldots,\alpha_\ell)$, we must have $\alpha_\ell \geqslant \alpha_{\ell-1} \geqslant \cdots \geqslant \alpha_1$, so
that $\rev(\Des(T))$ is a partition. Let $b \in \mathsf{QYT}(\lambda)$ be the quasi-Yamanouchi tableau corresponding to $T$
under the bijection in Lemma~\ref{lemma.QY}. If in the first column of $b$ a letter $1<j \leqslant \ell$ is missing, then there is a Dyck
pattern interval $I=[i,i+2m] \subseteq \alpha^{(j-1)} \cup \alpha^{(j)}$, which does not contain a rectangle $[i,i+2m+1]$, contradicting the
fact that all Dyck pattern intervals contain such a rectangle. Hence every letter $1\leqslant j \leqslant \ell$ appears in the first column of $b$,
which implies that $\ell$ is the number of parts of $\lambda$. Therefore $T \in \CS(\lambda)_{\ell(\lambda)}$. Recall from Theorem~\ref{theorem.B short}
that  $\CS(\lambda)_{\ell(\lambda)}$ is isomorphic to $B(\lambda)_{\ell(\lambda)}$. Every vertex in the crystal $B(\lambda)_{\ell(\lambda)}$ has an
outgoing edge except for the lowest weight vector of weight $\rev(\lambda)$. Hence under the isomorphism, every vertex 
$T\in \CS(\lambda)_{\ell(\lambda)}$ has a length-preserving outgoing edge unless $\Des(T)=\rev(\lambda)$. By assumption
$\Des(T)\neq \rev(\lambda)$, proving the claim.
\end{proof}

\section{Axiomatic characterization of the crystal skeleton}
\label{section.axioms}

In this section, we give an axiomatic characterization of the crystal skeleton. Graphs satisfying these
axioms are called $\CS$-graphs.
We begin in Section~\ref{section.GLn axioms} by stating the $GL_n$-version of the axioms. 
In Section~\ref{section.CS as CS-graph}, we show that crystal skeletons satisfy the $GL_n$-axioms.
We state dual versions of the $GL_n$-axioms in Section~\ref{section.axiom dual} and
show that $\CS$-graphs satisfy $S_n$ branching properties in Section~\ref{section.axiom branching}.
In Section~\ref{section.axiom uniqueness}, we show that the $GL_n$-axioms uniquely specify a $\CS$-graph
for each partition $\lambda$. This shows that the axioms characterize crystal skeletons. In Section~\ref{section.axiom Sn},
we give $S_n$ versions of the axioms and conclude in Section~\ref{section.axiom local} with local axioms.

\subsection{$GL_n$-axioms}
\label{section.GLn axioms}

Suppose $(V,E)$ is a finite, directed graph with vertex set $V$ and edge set $E=\{vw \mid v,w\in V\}$. Fix $n \in \ZZ_{\geqslant 1}$. 
Let $G$ be a vertex- and edge-labeled graph as follows:
\begin{itemize}
    \item The vertices are labeled by compositions $\alpha=(\alpha_1,\ldots,\alpha_\ell)$ of $n$, 
so that the labeled vertex set is $V_L=\{(v,\alpha) \mid v\in V\}$ with $\alpha \models n$. Sometimes we identify $\alpha$ with the partition 
$(\alpha^{(1)},\ldots,\alpha^{(\ell)})$ of $[n]$, where \[ \alpha^{(i)} = \{\alpha_1+\cdots + \alpha_{i-1}+1,\alpha_1+\cdots + \alpha_{i-1}+2,\ldots,
\alpha_1+\cdots+\alpha_i\}. \]
\item The edges are labeled by (odd-length) intervals  $I \subseteq [n]$, so that the labeled edge set is 
$E_L=\{(vw,I) \mid vw \in E\}$.
\end{itemize}

With $G$ as above, we define the Lusztig involution on $G$ as follows. 

\begin{definition}[Lusztig involution] \label{def.LIonCS}
The \defn{Lusztig involution} $\mathcal{L}_n$ on $G$ is defined by 
\begin{itemize}
    \item relabeling the vertices by replacing $(v,\alpha)$ by 
$(v,\rev(\alpha))$, and 
\item reversing all edge directions and changing the edge label 
\[ I=[a,b] \quad \quad \textrm{ to } \quad \quad I^{\LL}:=[n+1-b,n+1-a]. \]
\end{itemize}
\end{definition}
\begin{example} 
Suppose $n=6$. Then $\LL_6$ acts by
\[ \LL_6: \bigg( (v,(3,2,1)) \xrightarrow{[4,6]} (w,(3,1,2)) \bigg) \longmapsto \bigg((w,(2,1,3)) \xrightarrow{[1,3]} (v,(1,2,3)) \bigg). \]
\end{example}

We will also restrict our graphs via \defn{branched graphs} as follows.

\begin{definition}[Branched graph]
Define $G_{[1,n-1]} = (V_L',E_L')$, where 
\[ V_L'=\{ (v, \alpha \setminus \{ n \}) \mid (v, \alpha) \in V_L \} \quad \textrm{ and } \quad E_L' = \{(vw,I) \mid (vw,I)\in E_L, \  I \subseteq [n-1] \}.\]
\end{definition}

\begin{example}
In Figure~\ref{figure.CS321}, the graph $G_{[1,5]}$ is isomorphic to the portion of $\CS(3,2,1)$ shaded in gray. To obtain $G_{[1,5]}$ from the 
gray subgraph, keep the edge labels the same and replace the vertex $\Des(T)$ for $T \in \SYT(3,2,1)$ with $\Des(T_{[1,5]})$ as 
in Section~\ref{section.self similar}.
\end{example}

We are now ready to state the axioms for crystal skeletons. Graph isomorphisms are considered to preserve vertex and edge labels.
\begin{axiom}[Axioms for crystal skeletons]
\label{axioms}
Let $n$ be a positive integer and $G$ be a finite, connected, vertex- and edge-labeled graph with labeled vertex set 
$V_L$ and labeled edge set $E_L$ as above. 
We call $G$ a \defn{$\mathsf{CS}$-graph} if the following axioms hold:

\begin{enumerate}[label=\bf {A}\arabic{enumi}.,ref=\bf {A}\arabic{enumi}]\setcounter{enumi}{-1}
\item
\label{axiom.intervals} \emph{(Intervals)} 
Suppose $(v,\alpha) \xrightarrow{~I~} (w, \beta)$ is an edge in $G$. Then the interval $I \subseteq [n]$ satisfies
\begin{enumerate}[label=(\alph*)]
\item \label{axiom.interval types-decomp} $I=I^- \cup \{k\} \cup I^+$ where  $|I^-|=|I^+|>0$, 
\item \label{axiom.interval types-out} $I^- \cup \{k\} \subseteq \alpha^{(j)}$ and $I^+ \subseteq \alpha^{(j+1)}$ for some $1\leqslant j < \ell$, 
where $\ell$ is the length of $\alpha$.
\end{enumerate}
\[
\TIKZ[scale=.75, thick]{
\draw[thin, black!20] (-1.5, 0) to (8.5,0);
\draw[Bracket-Bracket] (0,0) to (3,0);
\draw[Bracket-Bracket] (4,0) to (7,0);
\begin{scope}[every node/.style={draw, thin, inner sep=1.75pt, fill=white, rounded corners}]
\node at (1.5,0) {\small $I^-\vphantom{I_-}$};
\node at (3.5,0) {\small $k\vphantom{I_-^-}$};
\node at (5.5,0) {\small $I^+\vphantom{I_-}$};
\end{scope}
\draw[rounded corners=2, sapphire] (-1+.2, .25) to (-1+.2, .5) to 
	node[pos=.15, inner sep=.5pt, fill=white]{$\ \!\cdots$} 
	node[above] {$\alpha^{(j)}$} (4-.15-.05, .5) to (4-.15-.05, .25);
\draw[rounded corners=2, sapphire] (8+.2, .25) to (8+.2, .5) to 
	node[pos=.15, inner sep=.5pt, fill=white]{$\ \!\cdots$} 
	node[above] {$\alpha^{(j+1)}$} (4-.15+.05, .5) to (4-.15+.05, .25);
}
\]

\item 
\label{axiom.outgoing}(\emph{Outgoing edges}) 
For each $(v, \alpha) \in V_L$ and each interval $I$ satisfying~\ref{axiom.intervals} exactly one of the following holds: Either
\begin{enumerate}[label=(\roman*)]
\item there is exactly one outgoing edge $(v, \alpha) \xrightarrow{~I~} \fbox{\small $\phantom{(v, \alpha)}$}$ \ labeled by $I$; or 
\item \label{axiom.outgoing.incoming} 
there is an incoming edge $(u,\gamma) \xrightarrow{~J~} (v, \alpha)$ with $J \subseteq I$ and $\gamma$ dominates $\alpha$.
\end{enumerate}

\item\label{axiom.labels} (\emph{Labels})
Let $(v,\alpha) \xrightarrow{~I~} (w,\beta)$ be an edge as in \ref{axiom.intervals}. Given $\alpha$ and $I$, then $\beta$ must be of the form 
\[\beta = (\alpha^{(1)}, \dots, \alpha^{(j-2)}, 
	\TIKZ[yscale=.45]{\draw[thin] (0,.5) to (0,0) to (2,0) to +(0,.5); \node[circle, draw, inner sep=1pt] at (1,.5) {$*$};}\ ,
	\alpha^{(j+2)}, \dots, \alpha^{(\ell)}),\]
	\def\circasterisk{\TIKZ{\node[circle, draw, inner sep=1pt] at (1,.5) {$*$};}}
where \circasterisk\  is one of
\begin{enumerate}[label=\bf \Roman{enumii}.,ref=\bf \Roman{enumii}]
\item \label{axiom.length pres}(\emph{length-preserving}) \ 
	$\circasterisk = (\alpha^{(j-1)}, 			\hspace{5pt}
				\alpha^{(j)} \setminus \{k\}, 	\hspace{5pt}
				\alpha^{(j+1)} \cup \{k\})$;
\item  \label{axiom.length inc} (\emph{length-increasing}) \ 
	$\circasterisk = (\alpha^{(j-1)}, 			\hspace{5pt}
				\alpha^{(j)} \setminus \{k\}, 	\hspace{5pt}
				I^+ \cup \{k\}, 				\hspace{5pt}
				\alpha^{(j+1)} \setminus I^+)$; 
	or
\item \label{axiom.length dec} (\emph{length-decreasing})  \ 
	$\circasterisk = (\alpha^{(j-1)} \cup I^-,  		\hspace{5pt}
				\alpha^{(j+1)} \cup \{k\})$, 
		given $\alpha^{(j)}=I^-\cup \{k\}$.
\end{enumerate}

\item\label{axiom.fans} (\emph{Fans})
Suppose $(v,\alpha) \xrightarrow{~I~} (w,\beta)$ is an edge in $G$ with $\beta$ satisfying axiom~\ref{axiom.labels}~\ref{axiom.length inc}.
Then this edge is part of a fan
\begin{center}
\TIKZ[>=latex,line join=bevel, yscale=.8]{
\node (N0) at (0,0) [draw,draw=none] {$(v,\alpha)$};
\node (N1) at (-7,3) [draw,draw=none] {$(w_1,\beta_1)$};
\node (N2) at (-4,3) [draw,draw=none] {$(w_2,\beta_2)$};
\node (N4) at (0,3) [draw,draw=none] {$(w_{m-1},\beta_{m-1})$};
\node (N5) at (4,3) [draw,draw=none] {$(w_m,\beta_m)$};
  \draw[->, bend left=10] (N0)  to node[opLabel, pos=.5](I1){$I_1$} (N1);
  \draw[->, bend left=10] (N0)  to node[opLabel, pos=.5](I2){$I_2$} (N2);
  \draw[->] (N0)  to node[opLabel, pos=.5](Ia-1){$I_{m-1}$} (N4);
  \draw[->] (N0)  to node[opLabel, pos=.5]{$I_m$} (N5);
  \draw[->] (N1)  to node[opLabel, pos=.5]{$J_1$} (N2);
  \draw[->] (N2)  to node[opLabel, pos=.5]{$\cdots$} (N4);
  \draw[->] (N4)  to node[opLabel, pos=.5]{$J_{m-1}$} (N5);
  \draw[->] (N5) to[bend left=30] node[opLabel, pos=.5]{$J_m$} (N0);
\path (I2) to node{$\cdots$} (Ia-1);
}
\end{center}
where $(w,\beta)$ is one of the $(w_j,\beta_j)$ and 
\begin{align*}
    I_1&=[i,i+2m], \quad I_2=[i+1,i+2m-1],  \ldots,  I_m=[i+m-1,i+m+1] \\
    J_1&=[i+2m-1,i+2m+1], \quad J_2=[i+2m-2,i+2m],\ldots,  J_{m}=[i+m,i+m+2]
\end{align*}
for some $i$ and $m$. The edges labeled $I_1,\ldots, I_m$ satisfy
axiom~\ref{axiom.labels}~\ref{axiom.length inc}, the edges labeled $J_1,\ldots,J_{m-1}$ satisfy axiom~\ref{axiom.labels}~\ref{axiom.length pres} and
the edge labeled $J_m$ satisfies axiom~\ref{axiom.labels}~\ref{axiom.length dec}.

\item\label{axiom.lusztig} (\emph{Strong Lusztig involution}) 
$G$ and $G_{[1,n-1]}$ are invariant under Lusztig involution:
\begin{enumerate}[label=(\alph*)]
\item \label{axiom.lusztig n}
$\mathcal{L}_n(G) \cong G$;
\item \label{axiom.lusztig n-1}
$\mathcal{L}_{n-1}(G_{[1,n-1]}) \cong G_{[1,n-1]}$.
\end{enumerate}

\item\label{axiom.crystal} (\emph{Top subcrystal}) 
Let $s = \min \{ \len(\alpha) \mid (v,\alpha) \in V_L\}$ be the minimal length of all vertex labels.
Let $G_s$ be the induced subgraph of $G$ with vertex set 
\[ \{(v,\alpha) \in V_L \mid \len(\alpha) =s\}. \] 
Then $G_s$ is isomorphic to the crystal graph $B(\lambda)_s$ for some partition $\lambda$ with $\ell(\lambda)=s$.
Under this isomorphism, the vertex labels become the weights in $B(\lambda)_s$ and an edge labeled $I$ satisfying axiom~\ref{axiom.intervals}
becomes $f_j$ in $B(\lambda)_s$.
\end{enumerate}
\end{axiom}

\begin{remark}
Note that the crystal $B(\lambda)_s$ appearing in axiom~\ref{axiom.crystal} has a unique vertex $(u,\lambda)$ with label $\lambda$.
Hence $G$ itself has a unique such vertex and so we will write $G_\lambda$ for the $\CS$-graph $G$ containing $B(\lambda)_s$.
\end{remark}

\subsection{Crystal skeleton as a $\CS$-graph}
\label{section.CS as CS-graph}

We first show that the crystal skeleton is indeed a $\CS$-graph.

\begin{theorem}
\label{theorem.CS is CS-graph}
    For any partition $\lambda$, the crystal skeleton $\CS(\lambda)$ is a $\CS$-graph with vertex labeling given by descent 
    compositions and edge labeling given by Dyck pattern intervals. 
\end{theorem}

\begin{proof}
Let $G = \CS(\lambda)$, with edge and vertex labelings as stated. We will show that $G$ satisfies all of the axioms of a $\CS$-graph.
By Theorem~\ref{theorem.Dyck pattern}, the edges in $G$ are labeled by Dyck pattern intervals $I$, which are intervals of odd length.
This implies that $G$ satisfies axiom~\ref{axiom.intervals}~\ref{axiom.interval types-decomp}. Axiom~\ref{axiom.intervals}~\ref{axiom.interval types-out} 
follows from Lemma~\ref{lemma:I.inside.alpha}.

Next we show that axiom~\ref{axiom.outgoing} holds. Let $(T,\alpha)$ be a vertex in $G$ with $T\in \SYT(\lambda)$ and 
$\alpha=\Des(T)$. Let $I=I^- \cup \{k\} \cup I^+$ be an interval satisfying Axiom~\ref{axiom.intervals} such that $I^- \cup \{k\} \subseteq \alpha^{(j)}$ 
and $I^+ \subseteq \alpha^{(j+1)}$. There is an edge labeled $I$ out of vertex $(T,\alpha)$ in $\mathsf{CS}(\lambda)$ if and only if 
$I$ is a Dyck pattern interval by Theorem~\ref{theorem.Dyck pattern}.

First consider the case that there is an edge with interval $I$. By Theorem~\ref{theorem.self similar}, $G$ restricted to only edges labeled by intervals
contained in $I$ is isomorphic to $\mathsf{CS}(m+1,m)$ with $m=|I^-|$ since $P(\pi|_I)$ has shape $(m+1,m)$  for $\pi=\row(T)$.
In $\mathsf{CS}(m+1,m)$, the vertex $(T,\alpha)$ has label $(m+1,m)$. Note that $\mathsf{CS}(m+1,m)$ has a unique vertex with label $(m+1,m)$
and this is the lexicographically largest label. Hence every incoming edge $(T',\alpha') \stackrel{J}{\longrightarrow} (T,(m+1,m))$ in
$\mathsf{CS}(m+1,m)$ has a label $\alpha'$ that is smaller in dominance order than $(m+1,m)$. This remains true when embedded
into the larger graph $G$. This confirms axiom~\ref{axiom.outgoing} in this case.

Next consider the case that there is no edge labeled $I$ out of $(T,\alpha)$. This means that the standard tableau $T$ does not contain a 
Dyck pattern $\pi|_I$, where $\pi=\row(T)$. Hence $P(\pi|_I)$ has shape $\mu$ with $\mu>_D (m+1,m)$ and $\mu_2<m$. 
Furthermore, $\mu_2>0$ since by assumption $I^- \cup \{k\} \subseteq \alpha^{(j)}$ and $I^+ \subseteq \alpha^{(j+1)}$ so that there is at least one 
descent. Hence in particular, $m>1$.
In $\mathsf{CS}(\mu)$, there is an edge labeled $[m+1,m+3]$ from the standard tableau with descent composition $(m+2,m-1)$ to that of
descent composition $(m+1,m)$. Since $(m+2,m-1)>_D (m+1,m)$, this proves axiom~\ref{axiom.outgoing} in this case.

Axiom~\ref{axiom.labels} holds by Theorem~\ref{theorem.descent composition}.

Axiom~\ref{axiom.fans} holds by Corollary~\ref{corollary.fans longer}.

Axiom~\ref{axiom.lusztig} \ref{axiom.lusztig n} holds by Corollary~\ref{corollary.lusztig}.

Axiom~\ref{axiom.lusztig} \ref{axiom.lusztig n-1} follows from Theorem~\ref{theorem.branching} and Corollary~\ref{corollary.lusztig}.

Axiom~\ref{axiom.crystal} holds by Theorem~\ref{theorem.B short}. 
\end{proof}

\subsection{Dual version of the axioms}
\label{section.axiom dual}

By the Lusztig involution, the axioms~\ref{axiom.intervals}-\ref{axiom.fans} have dual versions:
\begin{enumerate}[label=\bf {A}\arabic{enumi}'.,ref=\bf {A}\arabic{enumi}']\setcounter{enumi}{-1}
\item \label{axiom.intervals prime} \emph{(Intervals)} 
Suppose $(v,\alpha) \xrightarrow{~I~} (w, \beta)$ is an edge in $G$. Then the interval $I \subseteq [n]$ satisfies
\begin{enumerate}
\item[(b']) \label{ax:interval types-in} $I^- \subseteq \beta^{(j')}$ and $I^+ \cup \{k\} \subseteq \beta^{(j'+1)}$ for some $1\leqslant j'<\ell'$, where $\ell'$ is the length
of $\beta$.
\end{enumerate}
\[
\TIKZ[scale=.75, thick]{
\draw[thin, black!20] (-1.5, 0) to (8.5,0);
\draw[Bracket-Bracket] (0,0) to (3,0);
\draw[Bracket-Bracket] (4,0) to (7,0);
\begin{scope}[every node/.style={draw, thin, inner sep=1.75pt, fill=white, rounded corners}]
\node at (1.5,0) {\small $I^-\vphantom{I_-}$};
\node at (3.5,0) {\small $k\vphantom{I_-^-}$};
\node at (5.5,0) {\small $I^+\vphantom{I_-}$};
\end{scope}
\draw[rounded corners=2, sapphire] (-1+.2, .25) to (-1+.2, .5) to 
	node[pos=.15, inner sep=.5pt, fill=white]{$\ \!\cdots$} 
	node[above] {$\alpha^{(j)}$} (4-.15-.05, .5) to (4-.15-.05, .25);
\draw[rounded corners=2, sapphire] (8+.2, .25) to (8+.2, .5) to 
	node[pos=.15, inner sep=.5pt, fill=white]{$\ \!\cdots$} 
	node[above] {$\alpha^{(j+1)}$} (4-.15+.05, .5) to (4-.15+.05, .25);
\draw[rounded corners=2, plum] (-1-.2, -.25) to (-1-.2, -.5) to 
	node[pos=.15, inner sep=.5pt, fill=white]{$\ \!\cdots$} 
	node[below] {$\beta^{(j')}$} (3.15-.05, -.5) to (3.15-.05, -.25);
\draw[rounded corners=2, plum] (8-.2, -.25) to (8-.2, -.5) to 
	node[pos=.15, inner sep=.5pt, fill=white]{$\ \!\cdots$} 
	node[below] {$\beta^{(j'+1)}$} (3.15+.05, -.5) to (3.15+.05, -.25);
}
\]

\item (\emph{Incoming edges}) 
\label{axiom.incoming}
For each $(w, \beta) \in V_L$ and each interval $I$ satisfying~\ref{axiom.intervals prime} 
exactly one of the following holds: Either
\begin{enumerate}[(i)]
\item there is exactly one incoming edge $\fbox{\small $\phantom{(v, \alpha)}$} \xrightarrow{~I~} (w,\beta)$ \ labeled by $I$; or
\item there is an outgoing edge of the form $(w,\beta) \xrightarrow{~J~} (u, \gamma)$ with $J \subseteq I$ and $\gamma <_D \beta$.
\end{enumerate}

\item (\emph{Labels})
\label{axiom.labels incoming}
Let $(v,\alpha) \xrightarrow{~I~} (w,\beta)$ be an edge. Given $\beta$ and $I$, then $\alpha$ must be of the form
\[\alpha = (\beta^{(1)}, \dots, \beta^{(j-2)}, 
	\TIKZ[yscale=.45]{\draw[thin] (0,.5) to (0,0) to (2,0) to +(0,.5); \node[circle, draw, inner sep=1pt] at (1,.5) {$\diamond$};}\ ,
	\beta^{(j+3)}, \dots, \beta^{(\ell')}),\]
	\def\circdiamond{\TIKZ{\node[circle, draw, inner sep=1pt] at (1,.5) {$\diamond$};}}
where \circdiamond\ is one of
\begin{enumerate}[label=\bf \Roman{enumii}.,ref=\bf \Roman{enumii}]
\item (\emph{length-preserving}) \ 
	$\circdiamond = (\beta^{(j-1)}, 			\hspace{5pt}
				\beta^{(j)} \cup \{k\}, 		\hspace{5pt}
				\beta^{(j+1)} \setminus \{k\}, 	\hspace{5pt}
				\beta^{(j+2)})$;
\item (\emph{length-increasing}) \ 
	$\circdiamond = (\beta^{(j-1)}, 			\hspace{5pt}
				\beta^{(j)} \cup \{k\}, 		\hspace{5pt}
				I^+ \cup \beta^{(j+2)})$, 
			given $\beta^{(j+1)}=I^+\cup \{k\}$; or
\item (\emph{length-decreasing})  \ 
	$\circdiamond = (\beta^{(j-1)} \setminus I^-, 	\hspace{5pt}
				I^- \cup \{k\}, 				\hspace{5pt}
				\beta^{(j)} \setminus \{k\}, 	\hspace{5pt}
				\beta^{(j+1)}, 				\hspace{5pt}
				\beta^{(j+2)})$.
\end{enumerate}

\item (\emph{Fans})
\label{axiom.fans dual}
Suppose $(v,\alpha) \xrightarrow{~I~} (w,\beta)$ is an edge in $G$ with $\beta$ satisfying axiom~\ref{axiom.labels incoming}~\ref{axiom.length dec}.
Then this edge is part of a fan
\begin{center}
\raisebox{-0.8cm}{
\scalebox{0.8}{
\begin{tikzpicture}[>=latex,line join=bevel,]
\node (N0) at (0,0) [draw,draw=none] {$(w,\beta)$};
\node (N1) at (-6,3) [draw,draw=none] {$(v_1,\alpha_1)$};
\node (N2) at (-3,3) [draw,draw=none] {$(v_2,\alpha_2)$};
\node (N4) at (0,3) [draw,draw=none] {$(v_{m-1},\alpha_{m-1})$};
\node (M) at (-0.8,1.5) [draw,draw=none] {$\cdots$};
\node (M1) at (-1.5,3) [draw,draw=none] {$\cdots$};
\node (N5) at (4,3) [draw,draw=none] {$(v_m,\alpha_m)$};
  \draw[->] (N1)  to node[opLabel]{$I_1$} (N0);
  \draw[->] (N2)  to node[opLabel]{$I_2$} (N0);
  \draw[->] (N4)  to node[opLabel]{$I_{m-1}$} (N0);
  \draw[->] (N5)  to node[opLabel]{$I_m$} (N0);
  \draw[->] (N2)  to node[opLabel]{$J_1$} (N1);
  \draw[->] (N5)  to node[opLabel]{$J_{m-1}$} (N4);
  \draw[->] (N0) to[bend right=30] node[opLabel]{$J_m$} (N5);
\end{tikzpicture}
}}
\end{center}
where $(v,\alpha)$ is one of the $(v_j,\alpha_j)$ and $I_1=[i,i+2m]$, $I_2=[i+1,i+2m-1],\ldots,I_m=[i+m-1,i+m+1]$,
$J_1=[i-1,i+1]$, $J_2=[i,i+2],\ldots,J_{m}=[i+m-2,i+m]$ for some $i$ and $m$.
The edges labeled $I_1,\ldots, I_m$ satisfy axiom~\ref{axiom.labels incoming}~\ref{axiom.length dec}, the edges labeled $J_1,\ldots,J_{m-1}$ satisfy 
axiom~\ref{axiom.labels incoming}~\ref{axiom.length pres} and the edge labeled $J_m$ satisfies axiom~\ref{axiom.labels incoming}~\ref{axiom.length inc}.
\end{enumerate}

\subsection{Branching properties}
\label{section.axiom branching}

Recall that for a $\mathsf{CS}$-graph $G = (V_L,E_L)$ with labeled vertex set $V_L=\{(v,\alpha) \mid v\in V\}$ and labeled edge set
$E_L=\{(vw,I) \mid v,w\in V\}$ (where labels satisfy $\alpha \models n$ and $I\subseteq [n]$ as above), we defined:
\begin{itemize}
    \item $G_{[1,n-1]} = (V_L',E_L')$, with labeled vertex set $V_L'=\{ (v, \alpha \setminus \{ n \}) \mid (v, \alpha) \in V_L \}$, and 
    \item labeled edge set 
    $E_L' = \{(vw,I) \mid (vw,I)\in E_L, I \subseteq [n-1] \}$.
\end{itemize}
In addition, we define another graph:
\begin{itemize}
    \item $G_{[2,n]} = (V_L'', E_L'')$ with labeled vertex set $V_L''=\{(v, \alpha \setminus \{ 1 \}) \mid (v, \alpha) \in V_L \}$, and 
    \item labeled edge set 
    $E_L''=\{(vw,I-{\bf 1}) \mid (vw,I) \in E_L, I \subseteq [2,n]\}$,
    where $I - {\bf 1}$ decreases the interval bounds by 1 (e.g. $[i,j] - {\bf 1 } = [i-1,j-1]$). 
\end{itemize}

We now study properties of these subgraphs.

\begin{theorem}
\label{theorem.G1 and Gn}
    Suppose $G:=G_\lambda$ is a $\CS$-graph. Then: 
    \begin{enumerate}
    \item The graphs $G_{[1,n-1]}$ and $G_{[2,n]}$ are isomorphic.
    \item $G_{[1,n-1]}$ and $G_{[2,n]}$ are unions of $\CS$-graphs.
    \item We have
    \[
    	G_{[1,n-1]} \cong \bigcup_{\substack{\lambda^-\\ |\lambda/\lambda^-|=1}} G_{\lambda^-}.
    \]
    \end{enumerate}
    \label{smaller}
\end{theorem}

\begin{proof}
Recall that the Lusztig involution $\mathcal{L}_n$ reverses all labels $\alpha=(\alpha_1,\ldots,\alpha_\ell) \mapsto (\alpha_\ell,\ldots,\alpha_1)$, 
reverses all arrows, and maps the edge label $I$ to $\{n+1-a \mid a\in I\}$. Since $G$ is a $\CS$-graph, Axiom~\ref{axiom.lusztig} \ref{axiom.lusztig n}
holds so that $\mathcal{L}_n(G) \cong G$. Denote by $\psi^n$ (resp. $\psi^1$) the map that removes all edges labeled $I$ with $n\in I$ 
(resp. $1\in I$) and set $G^n = \psi^n(G)$ and $G^1=\psi^1(G)$. Note that $\psi^n = \mathcal{L}_n \circ \psi^1 \circ \mathcal{L}_n$ so that
\begin{equation}
\label{equation.LL}
	G^n = \psi^n(G) =  \mathcal{L}_n \circ \psi^1 \circ \mathcal{L}_n(G) \cong \mathcal{L}_n \circ \psi^1(G) = \mathcal{L}_n (G^1),
\end{equation}
where the isomorphism applies $\mathcal{L}_n(G) \cong G$.
Define the map $\varphi^n$ on vertex labels $\alpha=(\alpha_1,\ldots,\alpha_\ell)$ by $\alpha \mapsto (\alpha_1,\ldots,\alpha_{\ell-1},\alpha_\ell-1)$
and $\varphi^1$ by mapping vertex labels $\alpha \mapsto (\alpha_1-1,\alpha_2,\ldots,\alpha_\ell)$ and edge labels $I \mapsto I-{\bf 1}$.
Note that $G_{[1,n-1]} = \varphi^n(G^n)$ and $G_{[2,n]} = \varphi^1(G^1)$, so that using~\eqref{equation.LL} we have
\[
	G_{[1,n-1]} = \varphi^n(G^n) \cong \varphi^n \mathcal{L}_n(G^1) = \mathcal{L}_{n-1} \varphi^1(G^1) =  \mathcal{L}_{n-1}(G_{[2,n]}).
\]	
Since $\mathcal{L}_{n-1}$ is an involution, we thus have $G_{[2,n]} \cong \mathcal{L}_{n-1} (G_{[1,n-1]}) \cong G_{[1,n-1]}$ by
axiom~\ref{axiom.lusztig}~\ref{axiom.lusztig n-1}. This proves (1).

Next we show that $G_{[1,n-1]}$ is a union of $\CS$-graphs. The claim that $G_{[2,n]}$ is a union of $\CS$-graphs then follows from (1).

Let 
\[
	(v,\alpha) \stackrel{I}{\longrightarrow} (w,\beta)
\]
be an edge in $G$. This edge is also an edge in $G_{[1,n-1]}$ if $I \subseteq [1,n-1]$.
Axiom~\ref{axiom.intervals} is satisfied for $G_{[1,n-1]}$ since $I$ labels an edge in $G$, which satisfies axiom~\ref{axiom.intervals} by assumption.
If $(v,\alpha)$ is a vertex in $G$, then the corresponding vertex in $G_{[1,n-1]}$ is $(v,\alpha_{[1,n-1]})$, where 
$\alpha_{[1,n-1]} = (\alpha_1,\alpha_2,\ldots,\alpha_\ell-1)$. Note that $\alpha'>_D \alpha$ in $G$ implies that $\alpha'_{[1,n-1]} >_D \alpha_{[1,n-1]}$ in
$G_{[1,n-1]}$. Hence axiom~\ref{axiom.outgoing} holds for $G_{[1,n-1]}$.
If $n \not \in \alpha^{(j+1)}$, then~\ref{axiom.labels} is still satisfied in $G_{[1,n-1]}$ since the conditions are local on 
$\alpha^{(j)}, \alpha^{(j+1)}$.  Hence assume that $n \in \alpha^{(j+1)}$. If $\beta$ in axiom~\ref{axiom.labels} is in case~\ref{axiom.length pres} 
or~\ref{axiom.length dec} in $G$, then it is still in these cases in $G_{[1,n-1]}$ since the presence or absence of the element $n$ in $\alpha$ does not alter the 
form of $\beta$ given in axiom~\ref{axiom.labels}.
This is also true if $\beta$ is in case~\ref{axiom.length inc}, unless $\alpha^{(j+1)} = I^+ \cup \{ n \}$, in which case the edge $(vw, I)$ becomes an edge 
in case~\ref{axiom.length pres} in $G_{[1,n-1]}$, since $\alpha^{(j+1)} \setminus \{ n, I^+ \} = \emptyset$. Hence axiom~\ref{axiom.labels} holds for 
$G_{[1,n-1]}$.

The graph $G$ satisfies the fan axiom~\ref{axiom.fans}. If none of the intervals in the fan contain $n$, they will all still appear in $G_{[1,n-1]}$ and will be 
of the same type in axiom~\ref{axiom.labels}. Hence the fan exists in $G_{[1,n-1]}$ as desired. Note that $J_1=[i+2m-1,i+2m+1]$ 
contains the largest letter among all intervals in the fan in $G$. Now assume that $n\in J_1$, that is, $n=i+2m+1$. In this case the edge labeled $J_1$ 
is missing in $G_{[1,n-1]}$. Furthermore, the edge 
\[ (v,\alpha) \stackrel{I_1}{\longrightarrow} (w_1,\beta_1) \] satisfies $I_1=[n-2m+1,n-1]$ and
$\beta_1 = (\alpha^{(1)},\ldots,\alpha^{(\ell-2)},\alpha^{(\ell-1)} \setminus \{n-m\},\{n-m,\ldots,n-1\},\{n\})$. 
In particular, $\alpha^{(\ell)}=I_1^+\cup\{n\} = [n-m+1,n]$, so that the edge labeled $I_1$ becomes of type~\ref{axiom.length pres} in axiom~\ref{axiom.labels}
in $G_{[1,n-1]}$. This proves that $G_{[1,n-1]}$ satisfies axiom~\ref{axiom.fans}.

By assumption $G$ satisfies axiom~\ref{axiom.lusztig}~\ref{axiom.lusztig n-1}, that is $\mathcal{L}_{n-1}(G_{[1,n-1]}) = G_{[1,n-1]}$.
This is axiom~\ref{axiom.lusztig}~\ref{axiom.lusztig n} for $G_{[1,n-1]}$ since $n$ is replaced by $n-1$ for $G_{[1,n-1]}$.
By part (1), we have $G_{[1,n-1]} \cong G_{[2,n]}$. Acting by $\varphi^{n-1} \circ \psi^{n-1}$ on $G_{[1,n-1]}$ yields $G_{[1,n-2]}$.
Acting by $\varphi^{n-1} \circ \psi^{n-1}$ on $G_{[2,n]}$ yields
\[
	\varphi^{n-1} \psi^{n-1} G_{[2,n]} = \varphi^{n-1} \psi^{n-1} \varphi^1\psi^1 G =\varphi^1 \psi^1 \varphi^n \psi^n G = \varphi^1 \psi^1 G_{[1,n-1]}
	= G_{[2,n-1]}.
\] 
Hence $G_{[1,n-2]} \cong G_{[2,n-1]}$. 
By the same arguments as in the proof of part (1), this is equivalent to $\mathcal{L}_{n-2}(G_{[1,n-2]})
\cong G_{[1,n-2]}$ using $\mathcal{L}_{n-1}(G_{[1,n-1]}) = G_{[1,n-1]}$. This proves axiom~\ref{axiom.lusztig}~\ref{axiom.lusztig n-1}
for $G_{[1,n-1]}$.

By assumption, $G$ satisfies axiom~\ref{axiom.crystal}, hence $G_s$ is isomorphic to $B(\lambda)_s$ for some partition $\lambda$ of length $s$. 
The vertices of $B(\lambda)_s$ can be indexed by semistandard tableaux of shape $\lambda$ in the alphabet $\{1,2,\ldots,s\}$. The label of 
$t \in B(\lambda)_s$ in $G$ is its weight $\mathsf{wt}(t)$. Transitioning from $G$ to $G_{[1,n-1]}$ corresponds to removing the rightmost letter $s$ 
in $t$ since under standardization this letter corresponds to $n$ in the standard tableau. This results in a tableau $t^-$ of shape $\lambda^-$, where 
$\lambda/\lambda^-$ is the box of the corner cell containing the rightmost $s$. 

First consider the case that the length of $\lambda^-$ is still $s$. Then the crystal $B(\lambda^-)_s$ compared to $B(\lambda)_s$ restricted to the 
vertices with a letter $s$ in $\lambda/\lambda^-$ has potentially more edges $f_{s-1}$ (since the $s$ in $\lambda/\lambda^-$ is potentially bracketed 
with an $s-1$). Since $G$ satisfies axiom~\ref{axiom.outgoing}, there is an edge in $G$
corresponding to such an extra edge $f_{s-1}$ in $B(\lambda^-)_s$ compared to $B(\lambda)_s$, however it is of type~\ref{axiom.length inc} 
in axiom~\ref{axiom.labels} and not of  type~\ref{axiom.length pres}. It has $\alpha_{s+1}=1$. Hence $(G_{[1,n-1]})_s$ contains $B(\lambda^-)_s$ and 
the connected component of 
$G_{[1,n-1]}$ containing $B(\lambda^-)_s$ restricted to labels of length $s$ is isomorphic to it. This shows that this component of $G_{[1,n-1]}$ satisfies 
axiom~\ref{axiom.crystal}.

Next consider the case that the length of $\lambda^-$ is $s-1$. In this case $\lambda_s=1$. The vertices of the crystal $B(\lambda^-)_{s-1}$ are
the tableaux in $B(\lambda)_s$ with only one letter $s$ in row $s$. Since $B(\lambda^-)_{s-1}$ only has arrows labeled $1,2,\ldots,s-2$, these are
all contained in $B(\lambda)_s$ on the vertices in $B(\lambda^-)_{s-1}$. The connected component of $G_{[1,n-1]}$ containing $B(\lambda^-)_{s-1}$ 
restricted to labels of length $s-1$ is isomorphic to it. This shows that this component of $G_{[1,n-1]}$ satisfies axiom~\ref{axiom.crystal}.

Part (3) follows from (2) and the fact that the crystal $B(\lambda)_s$ in $G$ breaks into crystals $B(\lambda^-)s$ or $B(\lambda^-)_{s-1}$, where
$|\lambda/\lambda^-|=1$, when passing from $G$ to $G_{[1,n-1]}$.
\end{proof}

\begin{remark}
Note that under the assumption that Lusztig involution holds on $G_\lambda$ (axiom~\ref{axiom.lusztig}~\ref{axiom.lusztig n}), the proof of
Theorem~\ref{theorem.G1 and Gn} (1) shows that Theorem~\ref{theorem.G1 and Gn} (1) is equivalent to axiom~\ref{axiom.lusztig}~\ref{axiom.lusztig n-1}.
Hence Axiom~\ref{axiom.lusztig}~\ref{axiom.lusztig n-1} can be replaced by
\begin{enumerate}
\item[{\bf A4.}(b)] The graphs $G_{[1,n-1]}$ and $G_{[2,n]}$ are isomorphic.
\end{enumerate}
\end{remark}

\subsection{Uniqueness of $\CS$-graphs}
\label{section.axiom uniqueness}

In this section, we show that Axiom~\ref{axioms} uniquely determines $\CS$-graphs for a given $\lambda$, and hence 
$G_\lambda \cong \CS(\lambda)$.

We begin by considering  whether a $\CS$-graph $G$ can be reconstructed from $G_{[1,n-1]}$. To this end, we first make some observations.
As in the proof of Theorem~\ref{theorem.G1 and Gn} define the map $\varphi$ on compositions by
(this map is denoted by $\varphi^n$ in the proof of Theorem~\ref{theorem.G1 and Gn})
\[
	\varphi(\alpha_1,\ldots,\alpha_\ell) = (\alpha_1,\ldots,\alpha_{\ell-1},\alpha_\ell-1).
\]
We extend this to a map $\varphi \colon G \to G_{[1,n-1]}$ by acting on the labels of all vertices and removing all edges labeled by $I$ which contain $n$.
Consider an edge $(v,\alpha) \stackrel{I}{\longrightarrow} (w,\beta)$ in $G$ with $n\not \in I$ and consider its image under $\varphi$:
\begin{equation}
\label{equation.edges in G and Gn}
\begin{tikzcd}
	(v,\alpha) \arrow{r}{I} \arrow[swap]{d}{\varphi} & (w,\beta) \arrow{d}{\varphi} & \text{in $G$}\\
	(v,\alpha_{[1,n-1]}) \arrow{r}{I}& (w,\beta_{[1,n-1]}) & \qquad \quad \text{in $G_{[1,n-1]}$}
\end{tikzcd}
\end{equation}
where $\alpha_{[1,n-1]}:=\varphi(\alpha)$ and $\beta_{[1,n-1]}:=\varphi(\beta)$. By Theorem~\ref{theorem.G1 and Gn} both $G$ and
$G_{[1,n-1]}$ satisfy axiom~\ref{axiom.labels}, which implies
\begin{equation}
\label{equation.length type changes}
\begin{aligned}
	\len(\beta) &= \len(\alpha) + d  &&\text{with $d\in\{0,\pm1\}$,}\\
	\len(\beta_{[1,n-1]}) &= \len(\alpha_{[1,n-1]}) + d' &&\text{with $d'\in \{0,\pm1\}$.}
\end{aligned}
\end{equation}
Note that $d=0,1,-1$ if the edge is of type~\ref{axiom.length pres}, \ref{axiom.length inc}, \ref{axiom.length dec} in axiom~\ref{axiom.labels}, respectively.
The same is true for $d'$. Furthermore, by definition we have
\begin{equation}
\label{equation.length changes}
\begin{aligned}
	\len(\alpha) &= \len(\alpha_{[1,n-1]}) + \delta_\alpha && \text{with $\delta_\alpha\in \{0,1\}$,}\\
	\len(\beta) &= \len(\beta_{[1,n-1]}) + \delta_\beta && \text{with $\delta_\beta \in \{0,1\}$.}
\end{aligned}
\end{equation}

\begin{proposition}
\label{proposition.G from Gn}
Let $G$ be a $\CS$-graph and consider edges in $G$ and $G_{[1,n-1]}$ as in~\eqref{equation.edges in G and Gn}, define $d,d',\delta_\alpha,\delta_\beta$ as 
in~\eqref{equation.length type changes} and~\eqref{equation.length changes}, and set $\ell=\len(\alpha_{[1,n-1]})$.
Then the following is true:
\begin{enumerate}
\item \label{i:d}
$d=d'$ except possibly when $d'=0$, $\delta_\alpha=0$ and $I^+=\alpha_{[1,n-1]}^{(\ell)}$ in which case $d=1$ is possible.
\item \label{i:delta}
$\delta_\beta=\delta_\alpha$ except possibly when $d'=0$, $\delta_\alpha=0$ and $I^+=\alpha_{[1,n-1]}^{(\ell)}$ in which case $\delta_\beta=1$ is possible.
\end{enumerate}
\end{proposition}

\begin{remark}
\label{remark.G from Gn}
Proposition~\ref{proposition.G from Gn} can be reformulated as follows.
If the edge $(v,\alpha_{[1,n-1]}) \stackrel{I}{\longrightarrow} (w,\beta_{[1,n-1]})$
in $G_{[1,n-1]}$ is of type~\ref{axiom.length pres}, \ref{axiom.length inc} or~\ref{axiom.length dec} in axiom~\ref{axiom.labels}, then the corresponding
edge in $G$ is of the same type unless possibly when (1) the edge is of type~\ref{axiom.length pres} in axiom~\ref{axiom.labels}, (2) $\delta_\alpha=0$,
and (3) $n-1\in I$; when (1), (2) and (3) are satisfied, the edge can become of type~\ref{axiom.length inc} in $G$.
\end{remark}

\begin{proof}[Proof of Proposition~\ref{proposition.G from Gn}]
We first show that~(\ref{i:d}) implies~(\ref{i:delta}). Note that
\[
	\delta_\beta = \len(\beta) - \len(\beta_{[1,n-1]}) = \len(\alpha) + d -\len(\alpha_{[1,n-1]}) - d'
	= d-d' + \delta_{\alpha}.
\]
Hence $d=d'$ implies $\delta_\beta=\delta_\alpha$. Furthermore in the exceptional case $\delta_\beta = d-d'+\delta_\alpha=1-0+0=1$.

Next we prove~(\ref{i:d}). Both $G$ and $G_{[1,n-1]}$ are $\CS$-graphs by assumption and Theorem~\ref{theorem.G1 and Gn}.
Hence axiom~\ref{axiom.labels} must hold for both graphs. In particular, if $j+1<\ell$ for the edge in $G_{[1,n-1]}$ or $j+1=\ell$ and 
$I^+ \subsetneq \alpha_{[1,n-1]}^{(\ell)}$, then by the properties of the map $\varphi$, the edge must be of the same type in $G$. This implies
$d=d'$.

Now assume that $j+1=\ell$ and $I^+=\alpha_{[1,n-1]}^{(\ell)}$ for the edge in $G_{[1,n-1]}$. In this case type~\ref{axiom.length inc} 
in axiom~\ref{axiom.labels} reduces to type~\ref{axiom.length pres}, so that $d'=1$ is not possible.
When $d'=-1$, the edge in $G_{[1,n-1]}$ is of type~\ref{axiom.length dec} and $\beta_{[1,n-1]} >_D \alpha_{[1,n-1]}$. Since 
$\len(\beta_{[1,n-1]}) < \len(\alpha_{[1,n-1]})$ we have that $\beta>_D \alpha$ independent of whether we $\beta$ is obtained from $\beta_{[1,n-1]}$
by adding one to its last part (or $\delta_\beta=0$) or by creating a new part of size 1 (or $\delta_\beta=1$) and same for $\alpha$. This implies that the 
edge in $G$ is of type~\ref{axiom.length dec} since edges of type~\ref{axiom.length pres} and~\ref{axiom.length inc} satisfy $\beta<_D\alpha$.
Hence $d=-1=d'$ proving the claim.

Finally assume that $I^+=\alpha_{[1,n-1]}^{(\ell)}$ and $d'=0$. If $d=0$, then $d=d'$ and the claim holds. If $d=-1$, then
$\delta_\beta=\delta_\alpha-1$, which implies that $\delta_\beta=0$ and $\delta_\alpha=1$. Hence
$\alpha_{[1,n-1]} = (\alpha^{(1)},\ldots,\alpha^{(\ell)})$ and $\alpha=(\alpha^{(1)},\ldots,\alpha^{(\ell)},\{n\})$.
Since $d'=0$, we have $\beta_{[1,n-1]} =  (\alpha^{(1)},\ldots,\alpha^{(\ell-1)} \setminus \{k\}, \alpha^{(\ell)} \cup \{k\})$ by 
axiom~\ref{axiom.labels}~\ref{axiom.length pres}.
Since $\delta_\beta=0$, we have $\beta= (\alpha^{(1)},\ldots,\alpha^{(\ell-1)} \setminus \{k\}, \alpha^{(\ell)} \cup \{k,n\})$.
However, from $d=-1$ using axiom~\ref{axiom.labels} \ref{axiom.length dec} we obtain $\beta^{(\ell)}=\{n\} \neq \alpha^{(\ell} \cup \{k,n\}$ 
from the expression for $\alpha$, yielding a contradiction. Hence $d=-1$ is not possible.
Finally, $d=1$ implies that $\delta_\beta=\delta_\alpha+1$, so that $\delta_\alpha=0$ and $\delta_\beta=1$, proving the claim.
\end{proof}

\begin{theorem}
\label{theorem.unique branching}
Let $G$ be a $\CS$-graph. Then $G$ is uniquely determined by $G_{[1,n-1]}$.
\end{theorem}

\begin{proof}
By axiom~\ref{axiom.crystal}, $G$ contains $B(\lambda)_s$ as a subgraph. The crystal $B(\lambda)_s$ has a unique vertex of weight $\lambda$ and
hence $G$ has a unique vertex $(u_\lambda,\lambda)$ with label $\lambda$. By Theorem~\ref{theorem.G1 and Gn}, $G_{[1,n-1]}$ is a union of $\CS$-graphs
$G_{\lambda^-}$, where $\lambda/\lambda^-$ is a skew shape with one box. Each $G_{\lambda^-}$ has a unique vertex 
$(u_{\lambda^-},\lambda^-)$ with label $\lambda^-$. Furthermore, labeling the edges in $G$ that correspond to edges in the crystal in $B(\lambda)_s$ by 
$f_j$ if the weight changes as in axiom~\ref{axiom.labels}~\ref{axiom.length pres}, we have
\begin{equation}
\label{equation.connected by crystal}	
	(u_{\lambda^-},\tilde{\lambda}^-) = f_{s-1} f_{s-2} \cdots f_r (u_\lambda,\lambda)
\end{equation}
for $r$ the row index of $\lambda/\lambda^-$ and $\tilde{\lambda}^-=(\lambda_1,\ldots,\lambda_r-1,\ldots,\lambda_{s-1},\lambda_s+1)$.
Hence $G_{[1,n-1]}$ together with the sequence of edges $f^{(r)} := f_{s-1} f_{s-2} \cdots f_r$ is connected. This means that every vertex 
$(v,\alpha_{[1,n-1]})$ in $G_{[1,n-1]}$ can be reached by a path
\begin{equation}
\label{equation.path}
	(u_\lambda,\lambda_{[1,n-1]}) \stackrel{I_1}{\longrightarrow} (v^{(1)},\alpha^{(1)}_{[1,n-1]}) \stackrel{I_2}{\longrightarrow} \cdots 
	\stackrel{I_p}{\longrightarrow} (v^{(p)},\alpha^{(p)}_{[1,n-1]})=(v,\alpha_{[1,n-1]}),
\end{equation}
where either $n\not \in I_j$ or the edge labeled $I_j$ is in $B(\lambda)_s$.

To show that $G$ is uniquely determined by $G_{[1,n-1]}$, we need to show that:
\begin{enumerate}
\item All edges in $G$ labeled by an interval $I$ with $n\in I$ can be recovered.
\item The vertex labels in $G$ can be recovered.
\end{enumerate}
We begin by proving (2) by induction. The vertex $(u_\lambda,\lambda)$ in $G$ has vertex label $\lambda$ and label 
$\varphi(\lambda)=\lambda_{[1,n-1]}$ in
$G_{[1,n-1]}$. Suppose by induction that it is known for all vertices along the path in~\eqref{equation.path} how to lift the vertex labels
$\alpha^{(j)}_{[1,n-1]}$ to $\alpha^{(j)}$ in $G$ for $1\leqslant j<p$. That is, we have the path
\begin{equation}
\label{equation.induction path}
	(u_\lambda,\lambda) \stackrel{I_1}{\longrightarrow} (v^{(1)},\alpha^{(1)}) \stackrel{I_2}{\longrightarrow} \cdots 
	\stackrel{I_{p-1}}{\longrightarrow} (v^{(p-1)},\alpha^{(p-1)})
\end{equation}
in $G$. If the edge labeled $I_p$ in~\eqref{equation.path} is in $B(\lambda)_s$, then the 
length of the label $\alpha$ in $G$ and $\alpha_{[1,n-1]}$ are the same and hence $(v,\alpha_{[1,n-1]})$ in $G_{[1,n-1]}$ lifts to $(v,\alpha)$ in $G$. 
Otherwise, by Remark~\ref{remark.G from Gn} the lift from $(v,\alpha_{[1,n-1]})$ in $G_{[1,n-1]}$ to $(v,\alpha)$ in $G$ is determined
except when the last edge in~\eqref{equation.path} is of type~\ref{axiom.length pres}, $n-1\in I_p$,
and $\delta_{\alpha^{(p-1)}}=0$. We will deal with this exceptional case using Lusztig involution (axiom~\ref{axiom.lusztig}~\ref{axiom.lusztig n}) 
and the fans (axiom~\ref{axiom.fans}).

By axiom~\ref{axiom.lusztig}~\ref{axiom.lusztig n}, the path in~\eqref{equation.induction path} implies a path
\[
	(u_{\rev(\lambda)},\rev(\lambda)) \stackrel{I^\LL_1}{\longleftarrow} (v^{(1)},\rev(\alpha^{(1)})) \stackrel{I^\LL_2}{\longleftarrow} \cdots 
	\stackrel{I^\LL_{p-1}}{\longleftarrow} (v^{(p-1)},\rev(\alpha^{(p-1)})),
\]
where $(u_{\rev(\lambda)},\rev(\lambda))$ is the unique vertex with label $\rev(\lambda)$ in $G$. We also know by Lusztig involution
that there is an incoming edge labeled $I_p^\LL$. By Remark~\ref{remark.G from Gn} the edge 
\[ (v^{(p)},\rev(\alpha)) \stackrel{I_p^\LL}{\longrightarrow}
(v^{(p-1)},\rev(\alpha^{(p-1)})) \] and vertex labels are determined unless $n-1\in I_p^\LL$ or equivalently $2\in I_p$.

Hence let us now assume that the edge labeled $I_p$ in~\eqref{equation.path} is of type~\ref{axiom.length pres}, $\delta_{\alpha^{(p-1)}}=0$,
and $2,n-1\in I_p$. This implies that either $I_p=[2,n-1]$ or $I_p=[1,n-1]$. We want to determine whether $\delta_\alpha=0$ or $1$.
The case $\delta_\alpha=1$ would require the existence of a fan by axiom~\ref{axiom.fans}. We consider the two cases separately.

\medskip

\noindent \textbf{Case 1:} $I_p=[2,n-1]$. Since $|I_p|\geqslant 3$ is odd by axiom~\ref{axiom.intervals}, we know that $n$ is odd.
First assume that $n\geqslant 7$. The case $\delta_\alpha=1$ would require the existence of a fan by axiom~\ref{axiom.fans}, in particular
\begin{center}
\raisebox{-0.8cm}{
\scalebox{0.8}{
\begin{tikzpicture}[>=latex,line join=bevel,]
\node (N0) at (0,0) [draw,draw=none] {$(v^{(p-1)},\alpha^{(p-1)})$};
\node (N1) at (-4,3) [draw,draw=none] {$(v,\alpha)=(v^{(p)},\alpha^{(p)})$};
\node (N2) at (4,3) [draw,draw=none] {$(w,\beta)$};
  \draw[->] (N0)  to node[opLabel]{$I_p=[2,n-1]$} (N1);
  \draw[->] (N0)  to node[opLabel]{$[3,n-2]$} (N2);
  \draw[->] (N1)  to node[opLabel]{$[n-2,n]$} (N2);
\end{tikzpicture}
}}
\end{center}
Under the Lusztig involution, this becomes
\begin{center}
\raisebox{-0.8cm}{
\scalebox{0.8}{
\begin{tikzpicture}[>=latex,line join=bevel,]
\node (N0) at (0,0) [draw,draw=none] {$(v^{(p-1)},\rev(\alpha^{(p-1)}))$};
\node (N1) at (-4,3) [draw,draw=none] {$(v,\rev(\alpha))=(v^{(p)},\rev(\alpha^{(p)}))$};
\node (N2) at (4,3) [draw,draw=none] {$(w,\rev(\beta))$};
  \draw[->] (N1)  to node[opLabel]{$[2,n-1]$} (N0);
  \draw[->] (N2)  to node[opLabel]{$[3,n-2]$} (N0);
  \draw[->] (N2)  to node[opLabel]{$[1,3]$} (N1);
\end{tikzpicture}
}}
\end{center}
Here the edges on the right and the top and the vertex labels are determined by Remark~\ref{remark.G from Gn}. 
Hence the existence of this configuration under Lusztig involution would determine that $\delta_\alpha=1$. If it does not exist,
this requires $\delta_\alpha=0$.

Similarly for $n=5$, the case $\delta_\alpha=1$ would require by axiom~\ref{axiom.fans} that
\begin{center}
\raisebox{-0.8cm}{
\scalebox{0.8}{
\begin{tikzpicture}[>=latex,line join=bevel,]
\node (N0) at (0,0) [draw,draw=none] {$(v^{(p-1)},\alpha^{(p-1)})$};
\node (N1) at (0,3) [draw,draw=none] {$(v,\alpha)=(v^{(p)},\alpha^{(p)})$};
  \draw[->] (N0)  to node[opLabel]{$I_p=[2,4]$} (N1);
  \draw[->] (N1) to[bend left=60] node[opLabel]{$[3,5]$} (N0);
\end{tikzpicture}
}}
\end{center}
which under the Lusztig involution becomes
\begin{center}
\raisebox{-0.8cm}{
\scalebox{0.8}{
\begin{tikzpicture}[>=latex,line join=bevel,]
\node (N0) at (0,0) [draw,draw=none] {$(v^{(p-1)},\rev(\alpha^{(p-1)}))$};
\node (N1) at (0,3) [draw,draw=none] {$(v,\rev(\alpha))=(v^{(p)},\rev(\alpha^{(p)}))$};
  \draw[->] (N1)  to node[opLabel]{$[2,4]$} (N0);
  \draw[->] (N0) to[bend right=60] node[opLabel]{$[1,3]$} (N1);
\end{tikzpicture}
}}
\end{center}
By Remark~\ref{remark.G from Gn}, the vertex labels are determined under the arrow labeled $[1,3]$.
Hence the existence of this configuration under Lusztig involution would determine that $\delta_\alpha=1$. If it does not exist,
this requires $\delta_\alpha=0$.

\medskip

\noindent \textbf{Case 2:} $I_p=[1,n-1]$. Since $|I_p|\geqslant 3$ is odd by axiom~\ref{axiom.intervals}, we know that $n$ is even.
First assume that $n\geqslant 6$. The case $\delta_\alpha=1$ would require the existence of a fan by axiom~\ref{axiom.fans}, in particular
\begin{center}
\raisebox{-0.8cm}{
\scalebox{0.8}{
\begin{tikzpicture}[>=latex,line join=bevel,]
\node (N0) at (0,0) [draw,draw=none] {$(v^{(p-1)},\alpha^{(p-1)})$};
\node (N1) at (-4,3) [draw,draw=none] {$(v,\alpha)=(v^{(p)},\alpha^{(p)})$};
\node (N2) at (4,3) [draw,draw=none] {$(w,\beta)$};
  \draw[->] (N0)  to node[opLabel]{$I_p=[1,n-1]$} (N1);
  \draw[->] (N0)  to node[opLabel]{$[2,n-2]$} (N2);
  \draw[->] (N1)  to node[opLabel]{$[n-2,n]$} (N2);
\end{tikzpicture}
}}
\end{center}
Under the Lusztig involution, this becomes
\begin{center}
\raisebox{-0.8cm}{
\scalebox{0.8}{
\begin{tikzpicture}[>=latex,line join=bevel,]
\node (N0) at (0,0) [draw,draw=none] {$(v^{(p-1)},\rev(\alpha^{(p-1)}))$};
\node (N1) at (-4,3) [draw,draw=none] {$(v,\rev(\alpha))=(v^{(p)},\rev(\alpha^{(p)}))$};
\node (N2) at (4,3) [draw,draw=none] {$(w,\rev(\beta))$};
  \draw[->] (N1)  to node[opLabel]{$[2,n]$} (N0);
  \draw[->] (N2)  to node[opLabel]{$[3,n-1]$} (N0);
  \draw[->] (N2)  to node[opLabel]{$[1,3]$} (N1);
\end{tikzpicture}
}}
\end{center}
The edge labeled $[2,n-2]$ in the first picture under Remark~\ref{remark.G from Gn} determines the label $\beta$ in $G$, and the edge 
labeled $[1,3]$ in the second picture under Remark~\ref{remark.G from Gn} determines the label $\rev(\alpha)$ and hence $\alpha$ in $G$.
If these edges do not exist, then $\delta_\alpha=0$.

For $n=4$, by inspecting axiom~\ref{axiom.labels}, the case $\delta_\alpha=1$ can only happen if $\alpha^{(p-1)}=(2,2)$. By the fan 
axiom~\ref{axiom.fans}, $\delta_\alpha=1$ requires
\begin{center}
\raisebox{-0.8cm}{
\scalebox{0.8}{
\begin{tikzpicture}[>=latex,line join=bevel,]
\node (N0) at (0,0) [draw,draw=none] {$(v^{(p-1)},(2,2))$};
\node (N1) at (0,3) [draw,draw=none] {$(v,\alpha)=(v^{(p)},(1,2,1))$};
  \draw[->] (N0)  to node[opLabel]{$I_p=[1,3]$} (N1);
  \draw[->] (N1) to[bend left=60] node[opLabel]{$[2,4]$} (N0);
\end{tikzpicture}
}}
\end{center}
which is fixed under the Lusztig involution. By axiom~\ref{axiom.incoming}, there has to be an incoming edge into $(v^{(p-1)},(2,2))$
labeled $[2,4]$. For $\delta_\alpha=0$, this incoming edge has to come from a vertex with label $(3,1)$. Thus we can decide whether
$\delta_\alpha=0,1$ depending on whether $G$ has $3$ or $2$ vertices, respectively.

Hence we have recovered all vertex labels in $G$ from $G_{[1,n-1]}$, proving (2).

Since $G_{[1,n-1]}$ together with the crystal edges is connected and since all vertex labels in $G$ are known, the 
edges labeled $I$ with $n\in I$ can be recovered using the Lusztig involution in axiom~\ref{axiom.lusztig} \ref{axiom.lusztig n} unless $I=[1,n]$.

The case $I=[1,n]$ can only happen if $n$ is odd and for $(v,\alpha) \stackrel{I}{\longrightarrow} (w,\beta)$ we have $\alpha=(\frac{n+1}{2},\frac{n-1}{2})$.
This edge has to be of type~\ref{axiom.length pres}. If it was of type~\ref{axiom.length dec}, $\alpha$ would have to have at least three parts.
If it was of type~\ref{axiom.length inc}, by axiom~\ref{axiom.fans} there would have to be an edge out of $(w,\beta)$ with interval $[n-1,n+1]$, which is
not possible since intervals are subsets of $[n]$. By axiom~\ref{axiom.outgoing}~\ref{axiom.outgoing.incoming} there cannot be an incoming edge
into $(v,\alpha)$ from a vertex labeled $(\frac{n+3}{2},\frac{n-3}{2})$. Hence $\lambda$ of axiom~\ref{axiom.crystal} is
$\lambda=(\frac{n+1}{2},\frac{n-1}{2})$ in this case, which is a unique label. This implies by axiom~\ref{axiom.lusztig} \ref{axiom.lusztig n} that the label
$(\frac{n-1}{2},\frac{n+1}{2})$ is also unique and hence the edge labeled $I=[1,n]$ is uniquely determined.
\end{proof}

\begin{corollary}
\label{corollary.uniqueness}
    Suppose $G_\lambda$ is a $\CS$-graph. Then $G_\lambda$ is isomorphic to the crystal skeleton $\CS(\lambda)$.
\end{corollary}

\begin{proof}
By Theorem~\ref{theorem.CS is CS-graph}, $\CS(\lambda)$ is a $\CS$-graph. By Theorem~\ref{theorem.unique branching}, $G_\lambda$
is uniquely determined by $(G_\lambda)_{[1,n-1]}$. By Theorems~\ref{theorem.G1 and Gn} and~\ref{theorem.branching}, $G_\lambda$ and $\CS(\lambda)$
have the same branching properties. Hence we must have $G_\lambda \cong \CS(\lambda)$ by induction.
\end{proof}

\begin{corollary}
    Let $G_\lambda$ be a $\CS$-graph. Then $G_\lambda$ has a unique vertex labeled $\lambda$ and all other vertex labels $\alpha$ satisfy 
    $\alpha \leqslant_D \lambda$. That is, $\lambda$ is the unique largest label in dominance order.
\end{corollary}

\begin{proof}
By Corollary~\ref{corollary.uniqueness}, $G_\lambda$ is isomorphic to the crystal skeleton $\CS(\lambda)$. By Lemma~\ref{lemma.QY}, the 
descent composition associated to the vertex $T$ in $\CS(\lambda)$ is the weight $\alpha$ of the quasi-Yamanouchi tableau $Q$ in the same quasi-crystal
component as $T$. Recall that the descent compositions are precisely the labels of the vertices in $G_\lambda$ by Theorem~\ref{theorem.CS is CS-graph}.
In a crystal graph $B(\lambda)$, the weights of all vertices are smaller in dominance order than $\lambda$ (see for 
example~\cite[Chapter 4.4]{BumpSchilling.2017}). Since $Q$ is an element in the crystal graph $B(\lambda)$ underlying the crystal skeleton 
$\CS(\lambda)$, we hence have $\alpha \leqslant_D \lambda$.
\end{proof}

\subsection{$S_n$-axioms}
\label{section.axiom Sn}

We will now state $S_n$ analogues of the $\CS$-graph axioms.

\begin{axiom}[$S_n$-axioms for crystal skeletons]
\label{Sn axioms}
Let $n$ be a positive integer, $\lambda \vdash n$, and $G_\lambda$ be a finite, connected, vertex- and edge-labeled graph with labeled vertex set 
$V_L$ and labeled edge set $E_L$ as in Section~\ref{section.GLn axioms}. The $S_n$-axioms for $G_\lambda$ require:

\begin{enumerate}[label=\bf {S}\arabic{enumi}.,ref=\bf {S}\arabic{enumi}]\setcounter{enumi}{-1}
\item
\label{axiom.Sn intervals} \emph{(Intervals)} Same as axiom~\ref{axiom.intervals}.

\item 
\label{axiom.Sn outgoing}(\emph{Outgoing edges}) Same as axiom~\ref{axiom.outgoing}.

\item\label{axiom.Sn labels} (\emph{Labels}) Same as axiom~\ref{axiom.labels}.

\item\label{axiom.Sn fans} (\emph{Fans}) Same as axiom~\ref{axiom.fans}.

\item\label{axiom.Sn lusztig} (\emph{Lusztig involution}) We have $\mathcal{L}_n(G_\lambda) \cong G_\lambda$.

\item\label{axiom.Sn branching} (\emph{Branching}) We have
    \[
    	(G_\lambda)_{[1,n-1]} \cong \bigcup_{\substack{\lambda^-\\ |\lambda/\lambda^-|=1}} G_{\lambda^-},
    \]
    where $G_{\lambda^-}$ satisfy Axiom~\ref{Sn axioms}.

\item\label{axiom.Sn connectivity} (\emph{Connectivity}) The graph $G_\lambda$ has a unique vertex $(u_\lambda,\lambda)$ with label $\lambda$ and 
all other labels $\alpha$ satisfy $\alpha\leqslant_D \lambda$. Denote by $f_j$ the edge as in axiom~\ref{axiom.labels}~\ref{axiom.length pres}
and set $s=\len(\lambda)$. Then for each $\lambda^-$ such that $|\lambda/\lambda^-|=1$, we have
\[
	(u_{\lambda^-},\tilde{\lambda}^-) = f^{(r)}(u_\lambda,\lambda),
\]
where $r$ is the row index of $\lambda/\lambda^-$, $\tilde{\lambda}^-=(\lambda_1,\ldots,\lambda_r-1,\ldots,\lambda_{s-1},\lambda_s+1)$, and
$f^{(r)} := f_{s-1} f_{s-2} \cdots f_r$. Furthermore, $(G_\lambda)_{[1,n-1]}$ together with the sequence of edges $f^{(r)}$ is connected. 
\end{enumerate}
\end{axiom}

\begin{theorem}
A graph $G_\lambda$ satisfying Axiom~\ref{Sn axioms} is a $\CS$-graph (satisfying Axiom~\ref{axioms}).
Conversely, a $\CS$-graph satisfies Axiom~\ref{Sn axioms}.
\end{theorem}

\begin{proof}
Suppose $G_\lambda$ satisfies Axiom~\ref{Sn axioms}. Axioms~\ref{axiom.Sn intervals}-\ref{axiom.Sn fans} imply
axioms~\ref{axiom.intervals}-\ref{axiom.fans}. Axiom~\ref{axiom.Sn lusztig} implies Axiom~\ref{axiom.lusztig}\ref{axiom.lusztig n}.
Axiom~\ref{axiom.Sn branching} stipulates that each connected component $G_{\lambda^-}$ satisfies Axiom~\ref{Sn axioms}
for $\lambda^- \vdash n-1$, and hence in particular axiom~\ref{axiom.Sn lusztig} with $n$ replaced by $n-1$. This implies
axiom~\ref{axiom.lusztig}\ref{axiom.lusztig n-1}. 

We show that axiom~\ref{axiom.crystal} holds by induction on $n$. For $n=1$, $G_{(1)}$ is a 
single vertex as is $B((1))_1$. Hence axiom~\ref{axiom.crystal} holds. By the induction hypothesis, the connected components of
$(G_\lambda)_{[1,n-1]}$ in~\ref{axiom.Sn branching} satisfy axiom~\ref{axiom.crystal}. Thus all $G_{\lambda^-}$ in~\ref{axiom.Sn branching}
are $\CS$-graphs satisfying Axiom~\ref{axioms}. Analyzing the proof of Theorem~\ref{theorem.unique branching}, only 
axioms~\ref{axiom.intervals}-\ref{axiom.lusztig} and the properties in~\ref{axiom.Sn connectivity} are used. This implies that $G_\lambda$
is uniquely determined from the $\CS$-graphs in $(G_\lambda)_{[1,n-1]}$ in~\ref{axiom.Sn branching}. This proves that $G_\lambda$ is
a $\CS$-graph satisfying Axiom~\ref{axioms}.

For the converse, axioms~\ref{axiom.intervals}-\ref{axiom.lusztig} imply axioms~\ref{axiom.Sn intervals}-\ref{axiom.Sn lusztig}.
Axiom~\ref{axiom.crystal} implies Axiom~\ref{axiom.Sn connectivity}. The unique label $(u_\lambda,\lambda)$ is the highest weight
vertex in $G_s \cong B(\lambda)_s$. Connectivity follows by the proof of Theorem~\ref{theorem.unique branching}.
Axiom~\ref{axiom.Sn branching} holds by Theorem~\ref{theorem.G1 and Gn}.
\end{proof}

\subsection{Local axioms}
\label{section.axiom local}

As outlined in Section~\ref{section.commutation}, Stembridge~\cite{Stembridge.2003} provided a local characterization of crystal
bases associated to representations of simply-laced root systems, in particular type $A$. Here we provide analogous
local axioms (Axiom~\ref{local axioms}) for $\CS$-graphs.

Let $G$ be a finite, connected, vertex- and edge-labeled graph with labeled vertex set $V_L$ and labeled edge set $E_L$ as in 
Section~\ref{section.GLn axioms}. Recall that $G_s$ with $s=\min\{\len(\alpha) \mid (v,\alpha) \in V_L\}$ is the subgraph of $G$ 
containing only the vertices $(v,\alpha)\in V_L$ with $\len(\alpha)=s$. 

To state our local axioms, we begin by defining the \defn{crystal operators} $f_j(v,\alpha)$ and $e_i(v,\alpha)$ for $(v,\alpha)\in G_s$ and their 
\defn{string lengths} $\varphi_j(v,\alpha)$ and $\varepsilon_i(v,\alpha)$ (see~\eqref{equation.string lengths}). As the notation suggests, these maps
will be used to prove that $G_s$ is a crystal when Axiom~\ref{local axioms} is satisfied.

\begin{definition}\label{definition.fandphi}
Let $G$ be a finite, connected, vertex- and edge-labeled graph with labeled vertex set 
$V_L$ and labeled edge set $E_L$ as in Section~\ref{section.GLn axioms}. Suppose $G$ satisfies 
axioms~\ref{axiom.intervals},~\ref{ax:interval types-in},~\ref{axiom.outgoing},~\ref{axiom.incoming},~\ref{axiom.labels}, and~\ref{axiom.labels incoming}
and the relations in Theorem~\ref{theorem.commutation} and their duals for incoming edges. Then define the following:
\begin{enumerate}
    \item If there is no outgoing edge $(v,\alpha)\stackrel{I}{\longrightarrow} (w,\beta)$ in $G_s$ with $I \subseteq \alpha^{(j)} \cup \alpha^{(j+1)}$,
    then $f_j(v,\alpha):=\emptyset$, $\varphi_j(v,\alpha):=0$. Similarly, if there is no incoming edge 
    $(v,\alpha)\stackrel{I}{\longrightarrow} (w,\beta)$ in $G_s$ with $I \subseteq \beta^{(j)} \cup \beta^{(j+1)}$, then
    $e_j(w,\beta) := \emptyset$ and $\varepsilon_j(w,\beta):=0$.
    \item
    Otherwise, if there is an edge $(v,\alpha)\stackrel{I}{\longrightarrow} (w,\beta)$ in $G_s$ with $I \subseteq \alpha^{(j)} \cup \alpha^{(j+1)}$, define
    \begin{equation}
    \label{equation.f I}
    f_j(v,\alpha) :=  (w, \beta), \qquad e_j(w,\beta) = (v,\alpha),
   \end{equation}
   and
\begin{equation}
\label{equation.varphi I}
	\varphi_j(v,\alpha) := 	\min I - \min \alpha^{(j)} +1, \qquad
	\varepsilon_j(w,\beta) := \max \beta^{(j+1)} - \max I.
\end{equation}
\end{enumerate}
\end{definition}

We first check that $f_j$ is well-defined, and prove that whenever $(v,\alpha)\stackrel{I}{\longrightarrow} (w,\beta)$ is in $G_s$, 
there are certain edges in $G$ as well.

\begin{lemma}
\label{lemma.strings.in.G}
Let $G$ be a graph satisfying the hypotheses in Definition \ref{definition.fandphi}, with an edge $(v,\alpha)\stackrel{I}{\longrightarrow} (w,\beta)$
for $I=[i,i+2m]$ in $G$. Then
\begin{enumerate}
    \item If $I$ is length-preserving, $f_j(v,\alpha)$ is well-defined;
    \item  If $f_j^{k'}(v,\alpha) \in G_s$ for all $0\leqslant k' \leqslant k$ for some $0 \leqslant k < \varphi_j(v,\alpha)$, then there is an edge in $G$
    \begin{equation}
    \label{equation.f_j^k}
    	f_j^k(v,\alpha) \stackrel{I_{k}}{\longrightarrow} f_j^{k+1}(v,\alpha) 
    \end{equation}
   with $I_k = [i-k,i-k+2m]$. For $k=\varphi_j(v,\alpha)$ such an edge does not exist.
\end{enumerate}
Analogous statements hold for $e_j(w,\beta)$.
\end{lemma}

\begin{proof}
By assumption the edge $(v,\alpha)\stackrel{I}{\longrightarrow} (w,\beta)$ exists in $G$. Hence
$I^- \cup \{i+m\} \subseteq \alpha^{(j)}$ and $I^+ \subseteq \alpha^{(j+1)}$ for some $j$ by axiom~\ref{axiom.intervals}.
Observe the following:
\begin{enumerate}
\item
 If $I$ is length-preserving, we have $1\leqslant j<s$. Inspecting Theorem~\ref{theorem.commutation} case~\ref{case3}, 
 if  $(v,\alpha)\stackrel{I}{\longrightarrow} (w,\beta)$ in $G_s$, then $I$ must be maximal among all intervals 
 $J \subseteq \alpha^{(j)} \cup \alpha^{(j+1)}$ since if $I \subsetneq J$, then $I$ is length-increasing and thus not an edge in $G_s$. 
 Hence $f_j$ is well-defined.

\item The case $\varphi_j(v,\alpha)=0$ is vacuously true. Hence assume $\varphi_j(v,\alpha)\geqslant 1$.
By assumption, $f_j^k(v,\alpha) \in G_s$ with $k\geqslant 0$.
Furthermore, by assumption all edges $I_0,\ldots,I_{k-1}$ are in $G_s$ and hence length-preserving.
By axiom~\ref{axiom.labels}, this implies that the label of the vertex $f_j^k(v,\alpha)$ is
$\alpha-k \epsilon_j + k \epsilon_{j+1}$, where $\epsilon_j$ is the unit vector in $\mathbb{Z}^s$ with a 1 in position $j$.
Note that 
\[ 
	I_k^- \cup \{i-k+m\} \subseteq \bigg(\alpha-k \epsilon_j + k \epsilon_{j+1}\bigg)^{(j)} \quad \quad \textrm{and} 
	\quad \quad I_k^+ \subseteq \bigg(\alpha-k \epsilon_j + k \epsilon_{j+1}\bigg)^{(j+1)}
\]
since by assumption $0 \leqslant k<\varphi_j(v,\alpha) = \min I - \min \alpha^{(j)} +1=i-\min \alpha^{(j)}+1$.
Hence by axiom~\ref{axiom.outgoing} either the edge in~\eqref{equation.f_j^k} exists and the claim holds or there is an edge
$(v',\alpha') \stackrel{J}{\longrightarrow} f_j^k(v,\alpha)$ such that $J \subsetneq I_k$ and $\alpha' \geqslant_D \alpha-k\epsilon_j
+k\epsilon_{j+1}$. Suppose such an edge $J$ exists. By the dual of case~\ref{case3} of Theorem~\ref{theorem.commutation}, this implies 
that $J$ is length-decreasing which is not possible since $f_j^k(v,\alpha) \in G_s$ and hence has a label with minimal length.

Finally, if $k=\varphi_j(v,\alpha)$, we have $I_k \not \subseteq \alpha^{(j)} \cup \alpha^{(j+1)}$. Hence the edge in~\eqref{equation.f_j^k}
does not exist.
 \end{enumerate}
\end{proof}

We now give the set of local axioms for crystal skeletons.

\begin{axiom}[Local axioms for crystal skeletons]
\label{local axioms}
Let $G$ be a finite, connected, vertex- and edge-labeled graph with labeled vertex set 
$V_L$ and labeled edge set $E_L$ as in Section~\ref{section.GLn axioms}. The local axioms for $G$ require:

\begin{enumerate}[label=\bf {L}\arabic{enumi}.,ref=\bf {L}\arabic{enumi}]\setcounter{enumi}{-1}
\item
\label{axiom.local intervals} \emph{(Intervals)} Same as axioms~\ref{axiom.intervals} and~\ref{ax:interval types-in}.

\item 
\label{axiom.local outgoing}(\emph{Outgoing and incoming edges}) Same as axioms~\ref{axiom.outgoing} and~\ref{axiom.incoming}.

\item\label{axiom.local labels} (\emph{Labels}) Same as axioms~\ref{axiom.labels} and~\ref{axiom.labels incoming}.

\item\label{axiom.local commutations} (\emph{Commutation relations}) 
\begin{enumerate}
\item
If there is a vertex $(v,\alpha)\in V_L$ with two distinct outgoing edges labeled
$I$ and $J$, respectively, then commutation relations as in Theorem~\ref{theorem.commutation} hold.
\item
If there is a vertex $(v,\alpha)\in V_L$ with two distinct incoming edges labeled
$I$ and $J$, respectively, then the dual to the commutation relations as in Theorem~\ref{theorem.commutation} hold 
(see Remark~\ref{remark.commutations duality}).
\item
For any edge $(v,\alpha) \stackrel{I}{\longrightarrow} (w,\beta)$ with $|I|=3$ and $I$ length-increasing
(resp. $(w,\beta) \stackrel{J}{\longrightarrow} (v,\alpha)$ with $|J|=3$ and $J$ length-decreasing), we have
the following subgraph with $J=\{i+1 \mid i \in I\}$ (resp. $I=\{j-1 \mid j \in J\}$)
\begin{equation}
\label{equation.loop}
\raisebox{-0.5cm}{
\scalebox{0.8}{
\begin{tikzpicture}[>=latex,line join=bevel,]
\node (N0) at (0,0) [draw,draw=none] {$(v,\alpha)$};
\node (N1) at (3,0) [draw,draw=none] {$(w,\beta)$};
  \draw[->,thick] (N0)  to node[opLabel]{$I$} (N1);
  \draw[->,thick] (N1) to[bend left=60] node[opLabel]{$J$} (N0);
\end{tikzpicture}
}}.
\end{equation}
\end{enumerate}
\item\label{axiom.local string} (\emph{Strings in $G_s$}) 
\begin{enumerate}
\item
Suppose $(v,\alpha) \stackrel{I}{\longrightarrow} f_j(v,\alpha)$ is an edge in 
$G_s$ with $I=[i,i+2m] \subseteq \alpha^{(j)}\cup \alpha^{(j+1)}$. Then for $0\leqslant k < \varphi_j(v,\alpha)$, the edge 
\[ 
	f_j^k(v,\alpha) \stackrel{I_k}{\longrightarrow} f_j^{k+1}(v,\alpha)
\] 
with $I_k = [i-k,i-k+2m]$ (in $G$ by Lemma~\ref{lemma.strings.in.G}) is length-preserving.
\item Suppose $e_j(v,\alpha) \stackrel{I}{\longrightarrow} (v,\alpha)$ is an edge in 
$G_s$ with $I=[i,i+2m] \subseteq \alpha^{(j)}\cup \alpha^{(j+1)}$. Then for $0\leqslant k < \varepsilon_j(v,\alpha)$, the edge 
\[ 
	e_j^{k+1}(v,\alpha) \stackrel{J_k}{\longrightarrow} e_j^{k}(v,\alpha)
\] 
with $J_k = [i+k,i+k+2m]$ (in $G$ by Lemma~\ref{lemma.strings.in.G}) is length-preserving.
\end{enumerate}
\item\label{axiom.local edge types} (\emph{Edge types}) 
Let $\lambda$ be a maximal label in dominance order among all vertex labels. Let $(v,\alpha) \in V_L$. Then the following is true:
\begin{enumerate}
\item If $\alpha \neq \lambda$, then $(v,\alpha)$ has at least one incoming edge which is length-preserving or length-increasing.
\item If $\rev(\alpha) \neq \lambda$, then $(v,\alpha)$ has at least one outgoing edge which is length-preserving or length-decreasing.
\end{enumerate}
\end{enumerate}
\end{axiom}

 We first show that Axiom~\ref{local axioms} implies axiom~\ref{axiom.crystal}.

\begin{lemma}
\label{lemma.top crystal}
Axiom~\ref{local axioms} implies axiom~\ref{axiom.crystal} for $G_s$.
\end{lemma}

\begin{proof}
We will prove that $G_s$ is isomorphic to a crystal $B(\lambda)_s$ for some $\lambda$ by showing that $G_s$ satisfies the Stembridge 
axioms~\cite{Stembridge.2003}. By Axiom~\ref{axiom.local edge types}, the highest weight $\lambda$ is unique since by (a) any other maximal
label in dominance order of the same length would have to have a length-preserving incoming edge and would hence not be maximal.
By~\cite[Theorem~4.12]{BumpSchilling.2017}, this implies the isomorphism.

Suppose that $(v,\alpha)\stackrel{I}{\longrightarrow} (w,\beta)$ and $(v',\alpha') \stackrel{J}{\longrightarrow} (w,\beta)$ are in $G_s$
with $J \subsetneq I$. Then $J$ is length-decreasing by the dual of case~\ref{case3} of Theorem~\ref{theorem.commutation} which holds
by axiom~\ref{axiom.local commutations}. But since $(v',\alpha')\in G_s$ this is not possible. Hence the partial inverse of $f_j$ 
(called $e_j$) is well-defined.

Next we check that the Stembridge commutation relations~\eqref{equation.stembridge square} and~\eqref{equation.stembridge octagon} hold.
The corresponding relations for the raising operators $e_i$ can be checked similarly.
We consider a vertex $(v,\alpha)$ in $G_s$ such that there are edges in $G_s$
\[ 
	(v,\alpha) \stackrel{I}{\longrightarrow} f_k(v,\alpha) \quad \textrm{with} \quad I\subseteq \alpha^{(k)} \cup \alpha^{(k+1)} \quad \quad \textrm{and} 
	\quad \quad (v,\alpha) \stackrel{J}{\longrightarrow} f_\ell(v,\alpha) \quad \textrm{ with } \quad J\subseteq \alpha^{(\ell)} \cup \alpha^{(\ell+1)}. 
\]  
In particular, both $I$ and $J$ are length-preserving.
By axiom~\ref{axiom.local commutations}, the commutation relations in Theorem~\ref{theorem.commutation} must hold.
Inspection shows that case~\ref{case3} cannot apply since this would require $J$ to be length-increasing.
Cases~\ref{case2aii} and~\ref{case2bii} contain length-decreasing edges, which cannot happen since $G_s$ contains vertices with 
minimal length labels. Hence cases~\ref{case1}, \ref{case2ai}, \ref{case2bi}, or~\ref{case2biii} must apply, which are squares or octagons as 
in the Stembridge axioms. 

In cases~\ref{case1}, \ref{case2ai}, and~\ref{case2biii}, if $I$ and $J$ are length-preserving, then all edges in the square are length-preserving.
For the octagon in case~\ref{case2bi}, length-decreasing edges cannot appear since we are in $G_s$. If $J'$ is length-increasing,
case~\ref{case2ai} also applies, which is a square with all length-preserving edges. Hence either case~\ref{case2bi} or case~\ref{case2ai}
occurs with all edges length-preserving.

It remains to check that the conditions on the string lengths are satisfied. 
By Lemma~\ref{lemma.strings.in.G} and axiom~\ref{axiom.local string}, we have that $\varphi_j(v,\alpha)$ in~\eqref{equation.varphi I} 
satisfies $\varphi_j(v,\alpha) = \max\{k \in \mathbb{Z}_{\geqslant 0} \mid f_j^k(v,\alpha)\neq \emptyset\}$ as in~\eqref{equation.string lengths}.
Hence to check the Stembridge axioms, it remains to check that the conditions on the string length are satisfied, namely
\[
	\varphi_k(f_\ell(v,\alpha)) = \varphi_k(v,\alpha) + \{0,1\} \quad \text{and} \quad
	\varphi_\ell(f_k(v,\alpha)) = \varphi_\ell(v,\alpha) + \{0,1\},
\]
and
\[
	\varphi_k(f_\ell(v,\alpha)) = \varphi_k(v,\alpha) \quad \text{or} \quad 
	\varphi_\ell(f_k(v,\alpha)) = \varphi_\ell(v,\alpha) \quad \text{for squares,}
\]
and
\begin{equation}
\label{equation.octagon CS}
	\varphi_k(f_\ell(v,\alpha)) = \varphi_k(v,\alpha)+1 \quad \text{and} \quad 
	\varphi_\ell(f_k(v,\alpha)) = \varphi_\ell(v,\alpha)+1 \quad \text{for octagons.}
\end{equation}
In case~\ref{case1} we have $\varphi_k(f_\ell(v,\alpha)) = \varphi_k(v,\alpha)$ and $\varphi_\ell(f_k(v,\alpha)) = \varphi_\ell(v,\alpha)$
using~\eqref{equation.varphi I} since $I \cap J = \emptyset$. Hence the square applies. For case~\ref{case2} we have 
$I \subseteq \alpha^{(k)} \cup \alpha^{(k+1)}$ and $J \subseteq \alpha^{(k+1)}\cup \alpha^{(k+2)}$, so that $\ell=k+1$. In case~\ref{case2ai}, 
it can be explicitly checked using~\eqref{equation.varphi I} that
$\varphi_{k+1}(f_k(v,\alpha)) = \varphi_{k+1}(v,\alpha)$ and $\varphi_k(f_{k+1}(v,\alpha)) = \varphi_k(v,\alpha)+1$. Hence the square applies.
In case~\ref{case2biii}, it can be explicitly checked that 
$\varphi_{k+1}(f_k(v,\alpha)) = \varphi_{k+1}(v,\alpha)+1$ and $\varphi_k(f_{k+1}(v,\alpha)) = \varphi_k(v,\alpha)$. Hence the square applies.
Finally, in case~\ref{case2bi} the relations~\eqref{equation.octagon CS} hold, hence the octagon applies.
\end{proof}

We next prove that if a graph satisfies Axiom~\ref{local axioms}, then the graph must be invariant under Lusztig involution. 

\begin{lemma}
\label{lemma.lusztig local}
Axiom~\ref{local axioms} implies axiom~\ref{axiom.lusztig}~\ref{axiom.lusztig n}.
\end{lemma}

\begin{proof}
Since axiom~\ref{axiom.crystal} holds for $G$ by Lemma \ref{lemma.top crystal}, $G_s$ is a crystal, which enjoys $\LL_n(G_s) \cong G_s$. 
Let $G_{\leqslant \ell}$ (resp. $G_{=\ell}$) be the induced subgraph of $G$ with all vertices $(v,\alpha)$ such that $\len(\alpha)\leqslant \ell$
(resp. $\len(\alpha)=\ell$).
Assume by induction that $\LL_n(G_{\leqslant \ell}) = G_{\leqslant \ell}$; the statement is true for $\ell=s$ by Lemma \ref{lemma.top crystal}. We will show that 
$\LL_n(G_{\leqslant \ell}) = G_{\leqslant \ell}$ implies that $\LL_n(G_{\leqslant \ell+1}) = G_{\leqslant \ell+1}$ as well.
The proof will be done in multiple steps:
\newcommand{\Gc}{{G^{\mathsf{con}}}}
\begin{enumerate}
    \item[\emph{Step 1}.] Let $\Gc$ be the subgraph of $G$ consisting of $G_{\leqslant \ell}$ and all vertices $(v,\alpha)$ with $\len(\alpha)=\ell+1$ which 
    are connected by an incoming or outgoing edge to a vertex in $G_{\leqslant \ell}$. We show that $\LL_n(\Gc)=\Gc$. \medskip
    \item[\emph{Step 2}.] Suppose 
    \[ G_{\leqslant \ell} \subsetneq G' \subseteq G_{\leqslant \ell+1} \]
    is a subgraph of $G$ such that $\LL_n(G')=G'$. By step 1, $G'=\Gc$ has this property.
    We successively extend $G'$ by adding length-preserving edges entering or leaving the vertices $(v,\alpha)\in G_{=\ell+1} \cap G'$ in a specific order.
    We continue adding edges and vertices to $G'$ until we have $G' = G_{\leqslant \ell + 1}$. By axiom~\ref{axiom.local edge types}, all
    vertices in $G_{\leqslant \ell+1}$ can be reached.
 \end{enumerate}

\noindent \emph{Proof of Step 1}.
Let $(v,\alpha) \stackrel{I}{\longrightarrow} (w,\beta)$ be a length-increasing edge in $G$ with $\len(\alpha)=\ell$.
Let $(v',\rev(\alpha))$ be the image of $(v,\alpha)$ under $\LL_n$. Since $(v,\alpha) \in G_{\leqslant \ell}$ and 
$\LL_n(G_{\leqslant \ell}) = G_{\leqslant \ell}$, we have that $(v',\rev(\alpha)) \in G_{\leqslant \ell} $. 

By axiom~\ref{axiom.incoming}, there is an edge in $G$
\begin{equation} \label{eq:local axiom proof edge}
(w',\beta') \stackrel{I^\LL}{\longrightarrow} (v',\rev(\alpha)) 
\end{equation}
with $I^\LL=[n+1-b,n+1-a]$ if $I=[a,b]$. 

Since $(v',\rev(\alpha)) \in G_{\leqslant \ell}$, by definition the segment in \eqref{eq:local axiom proof edge} is in $\Gc$. Since 
$\LL_n(G_{\leqslant \ell}) =G_{\leqslant \ell}$ and $I$ is length-increasing, the vertex $(w',\beta')$ cannot 
be in $G_{\leqslant \ell}$. This implies that the edge labeled $I^\LL$ in \eqref{eq:local axiom proof edge} is length-decreasing and hence by 
axioms~\ref{axiom.labels} and~\ref{axiom.labels incoming} we have $\beta'=\rev(\beta)$. Thus \eqref{eq:local axiom proof edge} can be rewritten as 
\[
	(w',\rev(\beta)) \stackrel{I^\LL}{\longrightarrow} (v',\rev(\alpha)),
\]
which is precisely the image of $(v,\alpha) \stackrel{I}{\longrightarrow} (w,\beta)$ under $\LL_n$ by definition of $\LL_n$.

Conversely, if $(w,\beta) \stackrel{I}{\longrightarrow} (v,\alpha)$ is a length-decreasing edge in $G$ with $\len(\alpha)=\ell$, then by very similar 
arguments the image of this interval under $\LL_n$ is in $\Gc$. Hence $\LL_n(\Gc)=\Gc$. \medskip

\noindent \emph{Proof of Step 2.}
Suppose $G_{= \ell +1}$ contains a length-preserving edge $(v,\alpha) \stackrel{K}{\longrightarrow} (w,\beta)$.
 We will show the following:
\begin{enumerate}
    \item  
    Suppose $(v,\alpha)\in G'$. Then $\LL_n$ maps  $(v,\alpha)$ to $(v',\rev(\alpha)) \in G'$. 
    We prove there is an edge in $G$
    \begin{equation}\label{eq:step2} (w',\beta') \stackrel{K^\LL}{\longrightarrow} (v',\rev(\alpha)) \end{equation}
    that is length-preserving. By axioms~\ref{axiom.labels} and~\ref{axiom.labels incoming}, this forces $\beta' = \rev(\beta)$, which shows that the image 
    of $(v,\alpha) \stackrel{K}{\longrightarrow} (w,\beta)$ under $\LL_n$ is precisely the edge in \eqref{eq:step2}. We thus add 
    $(v,\alpha) \stackrel{K}{\longrightarrow} (w,\beta)$ and 
    $(w',\rev(\beta)) \stackrel{K^\LL}{\longrightarrow} (v',\rev(\alpha))$ to $G'$. \medskip

    \item Suppose $(w,\beta)\in G'$. Then $\LL_n$ sends $(w, \beta)$ to $(w', \rev(\beta)) \in G'$. We show there is an edge 
    \begin{equation}\label{eq:step2dual} (w',\rev(\beta)) \stackrel{K^\LL}{\longrightarrow} (v',\alpha') \end{equation}
    that is length-preserving. Again, axioms~\ref{axiom.labels} and~\ref{axiom.labels incoming} imply that $\alpha' = \rev(\alpha)$, and so the image of 
    $(v,\alpha) \stackrel{K}{\longrightarrow} (w,\beta)$ under $\LL_n$ is precisely the segment in \eqref{eq:step2dual}. We add these segments to $G'$ as well.
    \end{enumerate} \medskip

We claim that all of the vertices in $G_{=\ell +1}$ can be reached by increasing/decreasing edges to $G_{\leqslant \ell}$ and length-preserving edges 
within $G_{=\ell+1}$, and thus we will eventually obtain $G' = G_{\leqslant \ell+1}$ through this process.  Since $\LL_n(G') = G'$ at each step, this will imply 
that $\LL_n(G_{\leqslant \ell+1}) = G_{\leqslant \ell+1}$.

To see why the claim is true, note that by assumption $\ell +1 > s$, and so if $\len(\alpha) = \ell+1$, then $\alpha \neq \lambda$ and 
$\alpha \neq \rev(\lambda)$, which lie in $G_s$ by Lemma~\ref{lemma.top crystal}. It follows from axiom~\ref{axiom.local edge types} that every 
vertex in $G_{=\ell+1}$ contains an incoming edge that is length-preserving or increasing. Length-preserving edges decrease labels in dominance
order, hence $G_{=\ell+1}$ is acyclic. It follows that because $G_{=\ell+1}$ is finite, within each connected component of $G_{=\ell+1}$ there is a 
vertex $(w,\beta)$ with no incoming length-preserving edges. Thus by axiom~\ref{axiom.local edge types}, $(w,\beta)$ must have an incoming 
length-increasing edge from a vertex $(x,\gamma)$, which forces $(x,\gamma) \in G_{\leqslant \ell}$. This implies the claim.  

 \medskip
We add the edges to $G'$ in a specific order to be able to use the commutation relations in Theorem~\ref{theorem.commutation}. 
\begin{enumerate}
\item[\emph{Step 2a}.] We prove (1) under the assumption that there exists some edge $(v,\alpha) \stackrel{L}{\longrightarrow} (x,\gamma)$ in $G$ 
which is length-decreasing. Under duality $\LL_n$, this proves (2) when there is a length-increasing edge into $(w,\beta)$. \smallskip
\item[\emph{Step 2b}.] We prove (1) under the assumption that there exists some edge $(x,\gamma) \stackrel{L}{\longrightarrow} (v,\alpha)$ in $G$ 
which is length-increasing. Under duality $\LL_n$, this proves (2) when there is a length-decreasing edge out of $(w,\beta)$.
\end{enumerate}
Note that once Steps 2a and 2b are completed, $G'$ contains all length-preserving edges into and out of any vertex $(v,\alpha) \in \Gc \cap G_{=\ell+1}$. 
We proceed as follows: \smallskip
\begin{enumerate}
\item[\emph{Step 2c}.] We prove (1) when $(v,\alpha)$ has no edge connecting it to $G_{\leqslant \ell}$. 
By induction, we may assume that there is a length-preserving edge 
\[ (v,\alpha) \xrightarrow{L} (x,\gamma) \quad \quad \textrm{ or} \quad \quad (v,\alpha) \xleftarrow{L} (x,\gamma) \] 
such that $(x,\gamma)\in G'$, along with all its incoming and outgoing length-preserving edges. (In other words, $(v,\alpha)$ was added to $G'$ in 
Steps 2a or 2b or a previous application of Step 2c.) Under duality $\LL_n$, this also proves (2).
\end{enumerate}

\noindent \emph{Proof of Step 2a}.
We first prove (1) under the assumption that there is an edge $(v,\alpha) \stackrel{L}{\longrightarrow} (x,\gamma)$ which is length-decreasing. Note
that this edge and vertex $(x,\gamma)$ are in $G'$ by step 1.
Since $L$ is length-decreasing and $K$ is length-preserving, we must have a commutation relation by the cases in Theorem~\ref{theorem.commutation}.
Note that by Remark~\ref{remark.commutations duality} the commutation relations are dual under $\LL_n$. Hence if a commutation holds at $(v,\alpha)$,
then the dual commutation must hold at $(v',\rev(\alpha))$. We go through each of the commutation relations from Theorem~\ref{theorem.commutation} in turn.
\begin{itemize}
\item Suppose case~\ref{case1} holds. Then since $L$ is length-decreasing, all edges and vertices in the square are in $G'$ by induction and hence so are
their duals under $\LL_n$. Furthermore, the edge labeled $K$ out of $(x,\gamma)$ in $G$ is length-preserving as is the edge $K^\LL$ out of 
$(x',\rev(\gamma))$. This forces the edge $(w',\beta') \stackrel{K^\LL}{\longrightarrow} (v',\rev(\alpha))$ to be length-preserving. \smallskip
\item 
Case~\ref{case2ai} is not applicable since all edges in that case are length-preserving. \smallskip
\item  In case~\ref{case2aii}, we must have $I=K$ and $J=L$
since $J$ is length-decreasing. Since $J'$ is length-decreasing, we have that $(w,\beta)\in G'$ and hence 
$(w',\beta') \stackrel{K^\LL}{\longrightarrow} (v',\rev(\alpha))$ is length-preserving. \smallskip
\item  If case~\ref{case2bi} holds, we must have $J=K$ and $I=L$ since only $I$ can be length-decreasing. By the dual octagon, $L^\LL$ is length-increasing 
and $K^\LL$ is length-preserving as edges into $(v',\rev(\alpha))$. \smallskip
\item In case~\ref{case2bii}, we must have $J=L$ and $I=K$ since $J$ is length-decreasing. In particular, $I$ cannot be length-increasing in this case.
Hence in the dual $I^\LL=K^\LL$ cannot be length-decreasing.
If the edge $(w',\beta') \stackrel{K^\LL}{\longrightarrow} (v',\rev(\alpha))$ is length-increasing,
then $(w',\beta') \in G_{\leqslant \ell}$ and hence by induction its dual $(v,\alpha) \stackrel{K}{\longrightarrow} (w,\beta)$ is length-decreasing, contradicting
the assumption that it is length-preserving. Hence $K^\LL$ is length-preserving, proving (1). \smallskip
\item  The arguments for case~\ref{case2biii} are analogous to those in case~\ref{case1}. \smallskip
\item Case~\ref{case3} is not applicable since $J$ is length-increasing, but neither $K$ nor $L$
are length-increasing. 
\end{itemize}
 This concludes the proof of (1) when there is a length-decreasing edge out of $(v,\alpha)$. 
 
 \medskip
 
 \noindent \emph{Proof of Step 2b}.
Next we prove (1) under the assumption that there is an edge $(x,\gamma) \stackrel{L}{\longrightarrow} (v,\alpha)$ which is length-increasing. 
Let $L=[j,j+2m]$. We perform induction on $m$. If $|L|=3$ or equivalently $m=1$, then by~\eqref{equation.loop}, there is also an outgoing 
length-decreasing edge from $(v,\alpha)$. Hence by step 2a, (1) holds.
Therefore, we may assume that $|L|>3$ or equivalently $m>1$. By the fans (which are implied by the commutations in Theorem~\ref{theorem.commutation}
and~\eqref{equation.loop}), there is a length-preserving edge
\begin{equation}
\label{equation.y delta}
	(v,\alpha) \xrightarrow{[j+2m-1,j+2m+1]} (y,\delta),
\end{equation}
which is in $G'$ since both $(v,\alpha)$ and $(y,\delta)$ have incoming length-increasing edges.
By case~\ref{case3b} starting at vertex $(x,\gamma)$ with $I=L$ and $J=[j+1,j+2m-1]$, the edge $J'=[j,j+2m-2]$ out of $(v,\alpha)$ is
length-increasing. Hence $K$ cannot be $[j,j+2m-2]$ since it is length-preserving. We will now use the commutation relations of 
Theorem~\ref{theorem.commutation} on the edge $(v,\alpha) \xrightarrow{K} (w,\beta)$ and the edge in~\eqref{equation.y delta} to prove (1):
\begin{itemize}
\item
If $K \cap [j+2m-1,j+2m+1] = \emptyset$, case~\ref{case1} applies. By induction on $m$, all length-preserving edges leaving $(y,\delta)$
are in $G'$. Hence under $\LL_n$ they are still length-preserving. Since $K\neq [j,j+2m-2]$, we cannot be in the case that $\max K +1 =j+2m-1$.
Hence by case~\ref{case1}, the edge $K^\LL$ into $(v',\rev(\alpha))$ is length-preserving, proving (1). \smallskip
\item The case $[j+2m-1,j+2m+1]\subseteq K$ is not possible since it would imply that the edge in~\eqref{equation.y delta} is length-increasing,
a contradiction to the fact that it is length-preserving. \smallskip
\item If $K \cap [j+2m-1,j+2m+1] \neq \emptyset$, we are in case~\ref{case2}. If $I=K$ and $J=[j+2m-1,j+2m+1]$, then $I$ cannot be
length-increasing. Hence in the dual commutation, $I^\LL$ cannot be length-decreasing. We also cannot have that $I^\LL=K^\LL$ is length-increasing 
since then the edge 
\begin{equation}
\label{equation.dual edge}
	 (w',\beta') \stackrel{K^\LL}{\longrightarrow} (v',\rev(\alpha))
\end{equation}
was in $G'$ implying that 
\[ 
	(v,\alpha) \stackrel{K}{\longrightarrow} (w,\beta) 
\] 
is length-decreasing, contradicting the fact that this edge is length-preserving. Hence the edge in~\eqref{equation.dual edge} is
length-preserving, proving (1). 

If $I=[j+2m-1,j+2m+1]$ and $J=K$, then case~\ref{case2ai} is not applicable since $|I|=3$. Cases~\ref{case2aii}
and~\ref{case2bii} are not applicable since $J=K$ is length-preserving and not length-decreasing. In case~\ref{case2bi}, $J=K$ is required to be 
length-preserving and hence in its dual under $\LL_n$ it also has to be length-preserving, proving (1). The proof in case~\ref{case2biii} is similar
to case~\ref{case1}.
\end{itemize}

We have now established that all length-preserving incoming and outgoing edges from $(v,\alpha)\in G_{\leqslant \ell+1}$ with an edge connecting to 
$G_{\leqslant \ell}$ are in $G'$. 

\medskip

\noindent \emph{Proof of Step 2c}. 
Finally, we prove (1) when $(v,\alpha) \in G'$ has no edge connecting to $G_{\leqslant \ell}$. By induction, we may assume
that there is an edge 
\[ (v,\alpha) \xrightarrow{L} (x,\gamma) \quad \quad \textrm{ or} \quad \quad (v,\alpha) \xleftarrow{L} (x,\gamma) \] 
such that $(x,\gamma)\in G'$ with all 
incoming and outgoing length-preserving edges in $G'$.

First assume that $(v,\alpha) \xrightarrow{L} (x,\gamma)$ as above exists.
We are now going to consider the commutation relations as in Theorem~\ref{theorem.commutation} starting at vertex $(v,\alpha)$ with length-preserving
outgoing edges $K$ and $L$.
\begin{itemize}
    \item Cases~\ref{case2aii} and~\ref{case2bii} are not applicable since $J$ is length-decreasing. \smallskip
    \item Cases~\ref{case3a} and~\ref{case3b} are not applicable since $J$ is length-increasing. \smallskip
    \item In case~\ref{case2ai}, all edges are length-preserving and hence also in its dual, so (1) must hold. \smallskip
    \item In case~\ref{case2bi}, edge $I$ cannot be length-decreasing since otherwise there is an edge connecting $(v,\alpha)$ to $G_{\leqslant \ell}$, 
    a contradiction to our assumptions. For the same reason $I^\LL$ in the dual commutation cannot be length-increasing. Hence both $I$ and $J$ 
    are length-preserving and also in the dual, proving (1). \smallskip 
\item 
For case~\ref{case2biii} the edges $I$ and $I'$ (resp. $J$ and $J'$) are of the same type. In this setting $I=K$ and $J=L$ or vice versa. 
By induction hypothesis, the edges on one side of the square are in $G'$ and length-preserving. Hence the edges on the other side of the square
must also be length-preserving since by case~\ref{case2biii} are of the same type. This proves (1).
Case~\ref{case1a} can be proved in the same fashion.
\smallskip
\item 
Case~\ref{case1b} is not applicable since in this case $I$ is length-increasing and should be $K$ or $L$. But both $L$ and $K$ are length-preserving.
\end{itemize}

Next we assume that $(v,\alpha) \xleftarrow{L} (x,\gamma)$ as above exists. Hence we have
\begin{equation}
\label{equation.string in G'}
	(x,\gamma) \xrightarrow{L} (v,\alpha) \xrightarrow{K} (w,\beta),
\end{equation}
where both edges are length-preserving and $(x,\gamma),(v,\alpha)\in G'$.
When $K \cap L = \emptyset$, there is an edge $(x,\gamma) \xrightarrow{K} \bullet$ by axiom~\ref{axiom.outgoing}.
Hence case~\ref{case1} of Theorem~\ref{theorem.commutation} applies. In case~\ref{case1a}, all edges must be length-preserving since the 
edges in~\eqref{equation.string in G'} are length-preserving. Since the length-preserving edges out of $(x,\gamma)$ are in $G'$ by induction, 
they also must be length-preserving in the dual showing that $(w',\beta') \stackrel{K^\LL}{\longrightarrow} (v',\rev(\alpha))$ is length-preserving by
case~\ref{case1a}, proving (1). In case~\ref{case1b}, since the edges in~\eqref{equation.string in G'} are length-preserving 
we must have $K=I$. Since the edge $(x,\gamma) \xrightarrow{L} (v,\alpha)$ is in $G'$, it must be length-preserving in the dual forcing the dual of the edge 
labeled $K$ in~\eqref{equation.string in G'} to be length-preserving, proving (1). The case $L\subseteq K$ is not possible since $\gamma \geqslant_D \alpha$
and hence by axiom~\ref{axiom.outgoing} the edge $K$ in~\eqref{equation.string in G'} could not exist. If $K\subseteq L$, $K$ in~\eqref{equation.string in G'}
would be length-increasing by axiom~\ref{axiom.incoming}, contradicting that $K$ is length-preserving. 

Hence we may now assume that $K\cap L \neq \emptyset$, $K \not \subseteq L$ and $L \not \subseteq K$.
First assume that $\min L < \min K$. By axiom~\ref{axiom.outgoing}, there is an outgoing edge from $(x,\gamma)$ labeled $K=[k,k+2b]$ or
$[k+1,k+2b-1]$. By axiom~\ref{axiom.outgoing} this would not happen if there was an incoming edge $(y,\delta) \xrightarrow{M} (x,\gamma)$ with 
$M\subseteq K$ and $\delta \geqslant_D \gamma$. But then $(v,\alpha)$ would have such an incoming edge contradicting the fact the edge
$(v,\alpha) \xrightarrow{K} (w,\beta)$ exists.
We are now in case~\ref{case2} of Theorem~\ref{theorem.commutation} with $I=L$, $J=K$ or $J=[k+1,k+2b-1]$, and $J'=K$.
\begin{itemize}
\item
If case~\ref{case2ai} holds, all edges are length-preserving and hence also in the dual, proving (1).
\item
Case~\ref{case2aii} is not applicable since $J'$ is length-decreasing, contradicting our assumption that all edges in~\eqref{equation.string in G'}
are length-preserving.
\item
Suppose case~\ref{case2bi} applies. Then $J'=K=[k,k+2b]$ is length-preserving or length-increasing. Hence its dual can be length-decreasing 
or length-preserving. In the latter case (1) holds and we are done. If it is length-decreasing, then the dual of case~\ref{case2bi} stipulates
that the dual of case~\ref{case2ai} also holds. This would imply an edge $[k-1,k+2b+1]$ leaving $(v,\alpha)$. But then by the fans, the edge
$(v,\alpha) \xrightarrow{K} (w,\beta)$ is length-increasing, contradicting our assumptions that it is length-preserving.
\item
If case~\ref{case2bii} holds, then edge $J''$ is length-decreasing implying that $(w,\beta)$ is attached to $G_{\leqslant \ell}$ by an edge. In this case,
$(w,\beta)\in G'$ by a previous case, proving (1).
\item
In case~\ref{case2biii}, the edges labeled $J=J'=K$ have the same type and are length-preserving in this case. Since the edge labeled $K$
out of $(x,\gamma)$ is in $G'$, it must be length-preserving under $\LL_n$, proving (1) using case~\ref{case2biii} for the dual.
\end{itemize}

It remains to consider the case $\min L > \min K$. By the previous cases, all length-preserving edges $M$ out of $(v,\alpha)$ with 
$M \cap L = \emptyset$ or $\min L < \min M$ are in $G'$. If such an edge $M$ exists, then it can be used as the interval $L$ in the first case in the 
proof of Step 2c to conclude that $(w,\beta)\in G'$. Hence we may assume that no such edge $M$ exists. 
By Theorem~\ref{theorem.commutation} case~\ref{case3b} the edge $K$ is hence the only length-preserving edge out of $(v,\alpha)$ and all
potential other edges $J$ out of $(v,\alpha)$ satisfy $J \subsetneq K$ and are length-increasing. By the dual commutation relations of 
Theorem~\ref{theorem.commutation}, $(w',\beta') \stackrel{K^\LL}{\longrightarrow} (v',\rev(\alpha))$ is hence the only length-preserving edge into
$(v',\rev(\alpha))$ proving (1).
\end{proof}

Next, we show that if a graph $G$ satisfies Axiom~\ref{local axioms}, then so does each component of $G_{[1,n-1]}$.

\begin{lemma}\label{local axiom on restriction}
    Suppose $G$ satisfies Axiom~\ref{local axioms}. Then each connected component of $G_{[1,n-1]}$ does as well.
\end{lemma}

\begin{proof}
It is not hard to see that axioms~\ref{axiom.local intervals}-\ref{axiom.local labels} and axiom~\ref{axiom.local string} still hold
for $G_{[1,n-1]}$. 

Regarding axiom~\ref{axiom.local commutations}, the commutation relations
in Theorem~\ref{theorem.commutation} still hold as long as all intervals appearing in the relations do not contain $n$. Assuming that $I$ and
$J$ do not contain $n$, then no interval in cases~\ref{case1}, \ref{case2bi}, \ref{case2bii}, \ref{case2biii}, \ref{case3b} contains $n$.

In case~\ref{case2ai}, the interval $J'=[j-1,j+2\ell+1]$ contains $n$ if $j+2\ell+1=n$. By axiom~\ref{axiom.local outgoing}, the vertex from which
$J'$ originates also has an outgoing edge labeled $J=[j,j+2\ell] \subseteq J'$. Hence in the absence of $J'$ (when $n$ is removed) 
case~\ref{case2bi} holds which implies an octagon relation. Similarly, in case~\ref{case2aii} we have $n\in J'=[i+1,i+3+2\ell]$ if $n=i+3+2\ell$. 
By axiom~\ref{axiom.local outgoing}, the vertex from which $J'$ originates also has an outgoing edge labeled $J=[i+2,i+2+2\ell]$. Hence 
in the absence of $J'$ (when $n$ is removed) case~\ref{case2bii} holds which implies a pentagon relation. 
Furthermore, in case~\ref{case3a} we have $n\in J'=[i+2m-1,i+2m+1]$ if $n=i+2m+1$.
Removing $n$ changes $I$ from length-increasing to length-preserving. But then case~\ref{case3b} still applies.
Finally for~\eqref{equation.loop} if $n\in J$, then $I$ in $G_{[1,n-1]}$ is no longer length-increasing, so the loop does not need to exist
any longer. 

We next show $G_{[1,n-1]}$ satisfies axiom~\ref{axiom.local edge types}.
The partition $\lambda$ in axiom~\ref{axiom.local edge types} is the same as the weight of the crystal $B(\lambda)_s$ in axiom~\ref{axiom.crystal},
which we showed must hold by Lemma~\ref{lemma.top crystal}. Suppose $\rev(\alpha)\neq \lambda$. By axiom~\ref{axiom.local edge types}, there 
exists a length-preserving (resp. length-decreasing) edge $(v,\alpha) \xrightarrow{L} (w,\beta)$ in $G$. 
\begin{itemize}
    \item If $n\not \in L$, this edge is still length-preserving (resp. length-decreasing) in $G_{[1,n-1]}$.  \medskip 
    \item If $L=[i,n]$ with $|L|>3$,  there is a length-increasing edge $(v,\alpha) \xrightarrow{J} (x,\gamma)$ with $J=[i+1,n-1]$ in $G$. 
    This turns into a length-preserving edge in $G_{[1,n-1]}$ since the above forces $\gamma = (\gamma_1, \gamma_2, \ldots, \gamma_{\ell-1},1)$. \medskip
    \item If $L=[n-2,n]$ is length-decreasing, by~\eqref{equation.loop} 
	in axiom~\ref{axiom.local commutations} there exists a length-increasing edge 
	\[ (w,\beta) \xrightarrow{[n-3,n-1]} (v,\alpha) \] in $G$. Furthermore,
	$\alpha=(\alpha_1,\ldots,\alpha_{\ell-2},2,1)$ and $\beta=(\alpha_1,\ldots,\alpha_{\ell-2}+1,2)$. If $\rev(\beta)\neq\lambda$, by 
	axiom~\ref{axiom.local edge types} in $G$ there must be an edge $(w,\beta) \xrightarrow{K} (x,\gamma)$ which is not length-increasing.
	Employing the commutation relations of Theorem~\ref{theorem.commutation} with $I=K$ and $J=[n-3,n-1]$, only case~\ref{case1}
	or~\ref{case2biii} can apply since $J$ is length-increasing. In these cases, $I'=I=K$ out of $(v,\alpha)$ is not length-increasing and does not
	contain $n$. Hence it is still length-preserving or length-decreasing in $G_{[1,n-1]}$.  
	If $\rev(\beta)=\lambda$, we have $\beta=(1^k, 2^{\ell-1-k})$ and $\alpha=(1^k,2^{\ell-3-k},1,2,1)$.
	If $\ell-3-k>0$, $(v,\alpha)$ has an outgoing length-preserving edge $[n-5,n-3]$, which remains length-preserving in $G_{[1,n-1]}$.
	If $\ell-3-k=0$, the label $(\alpha_1,\ldots,\alpha_{\ell-1})=(1^{k+1},2)$ in $G_{[1,n-1]}$ is minimal in its component.
	\medskip
    \item The last case to consider is when $L=[n-2,n]$ is length-preserving in $G$. Note that in this case $\alpha=(\alpha_1,\ldots,\alpha_{\ell-1},1)$. 
    If this edge is the only length-preserving/decreasing edge out of $(v,\alpha)$, the label $(\alpha_1,\ldots,\alpha_{\ell-1})$ is minimal in its 
    connected component.
\end{itemize}
This proves axiom~\ref{axiom.local edge types}(b)
in $G_{[1,n-1]}$. Axiom~\ref{axiom.local edge types}(a) in $G_{[1,n-1]}$ follows similarly.
\end{proof}

Recall that a $\CS$-graph is a graph satisfying Axiom~\ref{axioms}. We are finally ready to prove that the local axioms characterize $\CS$-graphs as well. 

\begin{theorem}
A graph $G$ satisfying Axiom~\ref{local axioms} is a $\CS$-graph.
Conversely, a $\CS$-graph satisfies Axiom~\ref{local axioms}.
\end{theorem}

\begin{proof}
Suppose $G$ satisfies Axiom~\ref{local axioms}. Axioms~\ref{axiom.local intervals}-\ref{axiom.local labels} imply
axioms~\ref{axiom.intervals}-\ref{axiom.labels}. Axiom~\ref{axiom.fans} follows from the iteration of case~\ref{case3a} in 
Theorem~\ref{theorem.commutation} together with~\eqref{equation.loop}. Axiom~\ref{axiom.crystal} follows from 
Lemma~\ref{lemma.top crystal}. Axiom~\ref{axiom.lusztig}~\ref{axiom.lusztig n} follows from Lemma~\ref{lemma.lusztig local}.
It remains to prove axiom~\ref{axiom.lusztig}~\ref{axiom.lusztig n-1}.

We proceed by induction on $n$. For $n=1,2$, both Axiom~\ref{axioms} and Axiom~\ref{local axioms} imply that $G$ is a single vertex 
since no odd length intervals $I \subseteq [n]$ with $|I|\geqslant 3$ exist. Assume by induction that Axiom~\ref{local axioms}
implies Axiom~\ref{axioms} for graphs with index $n-1$ or smaller.

By Lemma~\ref{local axiom on restriction}, if $G$ satisfies Axiom~\ref{local axioms}, then the connected components of 
$G_{[1,n-1]}$ satisfy Axiom~\ref{local axioms}.
Thus, the connected components of $G_{[1,n-1]}$ satisfy Axiom~\ref{axioms} by induction. In particular,
axiom~\ref{axiom.lusztig}~\ref{axiom.lusztig n} for $G_{[1,n-1]}$ states that 
\[ 
	\LL_{n-1}(G_{[1,n-1]}) \cong G_{[1,n-1]},
\]
which is precisely axiom~\ref{axiom.lusztig}~\ref{axiom.lusztig n-1} for $G$.

We conclude by proving that a $\CS$-graph satisfies Axiom~\ref{local axioms}.
Axioms~\ref{axiom.intervals}-\ref{axiom.labels} together with Lusztig involution~\ref{axiom.lusztig} imply 
axioms~\ref{axiom.local intervals}-\ref{axiom.local labels}.
By Corollary~\ref{corollary.uniqueness}, a $\CS$-graph is isomorphic to $\CS(\lambda)$ for some partition $\lambda$.
Axiom~\ref{axiom.local commutations} hence follows from Theorem~\ref{theorem.commutation},
Corollaries~\ref{corollary.fans shorter} and~\ref{corollary.fans longer}, and their duals by Lusztig involution. 
Axiom~\ref{axiom.local string} follows from Theorem~\ref{theorem.B short} and Corollary~\ref{corollary.string}.
Axiom~\ref{axiom.local edge types} follows from Proposition~\ref{proposition.incoming and outgoing}.
\end{proof}

\section{Crystal skeleton for two-row partitions}
\label{section.two row}
\def\rcomp{\mathsf{rcomp}}
\def\cR{\mathcal{R}}
\def\cC{\mathcal{C}}

In this section, we illustrate our general findings for crystal skeletons in the two-row case, where the combinatorics is particularly elegant.

\subsection{Vertices and edges in the two-row crystal skeleton}
When $\lambda$ has two rows, the vertices and edges of $\CS(\lambda)$ can be translated into the language of paths. Recall that vertices of $\CS(\lambda)$ are indexed by $\SYT(\lambda)$.  

Fix $\lambda = (\lambda_1, \lambda_2)$. Then standard Young tableaux of shape $\lambda$ are in bijection with sequences 
\begin{equation}
\label{lattice path height}
	p = (p_1,\ldots , p_{\lambda_1+\lambda_2}) \in \{1, -1\}^{\lambda_1+\lambda_2} \quad \text{ satisfying } \quad 
	h_k(p):=\sum_{i=1}^k p_i \leqslant 0,
\end{equation}
for all $1 \leqslant k \leqslant \lambda_1+\lambda_2$. 
A tableau $T \in \SYT(\lambda)$ (and thus a vertex in $\CS(\lambda)$) corresponds to the sequence 
\[p_i = \begin{cases}
    -1 & \text{ if $i$ appears in the first row of $T$,}\\
    1 & \text{ if $i$ appears in the second row of $T$.}
    \end{cases}\]
We associate $p$ to the path from $(0,0)$ to $(\lambda_1+\lambda_2, -\lambda_1+\lambda_2)$, where at step $i$ the path proceeds in the direction 
$(1, p_i)$. See Example~\ref{ex:Two-row tableaux, paths, rectangular decompositions} below. We call $h_k=h_k(p)$ in~\eqref{lattice path height} 
the \defn{height} of path $p$ at the end of step $k$ and set $h_0=0$. 

It will be useful at times to think of a path as an arrangement of boxes, piled into the corner formed by the row-reading path ($\lambda_1$ 
downward steps, followed by $\lambda_2$ upward steps), where boxes never get piled above height $0$; this is illustrated in 
Example~\ref{ex:Two-row tableaux, paths, rectangular decompositions}. To recover a path from an arrangement of boxes, trace the upper face 
of the arrangement. 

Given a tableau $T$, denote the corresponding path by $p^T$; conversely denote the tableau associated to a path $p$ by $T_p$. The word 
$\row(p) := \row(T_p)$ is obtained from $p$ by reading (left to right) the locations of the upward steps ($i$ such that $p_i = 1$), followed by the 
locations of the downward steps. As a result, descents in $\row(p)$ occur precisely at (internal) local minima (where a down-step is followed by an 
up-step). 

\begin{example}
\label{ex:Two-row tableaux, paths, rectangular decompositions} 
Consider $\lambda = (15,10)$. The tableau
\[T = \TAB[.4]{{1,2,4,5,9,10,11,12,15,17,18,20,23,24,25},{3,6,7,8,13,14,16,19,21,22}}\]
corresponds to the sequence 
\[p = p^T = (-1,-1,1,-1,-1,1,1,1,-1,-1,-1,-1,1,1,-1,1,-1,-1,1,-1,1,1,-1,-1,-1)\]
drawn as
\[\begin{array}{c@{\qquad}c}
\text{path:} & \text{box arrangement:}\\[3pt]
\TIKZ[scale=.25]{
	\Path{1,1,2,1,1,2,2,2,1,1,1,1,2,2,1,2,1,1,2,1,2,2,1,1,1}
	\foreach \x in {1,2,4,5,9,10,11,12,15,17,18,20,23,24,25}{\node[below left, black!40, inner sep = .5pt] at (c\x) {$\scriptscriptstyle \x$};}
\foreach \x in {3,6,7,8,13,14,16,19,21,22}{\node[above left, black!40, inner sep = .5pt] at (c\x) {$\scriptscriptstyle \x$};}
	\node[v, label={above, inner sep = 1pt}:{$(0,0)$}] at (0,0) {};
	\node[v, label={right, inner sep = 1pt}:{$\scriptstyle \begin{matrix}(\lambda_1+\lambda_2,\quad \\\quad -\lambda_1+\lambda_2)\end{matrix}$}] at (25,-5) {};
}
&
\TIKZ[scale=.25*1.41, rotate=45]{
	\Part{10,10,10,8,7,7,6,4,4,4,4,1,1}
	\draw[very thick] (0,15) to (0,0) to (10,0);
	\draw[densely dotted] (0,15) to (10,5) to (10,0);
}
\end{array}\]
Reading the positions of the up-steps in $p$, followed by the positions of the down-steps, recovers 
\[\row(p) = 3\ 6\ 7\ 8\ 13\ 14\ 16\ 19\ 21\ 22\ 1\ 2\ 4\ 5\ 9\ 10\ 11\ 12\ 15\ 17\ 18\ 20\ 23\ 24\ 25 = \row(T).\]
The (right) descents of $\row(p)$ occur exactly at the 6 internal local minima of $p$. 
\end{example}

As before, if $I = [i,i+2m]$ is a Dyck pattern interval in $\pi = \row(p)$, we write 
\[
	I \cdot p :=  \cycle(\pi|_I) p,
\]
where $\cycle(\pi|_I)$ acts on $p$ according to its action on $T_p$ as in Section~\ref{section:cycles}.

Recall that descents in $\pi=\row(p)$ correspond to local minima in $p$. It follows that if $I = [i,i+2m]$ is a Dyck pattern interval such that $\row(p|_I)$ 
destandardizes to a word in $i$ and $i+1$, then the path $p|_I$ has exactly one (internal) local minima; the possible local paths are shown 
below in Lemma \ref{lem:two-row Dyck patterns}. We can thus succinctly describe the paths $p|_I$ that appear in $\CS(\lambda)$, as well as 
the corresponding adjacent vertex.

\begin{lemma}\label{lem:two-row Dyck patterns} 
Suppose $\lambda =(\lambda_1, \lambda_2) $ has two rows. There is an edge $p \xrightarrow{~I~} p'$ in $\CS(\lambda)$ if and only if $p$ and 
$p'$ agree on any $j \not \in I$, and on $I$ we have the following local behavior, where $w = \destd(\row(p)|_I)$ is the destandardization of the subword of $\row(p)$ over the interval $I$:

\begin{equation}\tag{A}\label{eq:two-row-A-move}
\raisebox{-0.5\height}{\includegraphics{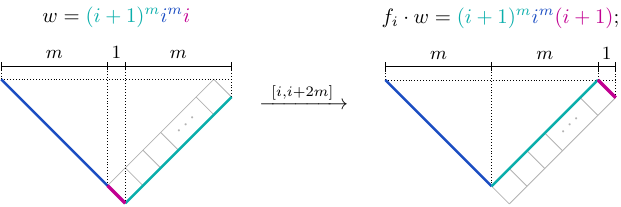}}
\end{equation}
or
\begin{equation}\tag{B}\label{eq:two-row-B-move}
\raisebox{-0.5\height}{\includegraphics{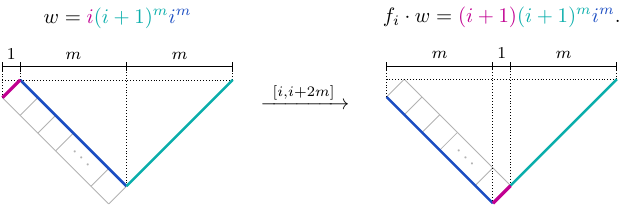}}
\end{equation}
\end{lemma}
We will refer to the case in \eqref{eq:two-row-A-move} as a type \eqref{eq:two-row-A-move} edge, and the case in \eqref{eq:two-row-B-move} 
as a type \eqref{eq:two-row-B-move} edge. As depicted above, in the box-arrangement model, we can think of an edge of 
type \eqref{eq:two-row-A-move} as adding a right strip of boxes, 
and an edge of type \eqref{eq:two-row-B-move} as removing a left strip of boxes. 

\subsection{Rectangular decompositions and strongly-connected components}

Throughout this paper, we have seen the importance of intervals $I= [a,b]$ where jeu de taquin takes $T|_I$ to a rectangular 
shape. For two-row partitions, the rectangles can only have one or two rows. In particular, we have the following cases:
\begin{itemize}
\item \emph{One row:} $\jdt(T|_I)$ forms a single row precisely when $\row(T)|_I$ has no descents. In terms of the path, this means $p^T$ has 
no local minima over $I$, that is, $p^T|_I = (1)^\ell(-1)^m$. 
\item \emph{Two rows:}  $\jdt(T|_I)$ forms a two-row rectangle $(m, m)$ precisely when $p^T|_I$ is a \defn{Dyck subpath}, meaning that $p^T|_I$ 
starts and ends at the same height and never exceeds that height:
\[h_{a-1} = h_b \quad \text{ and } \quad h_i \leqslant h_{a-1} \text{ for all } i \in [a,b]\]
(so that $b = a+2m-1$). 
\end{itemize}
We can then decompose any lattice path $p \in \CS(\lambda_1, \lambda_2)$ uniquely into a sequence of maximal rectangles as
\begin{align}\label{eq:DPdecomp}
p &= (-1)^{n_0} D_1 (-1)^{n_1} D_2 (-1)^{n_2} \cdots  D_\ell (-1)^{n_\ell}, \\
&\quad \text{ where }
\begin{array}{l}
n_1, \dots, n_{\ell-1} \geqslant 1,\\
n_0 + n_1 + \cdots + n_\ell = \lambda_1 - \lambda_2, \quad \text{and}\\
D_k \text{ is a Dyck subpath from height } -(n_0 + \cdots + n_{k-1}).
\end{array}\nonumber
\end{align}
\noindent We call \eqref{eq:DPdecomp} the (maximal) \defn{rectangular decomposition} of $p$ (or of $T_p$). One can also check that 
\[ |D_1| + \cdots + |D_\ell| = 2\lambda_2. \] 
Note that the $(-1)^{n_k}$ in \eqref{eq:DPdecomp} correspond to downsteps in $p$.
Define the \defn{rectangular composition}
\defn{$\rcomp(p)$} to be the set composition of $\{1, \dots, n\}$ corresponding to the intervals spanning the rectangles in \eqref{eq:DPdecomp}, so that 
\begin{equation}\label{defn:rcomp}
\rcomp(p)=\left(\alpha^{(0)}, \beta^{(1)}, \alpha^{(1)}, \dots, \beta^{(\ell)}, \alpha^{(\ell)}\right),
\end{equation}
where with $D_0=\emptyset$
\begin{align*}
\alpha^{(k)} 
	&= \left[1 + \sum_{r = 0}^{k-1} (n_r +  |D_r|), n_k + \sum_{r=0}^{k-1} (n_r +  |D_r|)\right], 
	&\text{i.e. the interval covering the terms $(-1)^{n_k}$, and}\\
\beta^{(k)} 
	&= \left[1 + n_k + \sum_{r=0}^{k-1} (n_r +  |D_r|), |D_k| + n_k + \sum_{r=0}^{k-1} (n_r +  |D_r|)\right], 
	&\text{i.e. the interval covering the terms in $D_k$.}
\end{align*}
The only terms that can be the empty set are $\alpha^{(0)}$ and $\alpha^{(\ell)}$, which occur when $n_0 = 0$ or $n_\ell = 0$, respectively. Note that $\rcomp(p)$ is not the same composition as the descent composition $\Des(T_p)$; see Section \ref{section:two row descent compositions and fans} for a discussion on descent compositions obtained from a path $p$.  

\begin{example}\label{ex:Two-row rectangular decompositions}
Continuing Example \ref{ex:Two-row tableaux, paths, rectangular decompositions}, the rectangular decomposition of $p$ from \eqref{eq:DPdecomp} is 
\[p = {(-1)^0}\viridian{D_1}{(-1)^2}\plum{D_2}{(-1)^3} \quad \text{ with } \quad \rcomp(p) = (
		\underset{\strut\alpha^{(0)}}{\strut\emptyset}, 
		\viridian{\underset{\strut\beta^{(1)}}{\strut[1,8]}}, 
		\underset{\strut\alpha^{(1)}}{\strut[9,10]}, 
		\plum{\underset{\strut\beta^{(2)}}{\strut[11, 22]}}, 
		\underset{\strut\alpha^{(2)}}{\strut[23,25]}),\]
as illustrated below:
\begin{equation}\label{eq:DPdecomp-ex}
 \TIKZ[scale=.25]{
\foreach \x in {9, 15, 21}{\node[black!50] at (\x,-5.5) {$\scriptstyle \vdots$};}
\begin{scope}
\clip (0,0) rectangle (25,-5.5);
\fill[maize!50] (25,-5) to ++(-3, 3) to ++(-12,0) to ++(-2,2) to (25,0) to (25,-5);
	\begin{scope}[black!20, thin]
	\draw (0,0) to (15, -15);
	\foreach \d in {15, ..., 10}{\draw (\d, -\d) to +(10,10);}
	\foreach \d [evaluate = \d as \h using \d-15, count = \c from 1] in {1, ..., 10}
		{\draw (\d, -\d) to (2*\d,0) to +(-\h,\h);}
	\end{scope}
\end{scope}
\node[v, viridian, very thick] (0) at (0,0){};
\node[v, very thick, black] (e) at (25,-5) {};
\draw [viridian, very thick] 
	(0) 
	to ++(2,-2) 
	to ++(1,1) 
	to ++(2,-2) 
	to ++(3,3) node[v, viridian](d1){}; 
\draw [plum, very thick] 
	(10,-2) node[v, plum](d2a){}
	to ++(2,-2) 
	to ++(2,2) 
	to ++(1,-1) 
	to ++(1,1) 
	to ++(2,-2) 
	to ++(1,1) 
	to ++(1,-1) 
	to ++(2,2) node[v, plum](d2b){}; 
\draw[|-|, viridian] (0,.5) to node[above, inner sep=1pt] (D1) {$\beta^{(1)}$}(8,.5);
\draw[|-|, plum] (10,.5) to node[above, inner sep=1pt] (D2) {$\beta^{(2)}$}(22,.5);
\draw [black, very thick] 
	(d1) 
	to 
		(d2a)
	(d2b) 
	to 
		(e);
\draw[viridian, thick, dotted] (0) to 
	(d1);
\draw[plum, thick, dotted] (d2a) to 
	(d2b) ;
 }
 \end{equation}
Compare this decomposition to that of $T$ into the maximal rectangles:
\[T = \TIKZ[scale = .4]{
\draw[thick, viridian, fill=viridian!20] (0,0) rectangle (4,2);
\draw (0,1) to (4,1);
	\foreach \x in {1,2,3}{\draw (\x,0) to (\x,2);}
\foreach \f [count = \x from 0] in {3,6,7,8}{\node at (\x+.5, 1.5){\small$\f$};}
\foreach \f [count = \x from 0] in {1,2,4,5}{\node at (\x+.5, .5){\small$\f$};}
\node at (5,1) {$\sqcup$};
\begin{scope}[shift={(6,0)}]
\draw[thick, black, fill=black!10] (0,0) rectangle (2,1); \draw (1,0) to (1,1);
	\node at (0+.5, .5){\small$9$}; \node at (1+.5, .5){\small$10$};
\end{scope}
\node at (9,1) {$\sqcup$};
\begin{scope}[shift={(10,0)}]
\draw[thick, plum, fill=plum!20] (0,2) to ++(6,0) to ++(0,-1) to ++(2,0) to ++(0,-1) to ++(-6,0) to ++(0,1) to ++(-2,0) to ++(0,1);
\draw (2,1) to ++(4,0);
	\foreach \x in {1,...,5}{\draw (\x,1) to (\x,2) (\x+2,1) to (\x+2,0);}
\foreach \f [count = \x from 0] in {13,14,16,19,21,22}{\node at (\x+.5, 1.5){\small$\f$};}
\foreach \f [count = \x from 2] in {11,12,15,17,18,20}{\node at (\x+.5, .5){\small$\f$};}
\end{scope}%
\begin{scope}[shift={(11,-4)}]
\draw[|->] (3,3.5) to node[right]{$\jdt$} (3,2.5);
\draw[thick, plum, fill=plum!20] (0,0) rectangle (6,2);
\draw (0,1) to ++(6,0);
	\foreach \x in {1,...,5}{\draw (\x,0) to (\x,2);}
\foreach \f [count = \x from 0] in {13,14,16,19,21,22}{\node at (\x+.5, 1.5){\small$\f$};}
\foreach \f [count = \x from 0] in {11,12,15,17,18,20}{\node at (\x+.5, .5){\small$\f$};}
\end{scope}%
\node at (19,1) {$\sqcup$};
\begin{scope}[shift={(20,0)}]
\draw[thick, black, fill=black!10] (0,0) rectangle (3,1); \draw (1,0) to (1,1) (2,0) to (2,1);
	\node at (0+.5, .5){\small$23$}; \node at (1+.5, .5){\small$24$}; \node at (2+.5, .5){\small$25$};
\end{scope}
}\]
In particular, $\beta^{(1)}$ and $\beta^{(2)}$ correspond to the teal and purple rectangles, respectively. The gray rectangles with a single row 
correspond to the black downsteps $(-1)^2$ (i.e. $\alpha^{(1)}$) and $(-1)^3$ (i.e. $\alpha^{(2)}$) in \eqref{eq:DPdecomp-ex} above.
\end{example}

As depicted in \eqref{eq:DPdecomp-ex}, one can think visually of the terms in $p$ that appear in $(-1)^{n_k}$ in \eqref{eq:DPdecomp} 
as \defn{eastern-exposed} down-steps, or equivalently, those exposed to a light shining horizontally from the right.

Returning to the results of Section \ref{sec:strongly-connected} on strongly-connected components, we can now interpret Theorem~\ref{theorem.rectangle} and 
Corollary~\ref{corollary.rectangle-covers} in the two-row case to say that if $p, q \in \CS(\lambda_1, \lambda_2)$ have the same 
rectangular composition, then they are in the same strongly-connected component. 

In particular, given a rectangular decomposition 
of $p$ as in \eqref{eq:DPdecomp}, one can systematically apply a sequence of type \eqref{eq:two-row-B-move} edges that will take $p$ 
to the path that replaces any $D_i$ with the \defn{trivial Dyck path} $(-1)^{m_i}(1)^{m_i}$ over the same interval, where $m_i = |D_i|$. One 
natural choice of sequence corresponds to removing one strip of boxes below the path at a time, moving right-to-left (see Example \ref{ex:C(rho) 
is connected}). A similar natural sequence of type \eqref{eq:two-row-A-move} edges will walk us back to $p$, adding one full strip of boxes at a 
time left-to-right. 

\begin{example}
\label{ex:C(rho) is connected}
Picking up from Example \ref{ex:Two-row rectangular decompositions}, one sequence of \eqref{eq:two-row-B-move} edges from 
$p|_{\beta^{(2)}}$ to the path that replaces $D_2$ by $(-1)^6(1)^6$ is
\[
	\includegraphics{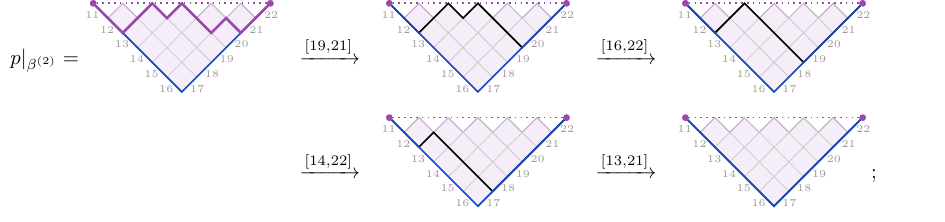}
\]
\noindent and one sequence of \eqref{eq:two-row-A-move} edges moving back to $p$ is
\[\includegraphics{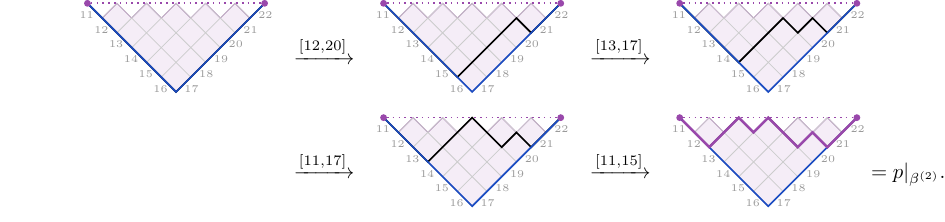}\]
\end{example}

In the two-row case, we can in fact characterize the strongly-connected components of $\CS(\lambda)$: the rectangular compositions \emph{exactly} 
index the strongly-connected components. In order to prove the characterization (Corollary \ref{thm:two-row connected components}), we first prove 
Lemma \ref{lem:two-row Dyck patterns on rcomp} below.

\begin{lemma}
\label{lem:two-row Dyck patterns on rcomp} 
Let $I = [i, i+2m]$ be a Dyck pattern interval in $p \in \CS(\lambda_1, \lambda_2)$. If $I$ is of type \eqref{eq:two-row-B-move}, 
then $\rcomp(I \cdot p) = \rcomp(p)$. If $I$ is of type \eqref{eq:two-row-A-move}, then $\rcomp(I \cdot p) \geqslant_D \rcomp(p)$
in dominance order (under the correspondence between set compositions and compositions).
 \end{lemma}

\begin{proof}
Let $I = [i, i+2m]$ be a Dyck pattern interval in $p$, and consider 
\[\rcomp(p)=\left(\alpha^{(0)}, \beta^{(1)}, \alpha^{(1)}, \dots, \beta^{(\ell)}, \alpha^{(\ell)}\right)\]
 as in \eqref{defn:rcomp}. 

If $I$ is of type \eqref{eq:two-row-B-move}, then since the height of both $p$ and $I \cdot p$ at the left endpoint is less than that of the right 
endpoint, we know $I \subseteq \beta^{(k)}$ for some $k$. Hence $\rcomp(p) = \rcomp(I\cdot p)$. 

On the other hand, if $I$ is of type \eqref{eq:two-row-A-move}, then $[i+1, i+2m] \subseteq \beta^{(k)}$ for some $k$, but we may have 
$i \in \beta^{(k)}$ or $i \in \alpha^{(k-1)}$. If $i \in \beta^{(k)}$, then $\rcomp(p) = \rcomp(I \cdot p)$. Otherwise, there are four cases: 
\begin{itemize}
\item If $\beta^{(k)} = I \setminus \{i\}$, then $D_k$ is precisely the Dyck path in $I$. Under the type~\eqref{eq:two-row-A-move} edge,
this Dyck path moves one step to the left. Otherwise, 
$\beta^{(k)} \supsetneq I \setminus \{i\}$, in which the initial trivial path moves one step to the left of $D_k$. For example, compare 
the effect of $[1,5]$ in Example \ref{ex:Two-row descent compositions} to that of $[10,14]$ in Example \ref{ex:Two-row rectangular decompositions}. 
\item If $n_{k-1} = 1$, then the portion of $D_k$ that moves left will merge with $D_{k-1}$. Otherwise, if $n_{k-1}>1$, then there is room for 
that trivial Dyck path to move left without merging. 
\end{itemize}
In any of these four cases, $\rcomp(p) <_D \rcomp(I \cdot p)$.
\end{proof}

One consequence of Lemma \ref{lem:two-row Dyck patterns on rcomp} is that if an edge $p \xrightarrow{I} p'$ changes the rectangular composition of a 
path, there is no sequence of edges from $p'$ to $p$. 

\begin{corollary}
\label{thm:two-row connected components}
Two paths (i.e. vertices) $p, q \in \CS(\lambda_1, \lambda_2)$ are in the same strongly-connected component of $\CS(\lambda_1, \lambda_2)$ 
if and only if $\rcomp(p) = \rcomp(q)$. 
Hence, the strongly-connected components of $\CS(\lambda_1, \lambda_2)$ are indexed by rectangular compositions.
\end{corollary}

\subsection{Descent compositions and fans}\label{section:two row descent compositions and fans}
\begin{figure}[t]
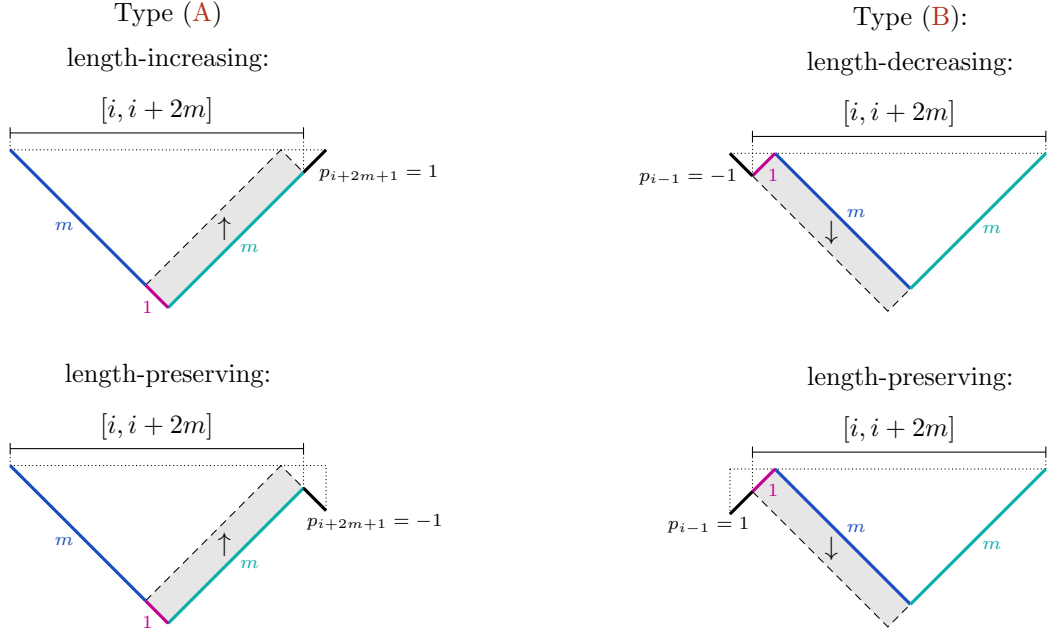

\[\TIKZ[scale=.3]{
\node at (7, 7) {Type \eqref{eq:two-row-A-move}};
\begin{scope}
\node at (7, 5) {length-increasing:};
	\draw[fill=black!10, densely dashed] (6,-5) to ++(6,6) to ++(1,-1) to ++(-6,-6) to ++(-1,1); 
	\node at (9.5,-2.5) {$\uparrow$};
\begin{scope}[densely dotted]
	\foreach \x/\y [count = \c from 1] in {0/1, 13/0}
		{\draw (\x,\y) to (\x,1.5);
		\draw[solid] (\x,1.5) to coordinate (m\c) (\x, 2);}
	\draw (0,1) to ++(14,0);
\end{scope}
\draw (m1) to node[above]{$[i, i+2m]$} (m2);
\begin{scope}[inner sep=1pt,  font=\scriptsize]
\draw[very thick, sapphire] (0,1) to node[below left]{$m$} ++(6,-6);
\draw[very thick, viridian] (7,-6) to node[below right]{$m$}  ++(6,6);
\draw[very thick, rose] (6,-5) to  node[below left]{$1$} ++(1,-1);
\draw[very thick, black] (13,0) to node[below right]{$p_{i+2m+1}=1$}  ++(1,1);
\end{scope}
\end{scope}
\begin{scope}[shift={(0,-14)}]
\node at (7, 5) {length-preserving:};
	\draw[fill=black!10, densely dashed] (6,-5) to ++(6,6) to ++(1,-1) to ++(-6,-6) to ++(-1,1); 
	\node at (9.5,-2.5) {$\uparrow$};
\begin{scope}[densely dotted]
	\foreach \x/\y [count = \c from 1] in {0/1, 13/0}
		{\draw (\x,\y) to (\x,1.5);
		\draw[solid] (\x,1.5) to coordinate (m\c) (\x, 2);}
	\draw (0,1) to ++(14,0) to ++(0,-2);
\end{scope}
\draw (m1) to node[above]{$[i, i+2m]$} (m2);
\begin{scope}[inner sep=1pt, font=\scriptsize]
\draw[very thick, sapphire] (0,1) to node[below left]{$m$} ++(6,-6);
\draw[very thick, viridian] (7,-6) to node[below right]{$m$}  ++(6,6);
\draw[very thick, rose] (6,-5) to node[below left]{$1$} ++(1,-1);
\draw[very thick, black] (13,0) to   ++(1,-1);
\node[below right] at (13,-1) {$p_{i+2m+1}=-1$};
\end{scope}
\end{scope}
}\hspace{1in} 
\TIKZ[scale=.3]{
\node at (7, 6) {Type \eqref{eq:two-row-B-move}:};
\begin{scope}
\node at (7, 4) {length-decreasing:};
	\draw[fill=black!10, densely dashed]  (0,-1) to ++(6,-6) to ++(1,1) to ++(-6,6) to ++(-1,-1); 
	\node at (3.5,-3.5) {$\downarrow$};
\begin{scope}[densely dotted]
	\draw (-1,0) to ++(14,0);
	\foreach \x/\y [count=\c from 1] in {0/-1, 13/0}
	{\draw[densely dotted] (\x,\y) to (\x,.5);
		\draw[solid] (\x,.5) to coordinate (m\c) (\x, 1);}
\end{scope}
\draw (m1) to node[above] {$[i, i+2m]$}(m2);
\begin{scope}[inner sep=1pt, font=\scriptsize]
\draw[very thick, sapphire] (1,0) to node[above right]{$m$} ++(6,-6);
\draw[very thick, viridian] (7,-6) to node[below right]{$m$} ++(6,6);
\draw[very thick, rose] (0,-1) to node[below right]{$1$} ++(1,1);
\draw[very thick, black] (0,-1) to node[below left]{$p_{i-1}=-1$} ++(-1,1);
\end{scope}
\end{scope}
\begin{scope}[shift={(0,-14)}]
\node at (7, 4) {length-preserving:};
	\draw[fill=black!10, densely dashed]  (0,-1) to ++(6,-6) to ++(1,1) to ++(-6,6) to ++(-1,-1); 
	\node at (3.5,-3.5) {$\downarrow$};
\begin{scope}[densely dotted]
	\draw (-1,0) to ++(14,0) (-1,0) to ++(0,-2);
	\foreach \x/\y [count=\c from 1] in {0/-1, 13/0}
	{\draw[densely dotted] (\x,\y) to (\x,.5);
		\draw[solid] (\x,.5) to coordinate (m\c) (\x, 1);}
\end{scope}
\draw (m1) to node[above] {$[i, i+2m]$}(m2);
\begin{scope}[inner sep=1pt, font=\scriptsize]
\draw[very thick, sapphire] (1,0) to node[above right]{$m$} ++(6,-6);
\draw[very thick, viridian] (7,-6) to node[below right]{$m$} ++(6,6);
\draw[very thick, rose] (0,-1) to node[below right]{$1$} ++(1,1);
\draw[very thick, black] (0,-1) to  ++(-1,-1);
\node[below left] at (0,-2) {$p_{i-1}=1$};
\end{scope}
\end{scope}
}
\]
\caption{The effect of each type of edge $p \xrightarrow{[i, i+2m]} p'$ on a path's descent composition obtained from combining 
Theorem \ref{theorem.descent composition} and \eqref{eq: two-row descent composition table}.}\label{two-row descent compositions}
\end{figure}

Suppose the local minima of a path $p$ occur at positions $m_1,m_2,\ldots,m_h$; that is, $(p_{m_i}, p_{m_i+1}) = (-1,1)$ for $1\leqslant i \leqslant h$.
Then for $n=\lambda_1+\lambda_2$ we have
\[
	\Des(T_p)=(m_1,m_2-m_1,\ldots,m_h-m_{h-1},n-m_h).
\]
Recall the behavior of descent compositions of adjacent vertices in $\CS(\lambda)$ described in Theorem \ref{theorem.descent composition}. 
We can specialize the conditions in Theorem \ref{theorem.descent composition} to the case that $\lambda$ has two rows, and translate them 
into conditions on paths. This is done below in \eqref{eq: two-row descent composition table}, where the edge labeled $I = [i,i+2m]$ is either of 
type \eqref{eq:two-row-A-move} or \eqref{eq:two-row-B-move}.
\begin{equation}\label{eq: two-row descent composition table}
\def\arraystretch{1.2}
\text{\begin{tabular}{|c|c|c|}
\hline
$I = [i, i+2m]$ 
	& $\{i-1\} \sqcup I$ 
	&  $I \sqcup \{i+2m + 1\}$\\\hline\hline
Type \eqref{eq:two-row-A-move}
	& always rectangle-free
	&{\begin{tabular}{cl}
		rectangle &if $p_{i+2m+1}=-1$\\
		rectangle-free &  if $p_{i+2m+1}=1$
		\end{tabular}}\\\hline
Type \eqref{eq:two-row-B-move}
	&{\begin{tabular}{cl}
		rectangle &if $p_{i-1}=-1$\\
		rectangle-free &  if $p_{i-1}=1$
		\end{tabular}}
	& always rectangle-free\\\hline
\end{tabular}}
\end{equation}
The effect of each case on the associated descent compositions is then given in Figure \ref{two-row descent compositions}.

\begin{example}
\label{ex:Two-row descent compositions} 
The path $p$ below has 
\[\row(p) = \plum{45}\pumpkin{7}\sapphire{9}\viridian{123}\plum{6}\pumpkin{8}, \quad \text{ and descent composition } 
\quad \alpha= (\viridian{3},\plum{3},\pumpkin{2}, \sapphire{1}).\]
The four Dyck patterns moving out from $p$ are two type \eqref{eq:two-row-A-move} edges, one length-preserving and one length
increasing; and two type \eqref{eq:two-row-B-move} edges, one length-preserving and one length-decreasing.
\[\raisebox{-0.5\height}{\includegraphics{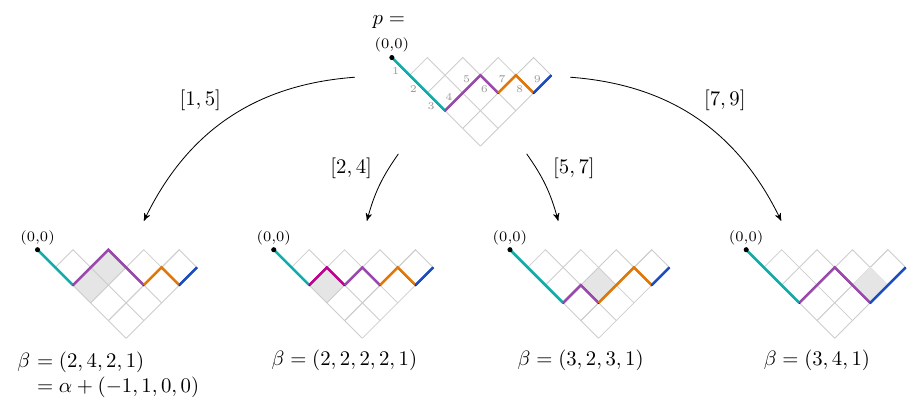}}\]
\end{example}

\subsubsection{Fans}\label{two-row: fans} 
The path model in the two-row case also allows us to be more concrete about the fans introduced in Section~\ref{section.fans}. Recall that a 
fan occurs in $\CS(\lambda)$ whenever there is a length-decreasing edge $I=[i,i+2m]$. In the language of paths, this means that there is a 
sequence length-preserving type \eqref{eq:two-row-A-move} edges moving from $I \cdot p$ to $p$, adding one box at a time. In particular, 
if $[i-2,i+2m+1]$ is not a rectangle in $I\cdot p$, then the local fan looks as follows.
\[
\TIKZ[scale=2]{
\node (p1) at (-3.5,0){$p$:
\TIKZ[scale=.2]{
\begin{scope}[black!30, thin] 
	\draw (0,-1) to ++(6,-6) to ++(1,1);
	\foreach \x in {1,2,4,5,6}{\draw (\x, -\x-1) to ++(1,1);}
	\node[rotate=-45] at (3.5, -3.5) {$\cdots$};
\end{scope}
\begin{scope}[densely dotted] 
	\draw (-1,0) to ++(14,0);
	\foreach \x/\y [count=\c from 1] in {0/-1, 13/0}
	{\draw[densely dotted] (\x,\y) to (\x,.5);
		\draw[solid] (\x,.5) to coordinate (m\c) (\x, 1);}
\end{scope}
\draw (m1) to node[above] {$I = I_1$}(m2);
\coordinate (n1) at (1,1.75); \node[above] at (n1) {$\phantom{J_1}$};
\begin{scope}[thick, inner sep=1pt, font=\scriptsize]
\draw (-1,0) 
	to ++(1,-1) 
	to ++(1,1)
	to ++(6,-6) 
	to ++(6,6);
\end{scope}
}};
\node (p2) at (-1.5,0){
\TIKZ[scale=.2]{
\begin{scope}[black!30, thin] 
	\draw (0,-1) to ++(6,-6) to ++(1,1);
	\foreach \x in {2,4,5,6}{\draw (\x, -\x-1) to ++(1,1);}
	\node[rotate=-45] at (3.5, -3.5) {$\cdots$};
\end{scope}
\begin{scope}[densely dotted] 
	\draw (-1,0) to ++(14,0);
	\foreach \x/\y [count=\c from 1] in {1/-2, 12/-1}
	{\draw[densely dotted, plum] (\x,\y) to (\x,.5);
		\draw[solid, plum] (\x,.5) to coordinate (m\c) (\x, 1);}
	\foreach \x/\y [count=\c from 1] in {-1/0, 2/-1}
	{\draw[densely dotted, viridian] (\x,\y) to (\x,1.5);
		\draw[solid, viridian] (\x,1.5) to coordinate (n\c) (\x, 2);}
\end{scope}
\draw[plum] (m1) to node[above] {$I_2$}(m2);
\draw[viridian] (n1) to node[above] {$J_1$}(n2);
\begin{scope}[thick, inner sep=1pt, font=\scriptsize]
\draw (-1,0) 
	to ++(2,-2) 
	to ++(1,1)
	to ++(5,-5) 
	to ++(6,6);
\end{scope}}};
\node (pa-2) at (1.5,0){
\TIKZ[scale=.2]{
\begin{scope}[black!30, thin] 
	\draw (0,-1) to ++(6,-6) to ++(1,1);
	\foreach \x in {4,5,6}{\draw (\x, -\x-1) to ++(1,1);}
\end{scope}
\begin{scope}[densely dotted] 
	\draw (-1,0) to ++(14,0);
	\foreach \x/\y [count=\c from 1] in {4/-5, 9/-4}
	{\draw[densely dotted, plum] (\x,\y) to (\x,.5);
		\draw[solid, plum] (\x,.5) to coordinate (m\c) (\x, 1);}
	\foreach \x/\y [count=\c from 1] in {2/-3, 5/-6}
	{\draw[densely dotted, viridian] (\x,\y) to (\x,1.5);
		\draw[solid, viridian] (\x,1.5) to coordinate (n\c) (\x, 2);}
\end{scope}
\draw[plum] (m1) to node[above, pos=.8] {$I_{m-1}$}(m2);
\draw[viridian] (n1) to node[above, pos=.2] {$J_{m-2}$}(n2);
\begin{scope}[thick, inner sep=1pt, font=\scriptsize]
\draw (-1,0) 
	to ++(5,-5) 
	to ++(1,1)
	to ++(2,-2) 
	to ++(6,6);
\end{scope}}};
\node (pa-1) at (3.5,0){
\TIKZ[scale=.2]{
\begin{scope}[black!30, thin] 
	\draw (0,-1) to ++(6,-6) to ++(1,1);
	\foreach \x in {5,6}{\draw (\x, -\x-1) to ++(1,1);}
\end{scope}
\begin{scope}[densely dotted] 
	\draw (-1,0) to ++(14,0);
	\foreach \x/\y [count=\c from 1] in {5/-6, 8/-5}
	{\draw[densely dotted, plum] (\x,\y) to (\x,.5);
		\draw[solid, plum] (\x,.5) to coordinate (m\c) (\x, 1);}
	\foreach \x/\y [count=\c from 1] in {3/-4, 6/-7}
	{\draw[densely dotted, viridian] (\x,\y) to (\x,1.5);
		\draw[solid, viridian] (\x,1.5) to coordinate (n\c) (\x, 2);}
\end{scope}
\draw[plum] (m1) to node[above, pos=.9] {$I_{m}$}(m2);
\draw[viridian] (n1) to node[above, pos=.1] {$J_{m-1}$}(n2);
\begin{scope}[thick, inner sep=1pt, font=\scriptsize]
\draw (-1,0) 
	to ++(6,-6) 
	to ++(1,1)
	to ++(1,-1) 
	to ++(6,6);
\end{scope}}};
\node[inner sep=1pt] (pa) at (0,-2){$I \cdot p$:
\TIKZ[scale=.2]{
\begin{scope}[black!30, thin] 
	\draw (5,-6) to ++(1,1) to ++(1,-1);
\end{scope}
\begin{scope}[densely dotted] 
	\draw (-1,0) to ++(14,0);
	\foreach \x/\y [count=\c from 1] in {4/-5, 7/-6}
	{\draw[densely dotted, viridian] (\x,\y) to (\x,.5);
		\draw[solid, viridian] (\x,.5) to coordinate (n\c) (\x, 1);}
\end{scope}
\draw[viridian] (n1) to node[above] {$J_{m}$}(n2);
\begin{scope}[thick, inner sep=1pt, font=\scriptsize]
\draw (-1,0) 
	to ++(7,-7) 
	to ++(7,7);
\end{scope}}};
\node[outer sep = 3pt] (d) at (0,0) {$\cdots$};
\begin{scope}[viridian]
	\draw[->] (p2) to node[below]{$J_1$} (p1);
	\draw[->] (d) to node[below]{$J_2$} (p2);
	\draw[->] (pa-2) to node[below]{$J_{m-2}$} (d);
	\draw[->] (pa-1) to node[below]{$J_{m-1}$} (pa-2);
	\draw[->, bend right] (pa) to node[below]{$J_m$} (pa-1);
\end{scope}
\begin{scope}[plum]
	\draw[->, bend right=10] (p1) to node[below]{$I_1$} (pa);
	\draw[->] (p2) to node[left]{$I_2$} (pa);
	\draw[->] (pa-2) to node[left]{$I_{m-1}$} (pa);
	\draw[->] (pa-1) to node[below]{$I_m$} (pa);
\end{scope}
}
\]
Similarly, if $I$ is a type  \eqref{eq:two-row-A-move} length-increasing edge, then there is an incoming fan consisting of a sequence of 
length-preserving~\eqref{eq:two-row-B-move} edges, \emph{removing} one box at a time. The following section will reveal concretely the dual relationships 
of these fans via Lusztig involution.

\subsection{Lusztig involution} 
Recall the technique for computing Lusztig involution on any crystal skeleton $\CS(\lambda)$ discussed in Remark~\ref{remark.evacuation}. This map 
has a nice specialization to paths, summarized below in Lemma \ref{lem:evacuation on 2 rows}.  

Given a tableau $T$ with $\row(T) = \pi = \pi_1 \pi_2 \dots \pi_n$, the word $\pi^\# = (n+1-\pi_n)\cdots(n+1-\pi_1)$ is the reading word of the tableau 
$T^\#$ obtained by rotating $T$ $180^\circ$ and taking the \defn{complement} of the fillings (i.e. replacing each filling $i$ with $n+1-i$).  

Performing jeu de taquin on $T^\#$ gives the same tableau as the RSK insertion tableau for $\pi^\#$:
\[ \evac(T) = \jdt(T^\#) = P(\pi^\#). \]
 If one follows the maximal rectangular decomposition of $T$ (in the sense of \eqref{eq:DPdecomp} and the end of 
 Example~\ref{ex:Two-row tableaux, paths, rectangular decompositions}) through the rotation-complement transformation, then jeu de taquin can be 
 performed efficiently by first resolving the rectangles individually.

\begin{example}\label{ex: evacuation on 2 row tableaux}
Returning to Example \ref{ex:Two-row tableaux, paths, rectangular decompositions}, the 
rotation and complement of $T$ is 
\[\TIKZ[scale = .4]{
\node[left] at (0,0) {$T^\# =$};
\node[left] at (0,-3) {$=$};
\node at (7.5,0) {
	\TIKZ[scale=-.4]{
	\draw[dotted] (15,1) to (15,2) to (10,2);
	\fill[black!10] (15,0) rectangle (12,1);
	\fill[plum!15] (4,2) to ++(6,0) to ++(0,-1) to ++(2,0) to ++(0,-1) to ++(-6,0) to ++(0,1) to ++(-2,0) to ++(0,1);
	\fill[black!10] (4,0) rectangle (6,1);
	\fill[viridian!15] (4,2) rectangle (0,0);
	\Tableau{{25,24,22,21,17,16,15,14,11,9,8,6,3,2,1},{23,20,19,18,13,12,10,7,5,4}}
	}};
\begin{scope}[shift={(0,-4)}]
\draw[thick, fill=black!10] (0,1) rectangle (3,2); 
\draw (1,2) to (1,1) (2,2) to (2,1);
\draw[dotted] (3,0) to (0,0) to (0,1);
	\node at (0+.5, 1.5){\small$1$}; \node at (1+.5, 1.5){\small$2$}; \node at (2+.5, 1.5){\small$3$};
\draw [|->] (1.5,-.5) to node[right]{$\jdt$} ++(0,-2);
\node at (4,1) {$\sqcup$};
\begin{scope}[shift={(5,0)}]
\draw[dotted] (2,0) to (0,0) to (0,1);
\draw[thick, plum, fill=plum!20] (0,2) to ++(6,0) to ++(0,-1) to ++(2,0) to ++(0,-1) to ++(-6,0) to ++(0,1) to ++(-2,0) to ++(0,1);
\draw (2,1) to ++(4,0);
	\foreach \x in {1,...,5}{\draw (\x,1) to (\x,2) (\x+2,1) to (\x+2,0);}
\foreach \f [count = \x from 2] in {4,5,7,10,12,13}{\node at (\x+.5, .5){\small$\f$};}
\foreach \f [count = \x from 0] in {6,8,9,11,14,15}{\node at (\x+.5, 1.5){\small$\f$};}
\draw [|->] (4,-.5) to node[right]{$\jdt$} ++(0,-2);
\end{scope}%
\node at (14,1) {$\sqcup$};
\begin{scope}[shift={(15,0)}]
\draw[dotted] (2,0) to (0,0) to (0,1);
\draw(2,1) to (2,2);
\draw[thick, fill=black!10] (0,2) rectangle (2,1); \draw (1,2) to (1,1);
	\node at (0+.5, 1.5){\small$16$}; \node at (1+.5, 1.5){\small$17$};%
\draw [|->] (1,-.5) to node[right]{$\jdt$} ++(0,-2);
\end{scope}
\node at (19,1) {$\sqcup$};
\begin{scope}[shift={(20,0)}]
\draw[thick, viridian, fill=viridian!20] 
(0,0) rectangle (4,2);
\draw (0,1) to (4,1);
	\foreach \x in {1,2,3}{\draw (\x,0) to (\x,2);}
\foreach \f [count = \x from 0] in {18,19,20,23}{\node at (\x+.5, .5){\small$\f$};}
\foreach \f [count = \x from 0] in {21, 22, 24, 25}{\node at (\x+.5, 1.5){\small$\f$};}%
\draw [|->] (2,-.5) to node[right]{$\jdt$} ++(0,-2);
\end{scope}
\end{scope}
\begin{scope}[shift={(0,-9)}]
\draw[thick, fill=black!10] (0,1) rectangle (3,0); 
\draw (1,0) to (1,1) (2,0) to (2,1);
\draw[dotted] (0,2) to (0,1);
	\node at (0+.5, .5){\small$1$}; \node at (1+.5, .5){\small$2$}; \node at (2+.5, .5){\small$3$};
\node at (4,1) {$\sqcup$};
\begin{scope}[shift={(5,0)}]
\draw[dotted] (6,0) to (8,0);
\draw[thick, plum, fill=plum!20] (0,0) rectangle (6,2);
	\foreach \x in {1,...,5}{\draw (\x,0) to (\x,2);}
	\draw (0,1) to (6,1);
\foreach \f [count = \x from 0] in {4,5,7,10,12,13}{\node at (\x+.5, .5){\small$\f$};}
\foreach \f [count = \x from 0] in {6,8,9,11,14,15}{\node at (\x+.5, 1.5){\small$\f$};}
\end{scope}%
\node at (14,1) {$\sqcup$};
\begin{scope}[shift={(15,0)}]
\draw[dotted] (0,2) to (0,1);
\draw[thick, fill=black!10] (0,0) rectangle (2,1); \draw (1,0) to (1,1);
	\node at (0+.5, 0.5){\small$16$}; \node at (1+.5, 0.5){\small$17$};
\end{scope}
\node at (19,1) {$\sqcup$};
\begin{scope}[shift={(20,0)}]
\draw[thick, viridian, fill=viridian!20] 
(0,0) rectangle (4,2);
\draw (0,1) to (4,1);
	\foreach \x in {1,2,3}{\draw (\x,0) to (\x,2);}
\foreach \f [count = \x from 0] in {18,19,20,23}{\node at (\x+.5, .5){\small$\f$};}
\foreach \f [count = \x from 0] in {21, 22, 24, 25}{\node at (\x+.5, 1.5){\small$\f$};}
\end{scope}
\end{scope}
}
\]
(compare to the end of Example \ref{ex:Two-row tableaux, paths, rectangular decompositions}).
Then 
\begin{align*}
\evac(T) = \jdt(T^\#) &= \jdt\left(
\TIKZ[scale=.4]{
		\fill[black!10] (0,1) rectangle (3,0) (9,1) rectangle (11,0);
		\fill[plum!15] (3,0) rectangle (9,2);
		\fill[viridian!15] (11,0) rectangle (15,2);
		\draw (0,1) to (15,1) (0,0) to (15,0) to (15,1);
		\foreach \i [count = \x from 0] in {1, 2, 3, 4, 5, 7, 10, 12, 13, 16, 17, 18, 19, 20, 23}{
			\draw (\x,1) to (\x,0);
			\node at (\x+.5, .5) {\small$\i$};}
		\foreach \x in {3,4,5,6,7,8,9,11,12,13,14,15}{\draw (\x,2) to (\x,1);}
		\draw (3,2) to (9,2) (11,2) to (15,2);
		\foreach \i/\x in {6/3, 8/4, 9/5, 11/6, 14/7, 15/8, 21/11, 22/12, 24/13, 25/14}{
			\node at (\x+.5, 1.5) {\small$\i$};}
		\draw[thick] 
			(3,1) to (0,1) to (0,0) to (3,0)
			(9,1) to (11,1)
			(9,0) to (11,0);
		\draw[thick,plum] (3,0) rectangle (9,2);
		\draw[thick,viridian] (11,0) rectangle (15,2);
	}
\right)\\
&=\TIKZ[scale=.4]{
\fill[black!10] (0,1) rectangle (3,0) (9,1) rectangle (11,0);
\fill[plum!15] (0,2) to ++(6,0) to ++(0,-1) to ++(3,0) to ++(0,-1) to ++(-6,0) to ++(0,1) to ++(-3,0) to ++(0,1);
\fill[viridian!15] (6,2) rectangle ++(4,-1) (11,1) rectangle ++(4,-1);
\Tableau{{1,2,3,4,5,7,10,12,13,16,17,18,19,20,23},{6,8,9,11,14,15,21,22,24,25}}
}.
\end{align*}
\end{example}

\begin{lemma}\label{lem:evacuation on 2 rows} Suppose $\lambda = (\lambda_1, \lambda_2)$.
Let $p$ be a vertex in $\CS(\lambda)$, and denote its rectangular decomposition and the associated set composition as 
\[p = (-1)^{n_0} D_1 (-1)^{n_1} D_2 (-1)^{n_2} \cdots  D_\ell (-1)^{n_\ell} \quad \text{and}\quad
\rcomp(p)=\left(\alpha^{(0)}, \beta^{(1)}, \alpha^{(1)}, \dots, \beta^{(\ell)}, \alpha^{(\ell)}\right).\]
Then the following are equivalent methods to compute $\evac(p)$. 
\begin{enumerate}
\item \emph{Rectangular decomposition.} Let $D_k^\#$ be the horizontal reflection of $D_k$ obtained from negating each term in the corresponding subsequence and reversing its order. Then
\[\evac(p) = (-1)^{n_\ell} D_\ell^\# \cdots (-1)^{n_2}D_2^\# (-1)^{n_1}D_1^\# (-1)^{n_0}.\]
\item \emph{Sequence.} For $i = 1, \dots, n$, the $i$th term in $\evac(p)$ is 
\[\evac(p)_i =  \epsilon\  p_{n+1 -i}, \quad \text{ where } \quad \epsilon = \begin{cases}
	-1 & \text{ if $n+1-i \in \beta^{(k)}$ for some $k$,}\\
	\phantom{-}1 & \text{ if $n+1-i \in \alpha^{(k)}$ for some $k$.}
	\end{cases}\]
\item \emph{Graph of the path.}  Horizontally reflect the path $p$ and place the new left endpoint at $(0,0)$. Then shine a light from the left and reflect 
each up-step that is exposed to light. The result is a graph of the path $\evac(p)$.
\end{enumerate}
\end{lemma}

\begin{proof}
Rotation of $T_p$ is the same as swapping the sign of each term of the sequence corresponding to $p$, and taking the complement is the same as 
reversing the order of the resulting sequence. Equivalently, $p^\#$ is the path obtained by horizontally reflecting $p$ and placing the left-endpoint at (0,0). 
In terms of the rectangular decomposition above, this means 
\[p^\# =      
  1^{n_\ell} D_\ell^\# \cdots 1^{n_2}D_2^\# 1^{n_1}D_1^\# 1^{n_0}.\]
Just as $T^\#$ may be a filling of a skew-shape rather than a partition, the path $p^\#$ may now exceed height 0. But since the boxes that changed rows 
from $T_p^\#$ to $\jdt(T_p^\#)$ were precisely those that started in 1-row rectangles in the maximal rectangular decomposition of $T_p$, the only steps that 
will change sign from $p^\#$ to $\jdt(p^\#)$ are those that corresponded to the eastern-exposed down-steps in $p$. Specifically, 
\[\evac(p) = \jdt(p^\#) = (-1)^{n_\ell} D_\ell^\# \cdots (-1)^{n_2}D_2^\# (-1)^{n_1}D_1^\# (-1)^{n_0}.\]
Tracking the resulting sequence will give (2). 
\end{proof}

\begin{example}\label{ex: evacuation on 2 rows}
The path corresponding to the tableau $T$ in Examples \ref{ex:Two-row tableaux, paths, rectangular decompositions} and \ref{ex: evacuation on 2 
row tableaux} reverses to
\[
p^\# = \TIKZ[scale=.25]{
\draw[densely dotted] (0,0) to +(25,0);
\begin{scope}[black!20, thin]
\clip (0,5) rectangle (25,-4.5);
\fill[maize!50](0,0) to ++(3,3) to ++(12,0) to ++(2,2) to ++(8,0) to ++(0,.5) to (0,5.5) to (0,0); 
	\draw (25,5) to ++(-15, -15);
	\foreach \d in {15, ..., 10}{\draw (25-\d, 5-\d) to +(-10,10);}
	\foreach \d [evaluate = \d as \h using \d-15, count = \c from 1] in {1, ..., 10}
		{\draw (25-\d, 5-\d) to (25-2*\d,5) to +(\h,\h);}
	\end{scope}
	\draw[thick] (0,0) 
		\foreach \x [count = \c from 1] in {2,2,2}{to coordinate (c\c) ++(1,2*\x - 3)};
	\draw[thick, plum] (3,3) node[v, plum] (I2b) {} 
		\foreach \x [count = \c from 4] in {1,1,2,1,2,2,1,2,1,1,2,2}{to coordinate (c\c) ++(1,2*\x - 3)};
	\draw[thick] (15,3) node[v, plum] (I2a) {} 
		\foreach \x [count = \c from 1] in {2,2}{to coordinate (c\c) ++(1,2*\x - 3)};
	\draw[thick, viridian] (17,5) node[v, viridian](I1b){} 
		\foreach \x [count = \c from 1] in {1,1,1,2,2,1,2,2}{to coordinate (c\c) ++(1,2*\x - 3)};
\node[v, thick, viridian, label={right, inner sep = 1pt}:{\scriptsize$\begin{matrix}(\lambda_1+\lambda_2,\quad \\\quad \lambda_1-\lambda_2)\end{matrix}$}] (I1a) at (25, 5){}; 
\draw[viridian, thick, dotted] (I1a) to (I1b);
\draw[plum, thick, dotted] (I2a) to (I2b);
\node[v, label={left}:{\scriptsize$(0,0)$}] at (0,0){};
},\phantom{p^\# = }\]
so that 
\[
\evac(p) = \jdt(p^\#) = \TIKZ[scale=.25]{
\draw[densely dotted] (0,0) to +(25,0);
\begin{scope}[black!20, thin]
\clip (0,0) rectangle (25,-8.5);
	\draw (0,0) to (15, -15);
	\foreach \d in {15, ..., 10}{\draw (\d, -\d) to +(10,10);}
	\foreach \d [evaluate = \d as \h using \d-15, count = \c from 1] in {1, ..., 10}
		{\draw (\d, -\d) to (2*\d,0) to +(-\h,\h);}
\end{scope}
	\draw[thick] (0,0) 
		\foreach \x [count = \c from 1] in {1,1,1}{to coordinate (c\c) ++(1,2*\x - 3)};
	\draw[thick, plum] (3,-3) node[v, plum] (I2b) {} 
		\foreach \x [count = \c from 4] in {1,1,2,1,2,2,1,2,1,1,2,2}{to coordinate (c\c) ++(1,2*\x - 3)};
	\draw[thick] (15,-3) node[v, plum] (I2a) {} 
		\foreach \x [count = \c from 1] in {1,1}{to coordinate (c\c) ++(1,2*\x - 3)};
	\draw[thick, viridian] (17,-5) node[v, viridian](I1b){} 
		\foreach \x [count = \c from 1] in {1,1,1,2,2,1,2,2}{to coordinate (c\c) ++(1,2*\x - 3)};
\node[v, thick, viridian, label={right, inner sep = 1pt}:{\scriptsize$\begin{matrix}(\lambda_1+\lambda_2,\quad \\\quad -\lambda_1+\lambda_2)\end{matrix}$}] (I1a) at (25, -5){}; 
\draw[viridian, thick, dotted] (I1a) to (I1b);
\draw[plum, thick, dotted] (I2a) to (I2b);
\node[v, label={left}:{\scriptsize$(0,0)$}] at (0,0){};
}.\phantom{\evac(p) = \jdt(p^\#) =}\]
\end{example}

One strength of this method of computing $\evac(p)$ is how it interfaces locally with the behavior of Dyck patterns. Recall from
Definition~\ref{definition.lusztig-on-CS} that every edge $T \xrightarrow{~I~} T'$ in $\CS(\lambda)$ has a corresponding edge 
$\LL(T') \xrightarrow{~I^\LL~} \LL(T)$, where if $I = [a,b]$, then $I^\LL = [n+1-b, n+1 - a]$. 

\begin{proposition}
Let $\lambda = (\lambda_1, \lambda_2)$ and let $p \xrightarrow{~I~} p'$ be an edge in $\CS(\lambda)$ with $I = [i,i+2m]$. As in \eqref{defn:rcomp}, write 
\[\rcomp(p) = \left(\alpha^{(0)}, \beta^{(1)}, \alpha^{(1)}, \dots, \beta^{(\ell)}, \alpha^{(\ell)}\right).\]
\begin{enumerate}[1.]
\item If $I \subseteq \beta^{(k)}$ for some $k$, then
\[p \xrightarrow{~I~} p' \text{ is of type  \eqref{eq:two-row-A-move}}\]
if and only if 
\[\LL(p') \xrightarrow{~I^\LL~} \LL(p) \text{ is of type  \eqref{eq:two-row-B-move}}.\]
\item Otherwise, $[i+1, i+2m] \subseteq \beta^{(k)}$ for some $k$ but $i \in \alpha^{(k-1)}$, in which case both 
\[p \xrightarrow{~I~} p' \quad \text{ and } \quad \LL(p') \xrightarrow{~I^\LL~} \LL(p)\]
are of type \eqref{eq:two-row-A-move}.
\end{enumerate}
\end{proposition}

\begin{proof}
This is a direct consequence of Lemmas \ref{lem:two-row Dyck patterns} and \ref{lem:evacuation on 2 rows}, and 
Lemma \ref{lem:two-row Dyck patterns on rcomp}. For edges of both type \eqref{eq:two-row-A-move} and \eqref{eq:two-row-B-move}, the subpath over 
the interval $[i+1, i+2m]$ is the trivial Dyck path, and hence must be contained in some $\beta^{(k)}$. Thus
\[\evac(p)_{n+1-j} = -p_{j} \quad \text{for each $j \in [i+1, i+2m]$},\]
i.e.\ $\evac(p)$ over  $[n+1-(i+1), n+1-(i+2m)]$ is a horizontal reflection of $p$ over $[i+1, i+2m]$.

If $i \in \beta^{(k)}$ as well, then the same equality holds for $j=i$, in which case $\LL(p)$ and $\LL(p')$ are locally just horizontal reflections of $p$ and 
$p'$, respectively: 
\[
\TIKZ[scale=3.2]{
\node (p) at (0,0){
	\TIKZ[scale=.25]{
	\begin{scope}[black!30, thin]
		\draw (6,-5) to ++(6,6); 
		\foreach \x in {1, 2, 3, 5, 6}{\draw (\x+6, -6+\x+1) to ++(1,-1);}
		\node[rotate=45] at (10.5,-1.5) {$\cdots$};
	\end{scope}
	\draw [densely dotted] 
		(-.75,1) node[left, inner sep=1pt]{\footnotesize$h_{a-1}$} to +(.75,0) 
		(0,1) to ++(13,0) to ++(0,-1);
	\draw[very thick] 
		(0,1) to ++(6,-6)
		(7,-6) to ++(6,6)
		(6,-5) to ++(1,-1);
	\draw[line width = 2pt, viridian] (0,1) to ++(1,-1);
	}};
\node (p') at (2,0){
\TIKZ[scale=.25]{
\begin{scope}[black!30, thin]
	\draw (7,-6) to ++(6,6); 
	\foreach \x in {1, 2, 3, 4, 6}{\draw (\x+5, -7+\x+1) to ++(1,-1);}
	\node[rotate=45] at (10.5,-1.5) {$\cdots$};
\end{scope}
\draw[densely dotted] 
		(-.75,1) node[left, inner sep=1pt]{\footnotesize$h_{a-1}$} to +(.75,0) 
		(0,1) to ++(13,0) to ++(0,-1);
\draw[very thick] 
	(0,1) to ++(6,-6)
	(6,-5) to ++(6,6)
	(12,1) to ++(1,-1);
\draw[line width = 2pt, viridian] (12,1) to ++(1,-1);
}};
\node (rp) at (0,-1){
	\TIKZ[scale=.25]{
	\begin{scope}[black!30, thin]
		\draw (1,0) to ++(6,-6); 
		\foreach \x in {0, 1, 2, 3, 5}{\draw (\x, -\x-1) to ++(1,1);}
		\node[rotate=-45] at (4.5,-4.5) {$\cdots$};
	\end{scope}
	\draw[densely dotted] 
		(0,-1) to (0,0) to ++(13,0);
	\draw[very thick] 
		(0,-1) to ++(6,-6)
		(7,-6) to ++(6,6)
		(6,-7) to ++(1,1);
	\draw[line width = 2pt, viridian] (12,-1) to ++(1,1);
}};
\node (rp') at (2,-1){
	\TIKZ[scale=.25]{
	\begin{scope}[black!30, thin]
		\draw (0,-1) to ++(6,-6); 
		\foreach \x in {1, 2, 3, 5, 6}{\draw (\x, -\x-1) to ++(1,1);}
		\node[rotate=-45] at (4.5,-4.5) {$\cdots$};
	\end{scope}
	\draw[densely dotted] 
		(0,-1) to (0,0) to ++(13,0);
	\draw[very thick, ] (1,0) to ++(6,-6) 
		(7,-6) to ++(6,6)
		(0,-1) to ++(1,1);
	\draw[line width = 2pt, viridian] (0,-1) to ++(1,1);
}};
\draw[xRightArrow] (p) to 
	node[above]{\scriptsize$[a,b]$}
	node[below] (A) {type \eqref{eq:two-row-A-move}}
	(p');
\draw[xRightArrow, shorten >=-5pt, shorten <=-5pt] (rp') to 
	node[above] (B) {\scriptsize$[n+1-b,n+1-a]$}
	node[below]{type \eqref{eq:two-row-B-move}}
	(rp);
\draw[line width=.6pt, , {Classical TikZ Rightarrow[length=1mm]}-{Classical TikZ Rightarrow[length=1mm]}, shorten >=5pt, shorten <=5pt] (A) to 
	node[left] {$\LL$} (B);
}
\]

On the other hand, suppose $i \notin \beta^{(k)}$. Then by Lemma \ref{lem:two-row Dyck patterns on rcomp}, $i$ is an eastern-exposed down-step  of 
$p$ ($i \in \alpha^{(k-1)}$) and $i+2m$ is an eastern-exposed down-step in $I \cdot p$, so that  
\[\evac(p)_{n+1-i} = p_i = -1 \quad \text{ and } \quad \evac(p')_{n+1-(i+2m)} = p'_{i+2m} = -1.\]
Then we see the following:
\[
\TIKZ[scale=3.2, xscale=1.2]{
\node (p) at (0,0){
	\TIKZ[scale=.25]{
	\fill[maize!50] (0,1) to ++(1,-1) to ++(12,0) to ++(0,1) to (0,1);
	\draw[densely dotted] (1,2.5) to (1,0);
	\draw[|-] (1,2.5) to node[above]{$\beta^{(k)}\strut$} ++(12.5,0) node[right]{$\cdots\dashv$};
	\draw[|-] (1,2.5) to  ++(-1,0) node[left]{$\vdash\cdots$};
	\node[above] at (-1,2.5){$\alpha^{(k-1)}\strut$};
	\begin{scope}[black!30, thin]
		\draw (6,-5) to ++(6,6); 
		\foreach \x in {1, 2, 3, 5, 6}{\draw (\x+6, -6+\x+1) to ++(1,-1);}
		\node[rotate=45] at (10.5,-1.5) {$\cdots$};
	\end{scope}
	\draw [densely dotted] 
		(-.75,1) node[left, inner sep=1pt]{\footnotesize$h_{a-1}$} to +(.75,0) 
		(0,1) to ++(13,0) to ++(0,-1);
	\draw[very thick] 
		(0,1) to ++(6,-6)
		(7,-6) to ++(6,6)
		(6,-5) to ++(1,-1);
	\draw[line width = 2pt, pumpkin] (0,1) to ++(1,-1);
	}};
\node (p') at (2,0){
\TIKZ[scale=.25]{
	\fill[maize!50] (12,1) to ++(1,-1) to ++(0,1);
	\draw[densely dotted] (12,2.5) to (12,0);
	\draw[|-] (12,2.5) to node[above]{$\beta'^{(k')}\strut$} ++(-12.5,0) node[left]{$\vdash\cdots$};
	\draw[|-] (12,2.5) to  ++(1,0) node[right]{$\cdots\dashv$};
	\node[above] at (14,2.5){$\alpha'^{(k')}\strut$};
\begin{scope}[black!30, thin]
	\draw (7,-6) to ++(6,6); 
	\foreach \x in {1, 2, 3, 4, 6}{\draw (\x+5, -7+\x+1) to ++(1,-1);}
	\node[rotate=45] at (10.5,-1.5) {$\cdots$};
\end{scope}
\draw[densely dotted] 
		(-.75,1) node[left, inner sep=1pt]{\footnotesize$h_{a-1}$} to +(.75,0) 
		(0,1) to ++(13,0) to ++(0,-1);
\draw[very thick] 
	(0,1) to ++(6,-6)
	(6,-5) to ++(6,6)
	(12,1) to ++(1,-1);
\draw[line width = 2pt, pumpkin] (12,1) to ++(1,-1);
}};
\node (rp) at (0,-1){
\TIKZ[scale=.25]{
\begin{scope}[black!30, thin]
	\fill[maize!50] (12,1) to ++(1,-1) to ++(0,1);
	\draw (7,-6) to ++(6,6); 
	\foreach \x in {1, 2, 3, 4, 6}{\draw (\x+5, -7+\x+1) to ++(1,-1);}
	\node[rotate=45] at (10.5,-1.5) {$\cdots$};
\end{scope}
\draw[densely dotted] 
		(0,1) to ++(13,0) to ++(0,-1);
\draw[very thick] 
	(0,1) to ++(6,-6)
	(6,-5) to ++(6,6)
	(12,1) to ++(1,-1);
\draw[line width = 2pt, pumpkin] (12,1) to ++(1,-1);
}};
\node (rp') at (2,-1){
\TIKZ[scale=.25]{
	\fill[maize!50] (0,1) to ++(1,-1) to ++(12,0) to ++(0,1) to (0,1);
	\begin{scope}[black!30, thin]
		\draw (6,-5) to ++(6,6); 
		\foreach \x in {1, 2, 3, 5, 6}{\draw (\x+6, -6+\x+1) to ++(1,-1);}
		\node[rotate=45] at (10.5,-1.5) {$\cdots$};
	\end{scope}
	\draw [densely dotted] 
		(-.75,1) 
		to +(.75,0) 
		(0,1) to ++(13,0) to ++(0,-1);
	\draw[very thick] 
		(0,1) to ++(6,-6)
		(7,-6) to ++(6,6)
		(6,-5) to ++(1,-1);
	\draw[line width = 2pt, pumpkin] (0,1) to ++(1,-1);
	}};
\draw[xRightArrow, shorten >=-5pt, shorten <=-5pt] (p) to 
	node[above]{\scriptsize$[i, i+2m]$}
	node[below] (A) {type \eqref{eq:two-row-A-move}}
	(p');
\draw[xRightArrow, shorten >=-0pt, shorten <=0pt] (rp') to 
	node[above] (B) {\scriptsize$[n+1-(i+2m),n+1-i]$}
	node[below]{type \eqref{eq:two-row-A-move}}
	(rp);
\draw[line width=.6pt, , {Classical TikZ Rightarrow[length=1mm]}-{Classical TikZ Rightarrow[length=1mm]}, shorten >=5pt, shorten <=5pt] (1,-.25) to 
	node[left] {$\LL$} (1,-.75);
}
\]
\end{proof}

\bibliographystyle{plain}
\bibliography{paper}{}

\end{document}